\documentclass[12pt,a4paper]{article}

\usepackage{hyperref}
\usepackage{breqn}
\usepackage{comment}

 \usepackage{graphicx}
 \usepackage{color}

\usepackage{pdfpages}			 %Allows you to include full-page pdfs.
\usepackage{wrapfig}			 %Lets you wrap text around figures.
\usepackage{upgreek}			 %Upright Greek characters.
\usepackage{dsfont}				 %Double-struck fonts.
\usepackage{simplewick}			 %For typesetting Wick contractions.
\usepackage{mathtools}		     %Can be used to fine-tune the maths presentation.	
\usepackage{framed}			     %For boxed text.
\usepackage{microtype}			 %pdfLaTeX will fix your kerning.
\usepackage{marvosym}			 %Include symbols (like the Euro symbol, etc.).
\usepackage{color}				 %Nice for scalable pdf graphics using InkScape.
\usepackage{transparent}	     %Nice for scalable pdf graphics using InkScape.
\usepackage[section]{placeins}			 %Lets you put in a \FloatBarrier to stop figures floating past this command.
\usepackage{mdframed,mdwlist}    %Use these for nice lists (less white space).
\usepackage{graphicx}            %Enhanced support for graphics.
\usepackage{float}               %Improved interface for floating objects. 
\usepackage{longtable}           %Allow tables to flow over page boundaries.
\usepackage{mathdots}            %Changed the basic LaTeX and plain TeX commands.
\usepackage{eucal}               %Font shape definitions to use the Euler script symbols in math mode.
\usepackage{array}               %Extending the array and tabular environments.
\usepackage{stmaryrd}            %The StMary’s Road symbol font.
\usepackage{amsthm}              %St Mary Road symbols for theoretical computer science. 
\usepackage{pifont}              %Access to PostScript standard Symbol and Dingbats fonts.
\usepackage{lipsum}              %Easy access to the Lorem Ipsum dummy text.
\usepackage{enumerate}           %Enumerate with redefinable labels. 
\usepackage[all]{xy}             %This is a special package for drawing diagrams.
\usepackage{amsmath}             %ATypesetting theorems (AMS style).
\usepackage{amssymb}             %Provided an extended symbol collection.
\usepackage[utf8]{inputenc}      %Allowed all displayable utf8 characters to be available as input.
\usepackage{fancyhdr}            %Extensive control of page headers and footers.
\usepackage{blindtext}           %Produced 'blind' text for testing.
\usepackage{tikz}                %To create graphic elements.
\usepackage[figuresright]{rotating}	%Allows large tables to be rotated to landscape.
\usetikzlibrary{shapes.geometric, arrows}
\usepackage{parskip} % no indent
\usepackage{upgreek}
\usepackage[english]{babel}
\usepackage{fancyhdr}
\usepackage{listings,color} 
\usepackage{mathrsfs}
\usepackage[T1]{fontenc}
\usepackage{dsfont}
\DeclareMathAlphabet{\mathcal}{OMS}{cmsy}{m}{n}

\usepackage{tabulary}
\usepackage{booktabs, mathrsfs}

\usepackage{graphicx,color}
\usepackage[font=footnotesize]{caption}
\usepackage{lscape}
\usepackage{mathtools,amssymb,array,amsthm,amsfonts}
\usepackage{bigfoot}
\usepackage[export]{adjustbox}
\usepackage{subcaption}
\usepackage{lmodern}
\usepackage{bbm}
\usepackage{scalerel,relsize}
\usepackage[table]{} %table
\usepackage{colortbl}
\usepackage{tabularx,ragged2e} % tabular tables

\usepackage[T1]{fontenc}

\usepackage{algorithm, bigfoot, amsfonts, mathtools, enumerate, multicol, graphicx} % to allow verbatim in footnote
\usepackage[noend]{algpseudocode}

\usepackage{amsthm}
\usepackage{listings}
\renewcommand\algorithmicthen{}
\algnewcommand{\IIf}[1]{\State\algorithmicif\ #1\ \algorithmicthen}
\algnewcommand{\EndIIf}{\unskip\ \algorithmicend\ \algorithmicif}

\let\oldReturn\Return
\renewcommand{\Return}{\State\oldReturn}
\newcommand{\vars}{\texttt}
\usepackage{mathptmx}
\usepackage{amsmath}
\usepackage{amssymb}
\usepackage{mathtools}

\usepackage{lineno}
\usepackage{amsthm}
\usepackage{amsmath}
\usepackage{amsfonts}
\usepackage{amssymb}
\usepackage{graphicx}
\usepackage{color}
\usepackage{epsf}
 \usepackage{ulem}
 \usepackage[T1]{fontenc}
\usepackage{breqn}

\theoremstyle{plain}

\usepackage{bm}
\newcommand{\vect}[1]{{\pmb{#1}}}

\definecolor{S1}{rgb}{0, 0.4470, 0.7410} % lines
\definecolor{S2}{rgb}{0.8500, 0.3250, 0.0980}
\definecolor{S3}{rgb}{0.9290, 0.6940, 0.1250}
\definecolor{S4}{rgb}{0.4940, 0.1840, 0.5560}
\definecolor{S5}{rgb}{0.4660, 0.6740, 0.1880}
\definecolor{S6}{rgb}{0.3010, 0.7450, 0.9330}
\definecolor{S7}{rgb}{0.6350, 0.0780, 0.1840}

\definecolor{t1}{RGB}{76.7550, 189.9750, 237.9150} % times
\definecolor{t2}{RGB}{255, 255, 0}
\definecolor{t3}{RGB}{255, 0 , 0}
\definecolor{t4}{RGB}{118.8300  171.8700   47.9400}
\definecolor{t5}{RGB}{255, 140, 0}
\definecolor{t6}{RGB}{72, 21, 170}
\definecolor{tn1}{RGB}{0, 200, 50} % n
\definecolor{tn2}{RGB}{20, 0 , 250}

\definecolor{DarkPastelGreen}{rgb}{0.01, 0.75, 0.24}
\definecolor{BV}{rgb}{0.54, 0.17, 0.89}
\definecolor{amber}{rgb}{1.0, 0.49, 0.0}
\definecolor{CeruleanBlue}{rgb}{0.01, 0.28, 1}

\definecolor{Turquoise}{rgb}{0.19, 0.84, 0.78}
\definecolor{codegreen}{rgb}{0,0.6,0}
\definecolor{codegray}{rgb}{0.5,0.5,0.5}
\definecolor{codepurple}{rgb}{0.58,0,0.82}
\definecolor{backcolour}{rgb}{0.95,0.95,0.92}
\definecolor{amethyst}{rgb}{0.6, 0.4, 0.8}
\definecolor{caribbeangreen}{rgb}{0.0, 0.8, 0.6}
\usepackage{tikz}
\usepackage{tkz-euclide}
\usetikzlibrary{calc,decorations.pathmorphing,patterns,shapes}

\usepackage{graphicx}
\usepackage{tikz}
\usepackage{tikz-3dplot}
\usetikzlibrary{fit,backgrounds}
\usepackage{lineno}

\hypersetup{%
    colorlinks=true,  	% no box around links
    linkcolor=red,        	% color of internal links (change box color with linkbordercolor)
    citecolor=black,      	% color of links to bibliography
    filecolor=black,      	% color of file links
    urlcolor=blue         	% color of external links
    }
\usepackage[T1]{fontenc}

\graphicspath{{./FinalFigs/}}
\DeclareMathOperator*{\argmin}{arg\,min}

\begin{document}
\title{A patch in time saves nine: Methods for the identification of localised dynamical behaviour and lifespans of coherent structures}
\date{\today}
\author{Chantelle Blachut\footnote{Corresponding author: chantelle.blachut@adelaide.edu.au}
\footnote{School of Computer and Mathematical Sciences, University of Adelaide, Adelaide, SA 5005, Australia}
, Cecilia Gonz\'alez-Tokman\footnote{School of Mathematics and Physics,
The University of Queensland,
St Lucia, QLD 4072, Australia} 
\\and  Gerardo Hern\'andez-Due\~nas \footnote{Institute of Mathematics, National Autonomous University of Mexico Campus Juriquilla, Blvd. Juriquilla 3001, Quer\'etaro, 76230, M\'exico} 
\footnote{
This version of the article has been accepted for publication, after peer review
but is not the Version of Record and does not reflect post-acceptance improvements, or any
corrections. The Version of Record is available online at: \url{https://doi.org/10.1007/s00332-023-09911-3}.
Chantelle Blachut and Cecilia Gonz\'alez-Tokman have been partially supported by the Australian Research Council Discovery Project scheme and the University of Queensland’s PWF
G. Hern\'andez-Due\~nas was supported, in part, by grants UNAM-DGAPA-PAPIIT IN112222 and Conacyt A1-S-17634. G. H-D would like to thank the hospitality of NorthWest Research Associates and the support of UNAM-PASPA-DGAPA during his sabbatical visit.}
}

\maketitle

% ********* Abstract: ********
\begin{abstract}
We develop a transfer operator-based method for the detection of coherent structures and their associated lifespans. Characterising the lifespan of coherent structures allows us to identify dynamically meaningful time windows, which may be associated with transient coherent structures in the localised phase space, as well as with time intervals within which these structures experience fundamental changes, such as merging or separation events. The localised transfer operator approach we pursue allows one to explore the fundamental properties of a dynamical system without full knowledge of the dynamics. The algorithms we develop prove useful not only in the simple case of a periodically driven double well potential model, but also in more complex cases generated using the rotating Boussinesq equations.
\end{abstract}


% ***************************************************

%%%%%%%%%%%%%%%%%%%%%%%%%%%%%%%%%%%%%%%%%%%%%%%%%%%%%%%%%%%%%%%%%%%%%%%%%
\FloatBarrier \section{Introduction}\label{Sec:Intro}
%%%%%%%%%%%%%%%%%%%%%%%%%%%%%%%%%%%%%%%%%%%%%%%%%%%%%%%%%%%%%%%%%%%%%%%%%
It is rare to have complete information regarding the evolution of a real world dynamical system. Whilst not knowing how the full system evolves is one obstacle to the effective numerical analysis of coherent structures, limiting one's investigation to localised regions of phase space could prove beneficial in the isolation of important dynamical phenomena. 
In this paper, we demonstrate how techniques from numerical ergodic theory allow one to analyse the local behaviour of coherent structures in complex dynamical systems. 
Of particular interest to this study is the exploration of {\textit{transient}} coherent structures, those structures characterised by a finite \textit{lifespan}. That is, the period over which a coherent structure persists. 

The ability to isolate dynamical behaviour using limited information allows one to quickly and effectively identify the existence of coherent structures and their associated lifespans. Furthermore, the identification of lifespans increases our ability to detect interesting dynamical behaviour, such as the merging or separation events that often characterise the birth or death of transient coherent structures.
We analyse this behaviour using the periodically driven double well potential of~\cite{BlachutChantelle2020Atot}, but we also consider the rotating Boussinesq equations~\cite{salmon1998lectures,smith2002generation} to investigate a variety of dynamical systems.

Among the models for geophysical flows involving complex interactions between dispersive waves and turbulence, the Boussinesq model may be used to study rotating stably stratified flows~\cite{smith2002generation}. In~\cite{Gerardo2014}, numerical simulations of a Boussinesq model were employed to study the effects of nonlinear wave-vortical interactions on the formation and evolution of coherent balanced structures, such as dipoles. Idealised models, like those generated from the Boussinesq equations, are especially useful for studying synoptic (large) scale transport in atmospheric and oceanic settings characterised by rotation and stable stratification. 

When using the Boussinesq model, the initial conditions employed in this analysis are characterised by the presence of evolving monopolar or dipolar like coherent structures. A cyclonic-anticyclonic pair of oppositely rotating vortices, whose interactions result in the propagation of this pair through space, is known as a dipole or modon~\cite{SternMelvin1975Mpop, MurakiDavidJ.2007Vdfs}. These pairs are capable of experiencing elastic collisions of various degrees as well as modon capture or fusion~\cite{MEUNIER2005431,BatchelorG.K.GeorgeKeith2000Aitf,McwilliamsJamesC.1982Ioiv}.
They are also known for their ability to characterise geophysical eddies which are capable of transporting heat, vorticity and momentum in a way that exhibits little interaction with their surrounding environment. Their ability to transport physical and dynamical properties over large distances has led to the characterisation of eddies as notable examples of coherent structures in the natural world~\cite{McwilliamsJamesC.1982Ioiv,FlierlGlennR.1983Tpso}. 

To identify coherent structures and their associated lifespans, our algorithms employ numerical approximations of a localised transfer (Perron-Frobenius) operator. The transfer operator approach was initially employed in the detection of persistent structures exhibiting a consistency in position through time; see e.g.~\cite{DellnitzMichael1999OtAo, MR2005610, FP_2009}. It was then employed to effectively identify coherent structures shifting through the configuration space~\cite{FroylandGary2007Doco,FroyPG2014,KPG_SR_SP_AV_2017,BalasuriyaSanjeeva2018GLcs,FSM_2010,FLS_2010,SantitissadeekornNaratip2010Ocsi,BlachutChantelle2020Atot}. In~\cite{FSM_2010} a non-global transfer operator was constructed by considering the action of the flow map $T_{\omega}: X_{\omega} \to X_{\sigma \omega}$ on some localised neighbourhood ${X}_{\omega} \subset X$, which could be much smaller than $X$. The sub-index $\omega\in \Omega$ denotes an initial state of the environment, which itself changes each time step, under the rule $\sigma: \Omega \to \Omega$.

Inspired by~\cite{FSM_2010}, Algorithms~\ref{alg:Seeding} through~\ref{alg:CSorNot} develop numerical methods for extracting dynamically useful information from the singular vectors of matrix compositions of non-global operators. The configuration space is partitioned into the pairwise disjoint collection of \textit{bins} {$\mathcal B$} $=\{ B_{1}, B_{2}, \ldots, B_{m} \}$ where $m=2^{\vars{depth}}$ for a given resolution $\vars{depth}$. We then seed an isolated area, or \textit{patch}, of given geometry and volume. For a given initial patching, $X_{\omega,0}$ constitutes a subset of the configuration space {$X=\cup_{i=1}^m B_i$} on which particles are seeded.  Non-global transition matrices are constructed using only local information regarding a given vector field. New bins $\{ B_{j}\}_{j \in J}$, where $J \subset \{ 1,2, \ldots, m\}$, are included only when an image of the directly preceding seeding is found within it. In the sequel, we denote $T_{\sigma^{\tilde n-1}\omega}\circ \dots \circ T_{\sigma\omega}\circ T_{\omega}(X_{\omega,0})$ by $X_{\omega,\tilde{n}}$.

We define $T_{\omega, 1}= T_\omega|_{X_{\omega,0}}$
and for
$\tilde{n} \in \{2,3 \ldots, n\}$, the evolution rule $T_{\omega, \tilde{n}}$ is given by $T_{\omega, \tilde{n}}= T_{\sigma^{\tilde n -1}\omega}|_{ 
X_{\omega,\tilde{n}-1}}$.
{The image of the} flow map {$T_{\omega,\tilde{n}}$} describes the terminal location of particles $x \in X_{\omega, \tilde{n}-1}$, {seeded in the initial patch $X_{\omega,0}$ and} initialised in the environment $\omega$, and {then} evolved {over} $\tilde{n}$ steps. The evolution process continues for {a total of} $n$ steps, where $n$ defines the number of matrices included in a particular composition. Each conditional composition is obtained as,
\begin{equation} \label{eq:TO_Ulam_illustration_Local}
P^{(n)}_{\omega} \coloneqq
P(\omega,1) P(\sigma \omega,2) \cdots
P(\sigma^{n-1} \omega,n).
\end{equation}
\noindent
Here $P(\omega,\tilde{n})$ are the conditional Ulam matrices,
\begin{equation} \label{eq:TO_Ulam_Local_P}
(P(\omega,\tilde{n}))_{i,j} = \frac{1}{Q} \sum\limits_{q=1}^{Q} \mathbbm{1}_{B_{j}}(T_{\omega,\tilde{n}} (x_{i,q}))
 \end{equation}
 where $i \in I$ and $j \in J$ depend not only on $\omega$ but also on the initial patch $X_{\omega,0}$ and $I$, like $J$, is a subset of the configuration space $X$, such that $\{ B_{i}\}_{i \in I}$ contains $X_{\omega,\tilde{n}-1}$ and  $\{ B_{j}\}_{j \in J}$ contains $X_{\omega,\tilde{n}}$. Here $Q$ is the number of test points in each bin. The matrix products defined in \eqref{eq:TO_Ulam_illustration_Local} describe the subsequent evolution of patched areas, seeded at the initial time, corresponding to the evolving environmental configurations $\omega, \sigma\omega, \dots, \sigma^{{n}-1}\omega$. 

Signals generated by the statistical properties of singular vectors and values of the matrices $\{ P^{(n)}_{\omega} \}$ identify modes associated with lifespans of coherent structures.
We identify the longest lived coherent structure, the lifespan associated with the minimal averaged \textit{equivariance mismatch}, as defined in Algorithm~\ref{alg:lifespan} line~\ref{alg:lifespan:pm}, and lifespans associated with maximum variance in the corresponding singular values, because as in~\cite{BlachutChantelle2020Atot}, 
persistent structures associated with large variations in singular values are expected to be associated with fundamental structural changes. An alternative layer identifies \textit{circular} coherent structures. The circularity of an object in $\mathbb{R}^2$ can be defined in terms of an isoperimetric quotient that compares volume contained by the structure to that of a disk with the same boundary length. The connection between structures characterised by minimal mixing or changes in boundary length, and isoperimetric analysis was first explored in~\cite{MR3404151}.  

In a recent work, Froyland and Koltai \cite{FK21} introduce an inflated dynamic Laplace operator and semi-material finite-time coherent sets (FTCSs) to investigate related questions regarding the number, lifetimes and evolution of coherent sets, and test these methods in settings different to ours. Our methods identify lifespans associated with dynamically meaningful coherent structures in a variety of systems. The first and simpler case explores the periodically driven double well potential initially described in~\cite{BlachutChantelle2020Atot}, whilst the more complex models generated by the Boussinesq equations in a rotating frame of reference are utilised to test the scope of our algorithms.

In Section \ref{Sec:algorithms}, we outline the algorithms that extract useful dynamical information from the singular vectors of matrix compositions of non-global operators. Details of the double well potential models and the rotating Boussinesq equations used to create the velocity datasets to test our methodology can be found in Section \ref{sec:models}. Section \ref{Sec:Results} is devoted to the discussion of the results. Using the double well potential model, seeds are located either in regions where a coherent structure is located at initial time, in areas where two structures merge, and in chaotic regions. In contrast, when using the Boussinesq equations the initial conditions consist of either dipole pairs, two monopoles that merge later on, or random noise that follows a Gaussian form on the initial vortical spectrum. Concluding remarks are left to Section \ref{Sec:Conclusions} and the details of the algorithms in Section \ref{Sec:algorithms} can be found in Appendix \ref{appd:algs}. 

%%%%%%%%%%%%%%%%%%%%%%%%%%%%%%%%%%%%%%%%%%%%%%%%%%%%%%%%%%%%%%%%%%%%%%%%%
\FloatBarrier \section{Algorithms}\label{Sec:algorithms}
%%%%%%%%%%%%%%%%%%%%%%%%%%%%%%%%%%%%%%%%%%%%%%%%%%%%%%%%%%%%%%%%%%%%%%%%%
%$\{ P_{t,\tilde{n}}^{(n)} \}_{ t_{i} \le t \le t_{F}-n, 1 \le \tilde{n} \le n}$ 
Our algorithms are aimed at insulating the dynamics of a particular subregion of phase space from the more complex noise that results when a large number of structures interact. The nature of these algorithms is outlined below, with further details provided in Appendix~\ref{appd:algs} and Section 4.3 of~\cite{thesis_Blachut}.
In this context, the environmental configuration $\omega$ is identified by the time $t$.

Algorithm~\ref{alg:Seeding} constructs our main tools, the collection of non-global Ulam matrices $\{ P(t,\tilde{n}) \}_{ t_{i} \le t \le t_{F}-n, 1 \le \tilde{n} \le n}$ and their respective products $\{ P_{t}^{(n)} \}_{ t_{i} \le t \le t_{F}-n}$, for $t_{i}$, $t_{F}$, $n$, $\mathcal{N}$ $\in \mathbb{Z}^{*}$, where $\mathbb{Z}^{*}=\{ 0\} \cup \mathbb{Z}^{+}$, $t_{i}$ is initial available time, $t_{F}$ is final available time and $\mathcal{N}$ is the chosen number of modes to explore.
This Algorithm partitions the configuration space and uniformly distributes $Q$ test points throughout each bin $B_{i}$ that has been added to the current, non-global collection at a given time step $t$ and period of evolution $\tilde{n}$. 
The utilisation of localised conditions, that is, an initial inclusion of bins with centres inside the patched region and the successive inclusion of those hit by evolved test points, requires one only integrate trajectories of interest. This leads to the development of the conditional localised flow maps $ \{ T_{t,\tilde{n}} \} _{ t_{i} \le t \le t_{F}-n, 1 \le \tilde{n} \le n}$ in which case $T_{t,\tilde{n}}: X_{t,\tilde{n}-1} \to X_{t,\tilde{n}}$ for some $X_{t,\tilde{n}} \subset X$.

Algorithm~\ref{alg:track_norm} allows one to track the evolution of modes via the singular value decomposition performed in
Algorithm~\ref{alg:Seeding}.  This algorithm builds on Algorithm~3 of~\cite{BlachutChantelle2020Atot} but utilises the collection of right singular vectors $\{V_{t}^{(n)} \}$. This allows for a thorough diffusion of the initial seeding when pairing vectors through time. 
The ordered collection of $\mathcal{N}$ vectors associated with the singular value paths defined using the path of modes tracked by the right singular vectors at $t$ are given by $\tilde{U}_{t}^{(n)}$ and $\tilde{V}_{t}^{(n)}$. Algorithm~\ref{alg:track_norm} also utilises a quasi-norm parameter $p$ to track structures through time.

For $u,v\in \mathbb{R}^d$, the formula ${\| u-v \|}_{p}={\left( \sum_{i=1}^{d} |u_{i} - v_{i} |^{p} \right)}^{(1/p)}$ defines a quasi-norm when $0<p<1$. Whilst in lower dimensions the Euclidean norm provides a natural choice by which to pair vectors through time, it is less clear what will be the most effective method by which to characterise the similarity of vectors in higher dimensions. The Manhattan distance (corresponding to $p=1$) has been found beneficial when sparsity is preferred~\cite{FroylandRossSakeralliou} but the utilisation of quasi-norms has led to mixed results~\cite{ICDT_DistanceHighDim,FLEXER2015281,IJCNN_FQNorms}. Smaller values of $p$ do not inherently circumvent the phenomenon of distance concentration, and the optimal choice appears to be highly application dependent and must be chosen empirically~\cite{FrancoisDamien2007TCoF}. In our case, $p$ is taken to be the largest value in $\mathcal{P}=\{ 0.1,0.2, \ldots, 1,2\}$ that returns the minimal value of ${\varsigma}_{z}$ averaged over time and modes. That is, $\displaystyle \max ( \argmin_{p \in \mathcal{P}} \overline{\varsigma}_{z} (p) )$ where ${\varsigma}_{z}$ denotes the equivariance mismatch for all lifespans identified according to either the conservative or relative threshold of Algorithm~\ref{alg:lifespan} for the $p$ dependent tracking defined by Algorithm~\ref{alg:Seeding} when all other inputs are held constant. In this case, equivariance mismatch $\varsigma$ is defined using the paired right singular vectors {{as outlined in Algorithm~\ref{alg:lifespan} Operation~\ref{alg:lifespan:pm}. 
The equivariance mismatch defines the distance between any
two vectors used to characterise the dynamics at the same point in time. This will be $0$ when two such vectors are the same (effective pairing) and $1$ when the two are orthogonal (mismatched pairing).
}}

Algorithm~\ref{alg:lifespan} discovers $\{z_{j,t}\}$, the lifespans of structures, for $j \in \{1, \ldots, \mathcal{N}\}$ and $t \in [t_{i}, t_{F}-n-2]$. This algorithm utilises the tracked singular vector pairs $\{\tilde{v}_{t,j}^{(n)} \}$ of Algorithm~\ref{alg:track_norm}. This allows one to identify when coherent structures experience \textit{birth} ($z_{\alpha}$) or \textit{death} ($z_{\omega}$). Two methods for the identification of lifespans are proposed. The first method is more conservative and looks to match neighbouring vectors if the angle between these vectors is less than $45^{\circ}$ in the Euclidean norm. That is, \texttt{threshold_c}$=\sqrt{2}\sin \frac{\pi}{8}$. The second method utilises a mode dependent threshold of $95\%$ change in the equivariance mismatch $\varsigma$ over two consecutive time steps, \texttt{threshold_p} $=0.95$. We take $95\%$ as the value for which to compare the efficacy of our algorithms on all models however, this is likely to be a model dependent parameter. Upper and lower bounds are also placed around this threshold. Our heuristic choices for \texttt{threshold_down} and \texttt{threshold_up} are equal to one quarter of the conservative threshold angle and $1$ minus this value, that is \texttt{threshold_down}$=\sqrt{2}\sin \frac{\pi}{32}$. This is done to control for sporadic results about the extremes that occur when equivariance mismatch is low (high) in absolute terms. In such cases the associated vectors are undoubtedly well (poorly) paired.

Algorithm~\ref{alg:3_lives} provides the layer that extracts lifespans likely to contain meaningful information. The lifespan of longest length, $z_{Eldest}$, provides meaningful information regarding the most persistent structures whilst the lifespan associated with the lowest average value of equivariance mismatch, $z_{MinEq}$, aims to locate structures that are well paired through time and thus more likely to be coherent. {{Those lifespans associated with the greatest variance in the associated singular values, $z_{MaxVarSV}$, are expected to characterise structures that experience dynamically meaningful changes whilst evolving in a coherent manner.}}

Algorithm~\ref{alg:CSorNot} provides an alternative layer, aimed at identifying those lifespans that contain more regular coherent structures for cases where it is not clear if the appropriate threshold or percentage has been utilised. Here we define \textit{regular coherent structures} as those structures, found utilising the associated paired singular vectors $\{ \tilde{v}_{t,j}^{(n)} \}$ and lifespans $\{z_{j,t}\}$, that contain at least one connected component with an isoperimetric ratio greater than a given threshold. This threshold is defined by the variable \texttt{iso_thresh} which is offered as input to Algorithm~\ref{alg:CSorNot}. 
According to one's discretion, and particular model requirements, further conditions, or layers, could also be introduced. Such layers could include, for example, a condition on the minimal lifespan length. Given that our goal is to examine the baseline efficacy of these algorithms, any further extensions are left for future research. 
%%%%%%%%%%%%%%%%%%%%%%%%%%%%%%%%%%%%%%%%%%%%%%%%%%%%%%%%%%%%%%%%%%%%%%%%%
\FloatBarrier \section{Models}\label{sec:models}
%%%%%%%%%%%%%%%%%%%%%%%%%%%%%%%%%%%%%%%%%%%%%%%%%%%%%%%%%%%%%%%%%%%%%%%%%
The algorithms developed in Section~\ref{Sec:algorithms} are tested on two classes of models. The periodically driven double well potential of Section~\ref{SSec:DWP_model} provides a simple scenario to test and validate the approach, whilst the more complex dynamics of the Boussinesq models described in Section~\ref{SSec:RBE_model} allow us to analyse the efficacy of our algorithms in a more realistic environment.

%%%%%%%%%%%%%%%%%%%%%%%%%%%%%%%%%%%%%%%%%%%%%%%%%%%%%%%%%%%%%%%%%%%%%%%%%
\FloatBarrier \subsection{Double well potential}\label{SSec:DWP_model}
%%%%%%%%%%%%%%%%%%%%%%%%%%%%%%%%%%%%%%%%%%%%%%%%%%%%%%%%%%%%%%%%%%%%%%%%%
The double well potential model described in~\cite{BlachutChantelle2020Atot} is a non-autonomous system of differential equations modelling simple mergers and separations of structures through time,

\begin{equation} \label{eqn:SimpleModel}
\left\{
  \begin{aligned}
    \dot{x}(t) & = y(t) \\
    \dot{y}(t) & = x(t)\left(\frac{x(t)}{2} + a(t)\right)\left(a(t) - \frac{x(t)}{2}\right),
  \end{aligned}
  \right.
\end{equation}
where
\begin{equation} \label{eqn:SimpleModel_Forcing}
a(t) =
\left\{
  \begin{aligned}
    & 1 & \text{if} \quad & 0 \le t\;(\bmod{100}) \le 10 \\
    & \cos^{2}\left((t-10)\frac{\pi}{60}\right) & \text{if} \quad & 10 \le t\;(\bmod{100}) \le 40 \\
    & 0 & \text{if} \quad & 40 \le t\;(\bmod{100}) \le 60 \\
    & \cos^{2}\left((t-30)\frac{\pi}{60}\right) & \text{if} \quad & 60 \le t\;(\bmod{100}) \le 90 \\
    & 1 & \text{if} \quad & 90 \le t\;(\bmod{100}) \le 100.
  \end{aligned}
\right.
\end{equation}

{The periodic nature of $\alpha(t)$ dictates that this vector field is the same at $t=0$ and $t=100$. For these values of $t$, the corresponding phase space is characterised by two distinctly coherent structures. These two structures merge to one by $t=50$. They then separate as $t$ increases to $100$. Figure~\ref{fig:DWP_VF} shows how the vector field evolves as $t$ shifts from $50$ to $100$. In this way, the vector field is characterised by a periodic merging and separating of the two structures.
\begin{figure}[htb!]
	\centering
	\begin{minipage}[b]{1\textwidth}
	\includegraphics[width=\columnwidth]{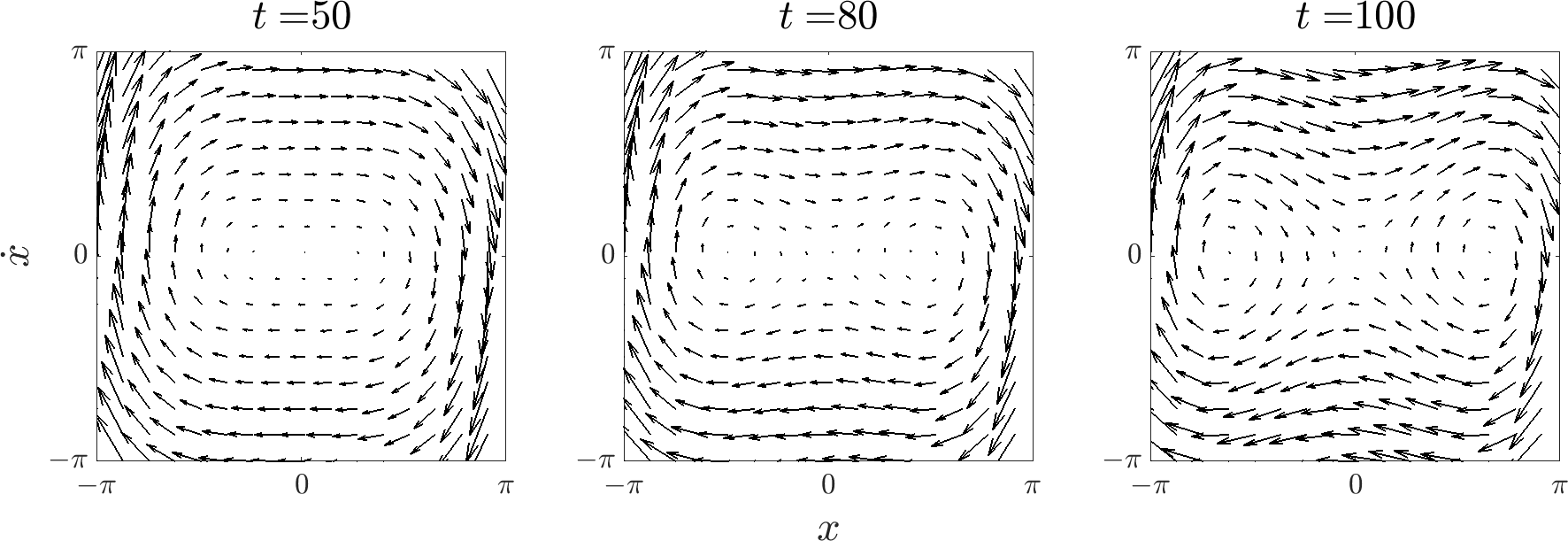}
	\end{minipage}
\caption{Vector field for the double well potential model at selected time instances.}
\label{fig:DWP_VF}
\end{figure}
}
%%%%%%%%%%%%%%%%%%%%%%%%%%%%%%%%%%%%%%%%%%%%%%%%%%%%%%%%%%%%%%%%%%%%%%%%%
\FloatBarrier \subsection{Rotating Boussinesq equations}\label{SSec:RBE_model}
%%%%%%%%%%%%%%%%%%%%%%%%%%%%%%%%%%%%%%%%%%%%%%%%%%%%%%%%%%%%%%%%%%%%%%%%%

Section~\ref{Sec:Intro} discussed the utility of the Boussinesq model to studies of transport in oceanic and atmospheric settings characterised by rotation and stable stratification. As discussed in~\cite{Remmel2010,Gerardo2014,majda2003introduction}, the Boussinesq equations for an inviscid, non-diffusive setting with stably stratified flows rotating about the vertical $\hat{\vect{z}}$-axis is given by the following equations \cite{salmon1998lectures,majda2003introduction}
\begin{equation}
\label{eqn:RBE1}
\begin{array}{rcl}
\frac{D {\vect{u}}}{Dt} + {f} \hat{\vect{z}} \times \vect{u} + N \theta \hat{\vect{z}} & = & - \nabla p_{e}, \\ \\
\frac{D\theta}{Dt} - N \vect{u}\cdot \hat{\vect{z}} & = & 0,\\ \\
\nabla \cdot \vect{u} & = & 0,
\end{array}
\end{equation} 
%, particle position is given by $\mathbf{x}=(x,y,z)$
which model vertically stratified incompressible flows. In this case the multi-dimensional vector field required by Algorithm~\ref{alg:Seeding} is given by the three dimensional vector field $\vect{u}=(u,v,w)$. As usual, $\frac{D}{Dt}={\partial{}}_{t}+\vect{u} \cdot \nabla $ is the material derivative. The Coriolis parameter is denoted by ${f}$. One notes that this is twice the frame rotation rate. The density $\rho$ is decomposed into background and fluctuating components as $\rho=\overline{\rho} + \rho'$. We assume that the background state $\bar \rho =\rho_o - \alpha z $ is linear with respect to height. Here, $\alpha > 0$ is constant for uniform stable stratification. In the derivation of the model, it is assumed that $|\rho'|,|\alpha z| \ll \rho_0$, which is valid for flows where the depth of the fluid motion is small compared to the density scale height \cite{Gerardo2014}. Given the above reference values, the variable $\theta=(\alpha\rho_0/g)^{-1/2}\rho'$ is a rescaled density fluctuation, given in units of velocity. The Brunt-V\"ais\"al\"a or buoyancy frequency is denoted by $N=(g\alpha/\rho_{0})^{1/2}$. The gravitational constant is given by $g$. The direction of gravity is $\hat{\vect z} = (0,0,1)^T$. Effective pressure is denoted $p_{e}$. This simply refers to a rescaling of pressure by $\rho_{0}$.

{
We note that although the Boussinesq equations consider three dimensional motions, we are considering parameter regimes where the dynamics are dominated by horizontal displacements. Specifically, the evolving dipoles are coherent structures that remain localized near a fixed layer in space. The vertical motions are weak and the horizontal displacements
mainly occur either around the poles or in the direction of the jet streak. A similar situation occurs for the merging monopoles in Section~\ref{SSec:MMpole}. The numerical test that takes random initial conditions considers a Coriolis term for rotation 10 times larger than the Brunt-V\"ais\"al\"a frequency for stratification. As a result, the fluid is under a rotation dominated turbulence. The 3D structure of the flow exhibits vertically coherent vortices. The algorithms used in this work utilise the velocity data obtained with the 3D model as a post-process. Due to the parameter regime considered here, the overall dynamics in these merging events are well-captured by the 2D data collected at appropriate $z$ levels. More evidence will be provided in Section~\ref{SSec:RICs}.}

%%%%%%%%%%%%%%%%%%%%%%%%%%%%%%%%%%%%%%%%%%%%%%%%%%%%%%%%%%%%%%%%%%%%%%%%%
\FloatBarrier \FloatBarrier \subsubsection{A note on computing the streamfunction}\label{SSec:RBE_model_StreamFunction}
%%%%%%%%%%%%%%%%%%%%%%%%%%%%%%%%%%%%%%%%%%%%%%%%%%%%%%%%%%%%%%%%%%%%%%%%%
Streamlines, or $\psi$ contours, of a flow, describe a family of trajectories evolving parallel to the velocity field at a given time. In the same way that linear potential vorticity can be utilised to recover the streamfunction from Charney's celebrated quasi-geostrophic derivations~\cite{MR0051104}, the projection of the Boussinesq solution into the vortical modes can be characterised by the streamfunction as,
\[
 \psi = \left[ \partial_x^2+\partial_y^2+\frac{f^2}{N^2} \partial_z^2 \right]^{-1} \left( \partial_x v - \partial_y u -\frac{f}{N} \partial_z \theta \right).
\]

One notes, that in general, such detailed information regarding natural flows is not available.
As such, these are utilised only to confirm findings and do not appear in the formal algorithms developed in Section~\ref{Sec:algorithms}. 

%%%%%%%%%%%%%%%%%%%%%%%%%%%%%%%%%%%%%%%%%%%%%%%%%%%%%%%%%%%%%%%%%%%%%%%%%
\FloatBarrier \FloatBarrier \subsubsection{Numerical scheme and parameter values}\label{SSec:RBE_model_NumericSch}
%%%%%%%%%%%%%%%%%%%%%%%%%%%%%%%%%%%%%%%%%%%%%%%%%%%%%%%%%%%%%%%%%%%%%%%%%

The numerical scheme that we employ is developed to approximate solutions of the rotating Boussinesq equations. These are constructed by employing the three dimensional periodic pseudo-spectral method discussed in~\cite{Gerardo2014} and further detailed in~\cite{smith2002generation}. There the authors employ a 2/3 dealiasing rule alongside third-order Runge-Kutta integration in time. Hyperdiffusion/hyperviscosity (higher order) damping of the form $\nu\vect{\nabla}^{16}$ is also utilised to effectively disperse energy in the smaller scales and better resolve the larger scale dynamics. In this case, the corresponding coefficient is salvaged from the energy of the highest available wavenumber shell. 

To maintain simplicity when discussing results, we continue as in Section~\ref{Sec:algorithms}, by considering normalised time units, or inertial periods. An inertial period is defined as
\begin{equation}
\label{eq:InertialPeriod}
\tau = \frac{L}{U},
\end{equation}
where $L$ and $U$ are length and velocity scales to be defined in each example. This time scale is the value employed in the fourth-order Runge-Kutta scheme utilised, in Algorithm~\ref{alg:Seeding}, to integrate the time-dependent vector fields generated by the numerical models constructed as per Sections~\ref{SSec:DDip},~\ref{SSec:MMpole} and~\ref{SSec:RICs}.

%%%%%%%%%%%%%%%%%%%%%%%%%%%%%%%%%%%%%%%%%%%%%%%%%%%%%%%%%%%%%%%%%%%%%%%%%
\FloatBarrier \FloatBarrier \subsubsection{Dancing dipoles}\label{SSec:DDip}
%%%%%%%%%%%%%%%%%%%%%%%%%%%%%%%%%%%%%%%%%%%%%%%%%%%%%%%%%%%%%%%%%%%%%%%%%

The first dataset generated to test our methodology describes the evolution of two initially balanced modon (dipolar) eddy pairs whose behaviour evolves over time. The basis for the dancing dipoles can be found in ~\cite{FlierlGR1987IEMi}. The streamfunction at $t=0$ associated with a pair of dipoles of different strengths, $\beta_1$ and $\beta_2$, is given by the equation
\begin{equation}
\begin{array}{lcl}
\left[ \frac{\partial^2}{\partial x^2}+\frac{\partial^2}{\partial y^2}+\frac{f^2}{N^2}\frac{\partial^2}{\partial z^2} \right] \tilde \psi & = & 
\beta_1~\delta(\vect x-\vect{x}_1^+)-\beta_1~\delta(\vect x-\vect{x}_1^-)+\\ \\
&& \beta_2~\delta(\vect x-\vect{x}_2^+)-\beta_2~\delta(\vect x-\vect{x}_2^-),
\end{array}
\end{equation}
with triply periodic boundary conditions in the domain $[0,2\pi]\times [0,2\pi]\times [0,2\pi]$. Each vortex strength $\pm \beta_k$ is associated to a pole $\vect{x}_k^{\pm}$ and $\delta$ is the Dirac delta function. 

Following~\cite{Gerardo2014}, the numerical implementation of the streamfunction utilises Gaussian functions to approximate the Dirac delta functions. This is done in order to smooth out the singularities near the poles. That is, the {\it modified} streamfunction satisfies 
\begin{equation}
\label{eq:psi}
\psi = \left[ \frac{\partial^2}{\partial x^2}+\frac{\partial^2}{\partial y^2}+\frac{f^2}{N^2}\frac{\partial^2}{\partial z^2} \right]^{-1} D(\vect x),
\end{equation}
where
\begin{dmath}
\label{eq:di}
D(\vect x) = \frac{1}{(2\pi \gamma)^{3/2}}\left( 
\beta_{1}~e^{-\|\vect x-\vect x_1^+\|_{2}^{2}/2\gamma}-\beta_{1}~e^{-\|\vect x-\vect x_1^{-}\|_{2}^2/2\gamma} 
+
\beta_{2}~e^{-\|\vect x-\vect x_2^+\|_{2}^2/2\gamma}-\beta_{2}~e^{-\|\vect x-\vect x_2^{-}\|_{2}^2/2\gamma}
\right)
\end{dmath}
\noindent
and $\gamma=1/128$ in these experiments. 

The initial conditions correspond to a flow in geostrophic balance, and are given by
\[
u=-\partial_y \psi, \quad v = \partial_x \psi, \quad w =0 \quad \text{and} \quad \theta = -\frac{f}{N} \partial_z \psi.
\]

The left two poles are located at $\vect x_1^\pm = (\pi/2,\pi\pm a/2,\pi\pm h/2)$ and the right two poles at $\vect x_2^\pm = (\pi,\pi\pm a/2,\pi\pm h/2)$. In this case $a=0.5$ is the separation in the meridional direction between the two poles in each dipole, whilst $h=0.5$ is the corresponding difference in elevation. The vortex centres (poles) have strengths $\beta_1=20$ and $\beta_2=10$ for the left and right dipoles respectively.

The theoretical speed at which each dipole will move in the $\hat{\vect x}$-direction, under the quasi-geostrophic dynamics, is given by
\begin{equation}
\label{eq:SpeedC}
c_k=\frac{N \beta_k a}{4\pi f}\left( a^2+\frac{N^2}{f^2}h^2 \right)^{-3/2} \quad \text{where} \quad k = 1,2.
\end{equation}
\noindent
In this particular case, the Coriolis and buoyancy frequencies are set at $f=94.1$ and $N=9.41$ with the left dipole moving at twice the speed of the right.

\begin{figure}[htb!]
	\centering
	\begin{minipage}[b]{1\textwidth}
	\includegraphics[trim=15 15 105 27,clip,width=\columnwidth]{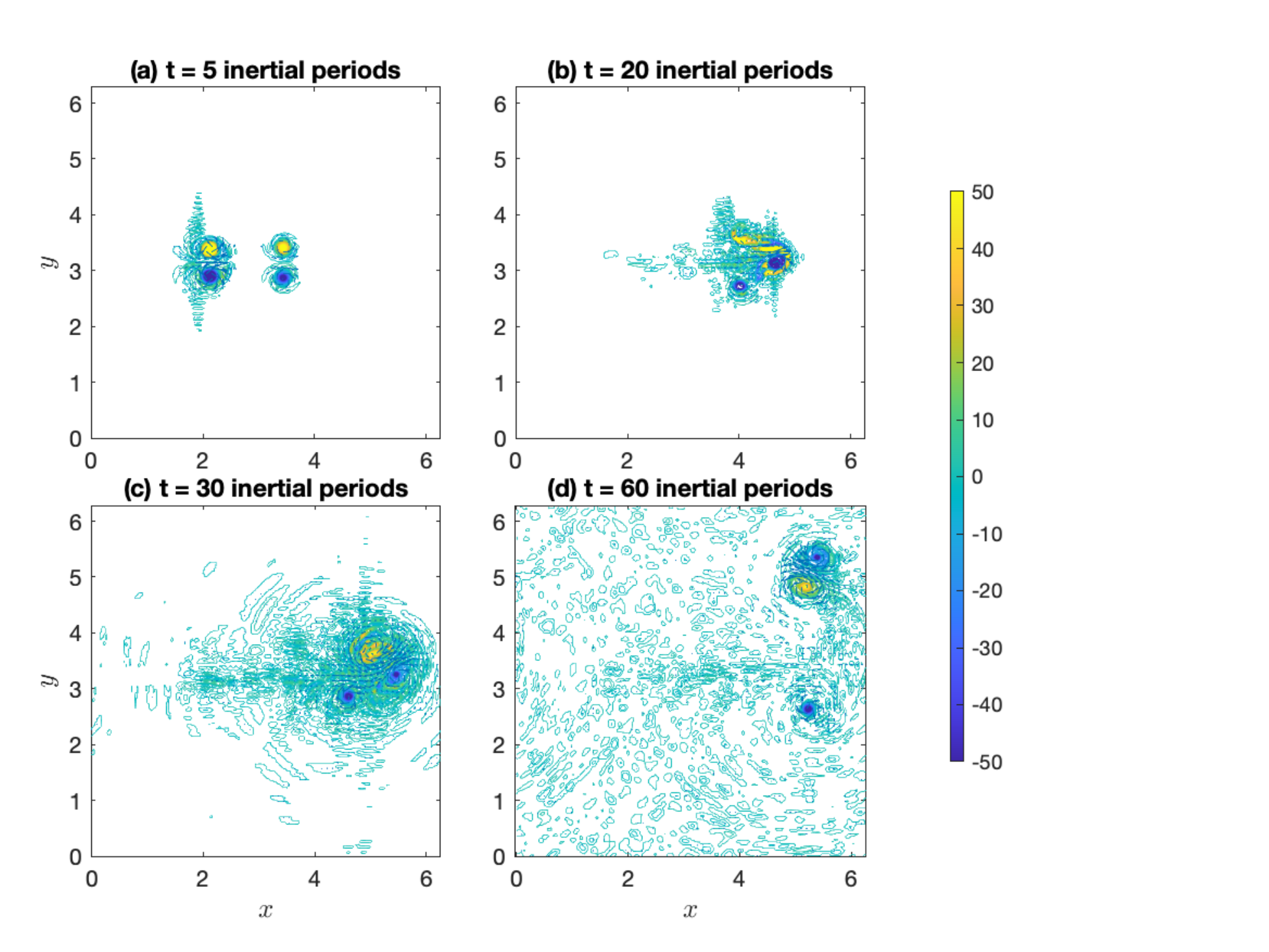}
	\end{minipage}
\caption{Horizontal contours of vertical vorticity and vector field indicating dipole's evolution at times $t=5 \tau$ (a), $t = 20 \tau$ (b), $t = 30\tau$ (c), and $t = 60\tau$ (d) for height $z=\pi$ where the arrows indicate the velocity field.}
\label{fig:TwoDipoles}
\end{figure}

Although the Boussinesq system is solved in a 3D domain, the numerical analysis to detect coherent structures is limited to two dimensional slices of the time dependent velocity field at a height of $z=\pi$ at times $t=0,1\tau,\ldots,100 \tau$. Figure~\ref{fig:TwoDipoles} shows the evolution of the two dipoles at times 5, $20$, $30$ and $60$ inertial periods~\footnote{Evolution of this vector field is shown in the supplementary file DDipolesVectorField.avi. Further data is available on request.}. The horizontal contours correspond to vertical vorticity at $z=\pi$. The velocity field is indicated by arrows. 

%%%%%%%%%%%%%%%%%%%%%%%%%%%%%%%%%%%%%%%%%%%%%%%%%%%
%%%%%%%%%%%%%%%%%%%%%%%%%%%%%%%%%%%%%%%%%%%%%%%%%%%
Here $L=1.05$ is the lengthscale, which is computed as twice the difference in the meridional positions of each pole, and $U=9.86$ is the velocity scale which is defined as the $L^\infty$ norm of the velocity field, which gives $\tau = 0.106$, according to equation \eqref{eq:InertialPeriod}. 

One notes that the above parameters correspond to a strong rotation regime since the Froude and Rossby numbers are
\[
\text{Fr} = \frac{U}{N L} = 1, \text{Ro} = \frac{U}{f L} = 0.1.
\]

This model is characterised by an overtaking collision whereby the stronger dipole overtakes a weaker one. A merger occurs between the cyclonic upper vortices whilst pseudo-merging and splitting events characterise the evolution of the anti-cyclonic lower halves. We know that the time spent in a collision state will determine the extent of structural changes that a dipole experiences for a fixed relative speed. The more destructive structural changes are associated with increasingly inelastic collisions~\cite{McwilliamsJamesC.1982Ioiv}. 

%%%%%%%%%%%%%%%%%%%%%%%%%%%%%%%%%%%%%%%%%%%%%%%%%%%%%%%%%%%%%%%%%%%%%%%%%
\FloatBarrier \FloatBarrier \subsubsection{Merging monopoles}\label{SSec:MMpole}
\begin{figure}[htb!]
	\centering
	\begin{minipage}[b]{1\textwidth}
	\includegraphics[trim=15 15 105 27,clip,width=1\columnwidth]{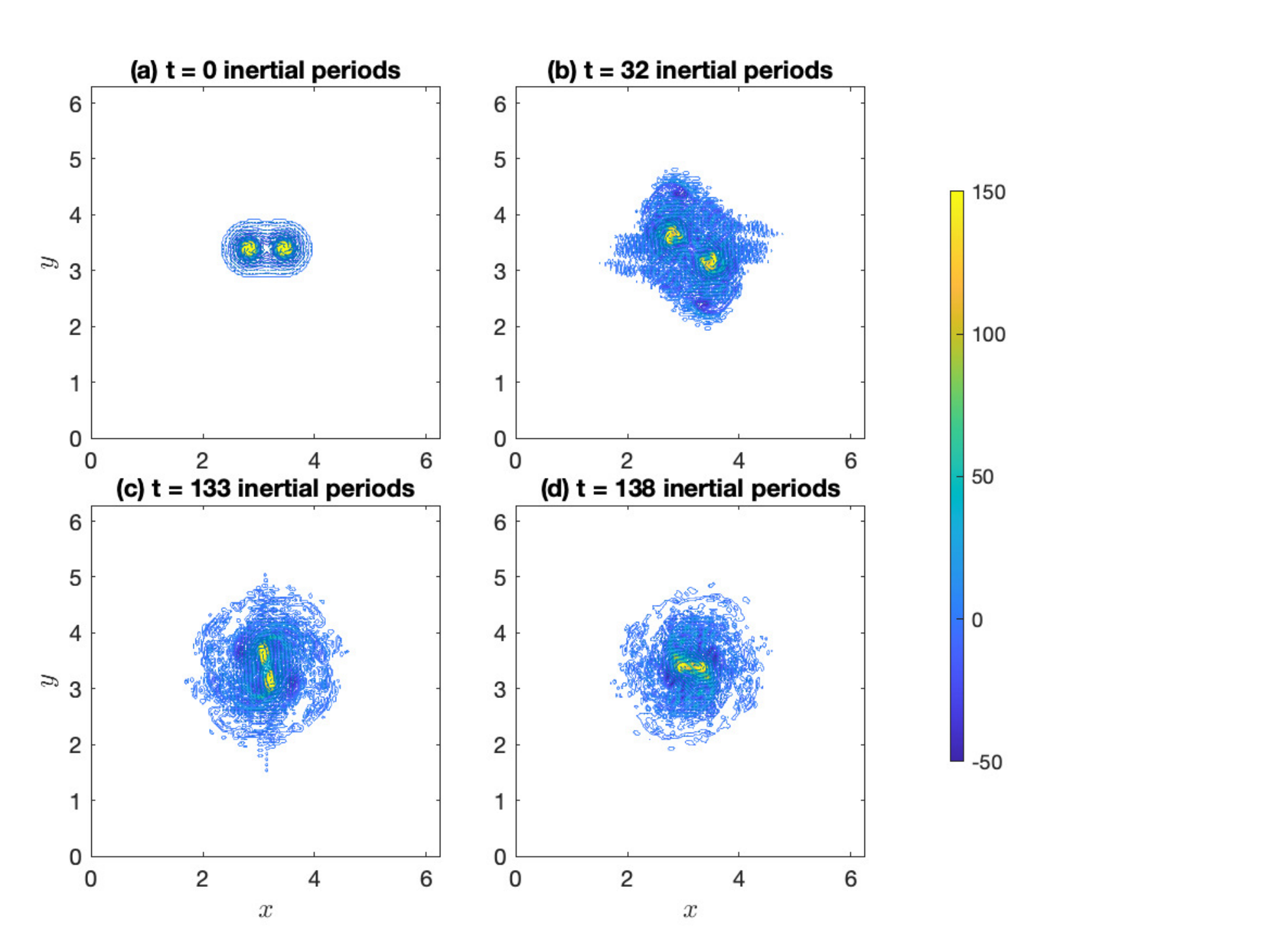}
	\end{minipage}
\caption{Time evolution of the two monopoles at $t=0\tau$ (top left), $t=32\tau$ (top right), $t=133\tau$ (bottom left), $t=138\tau$ for height $z=\pi$. Here $\tau = 0.27$ according to equation \eqref{eq:InertialPeriod}.}
\label{fig:Monopoles}
\end{figure}

%%%%%%%%%%%%%%%%%%%%%%%%%%%%%%%%%%%%%%%%%%%%%%%%%%%%%%%%%%%%%%%%%%%%%%%%%
We also consider the merging of two equal strength monopoles both of which are rotating in a counterclockwise direction. Figure \ref{fig:Monopoles} illustrates the evolution of these two poles at $0$, $32$, $133$ and $138$ inertial periods~\footnote{Supplementary file MMonopolesVectorField.avi illustrates the full time frame.}.

To construct this model one again considers the domain $[0,2\pi]\times [0,2\pi]\times [0,2\pi]$ but only introduces two poles of positive vorticity. The centre of each pole is placed near to the other. The goal of this is to encourage their interaction and eventual merging. In three dimensions these two poles are located at
\[
\vect{x}_0^- =  (0.9 \pi, \pi+0.25, \pi)  , \; \; \;  \vect{x}_0^+= (1.1 \pi, \pi+0.25, \pi)
\]
and the associated streamfunction is given by
\begin{equation}
\left[ \frac{\partial^2}{\partial x^2}+\frac{\partial^2}{\partial y^2}+\frac{f^2}{N^2}\frac{\partial^2}{\partial z^2} \right] \psi = 
\beta~\delta(\vect x-\vect{x}_0^+)+\beta~\delta(\vect x-\vect{x}_0^-).
\end{equation}
with the Dirac delta functions again approximated by the corresponding Gaussian. In this case one sets $\beta=10$ and $f=N=36.6$. We utilise the two dimensional flow for time slices at the height $z=\pi$ for times $t=0,1\tau,\ldots,100 \tau$ in our numerical analysis. In this case, we take the lengthscale as $L = 0.628$, and the velocity scale is $U=23.0$. We note that $f$ and $N$ are smaller, which increases the Froude and Rossby numbers, which corresponds to a regime farther away from strong rotation/stratification when compared to the previous case. This could serve as a sensitivity analysis where the algorithms proposed in this work are tested in different parameter regimes. The corresponding timescale is quite small. Taking into account the time it takes for the monopoles to merge, we instead consider a larger timescale $\tau = 0.27$.

%%%%%%%%%%%%%%%%%%%%%%%%%%%%%%%%%%%%%%%%%%%%%%%%%%%%%%%%%%%%%%%%%%%%%%%%%
\FloatBarrier \FloatBarrier \subsubsection{Random initial conditions}\label{SSec:RICs}
\begin{figure}[htb!]
	\centering
	\begin{minipage}[b]{1\textwidth}
	\includegraphics[trim=15 15 118 27,clip,width=1\columnwidth]{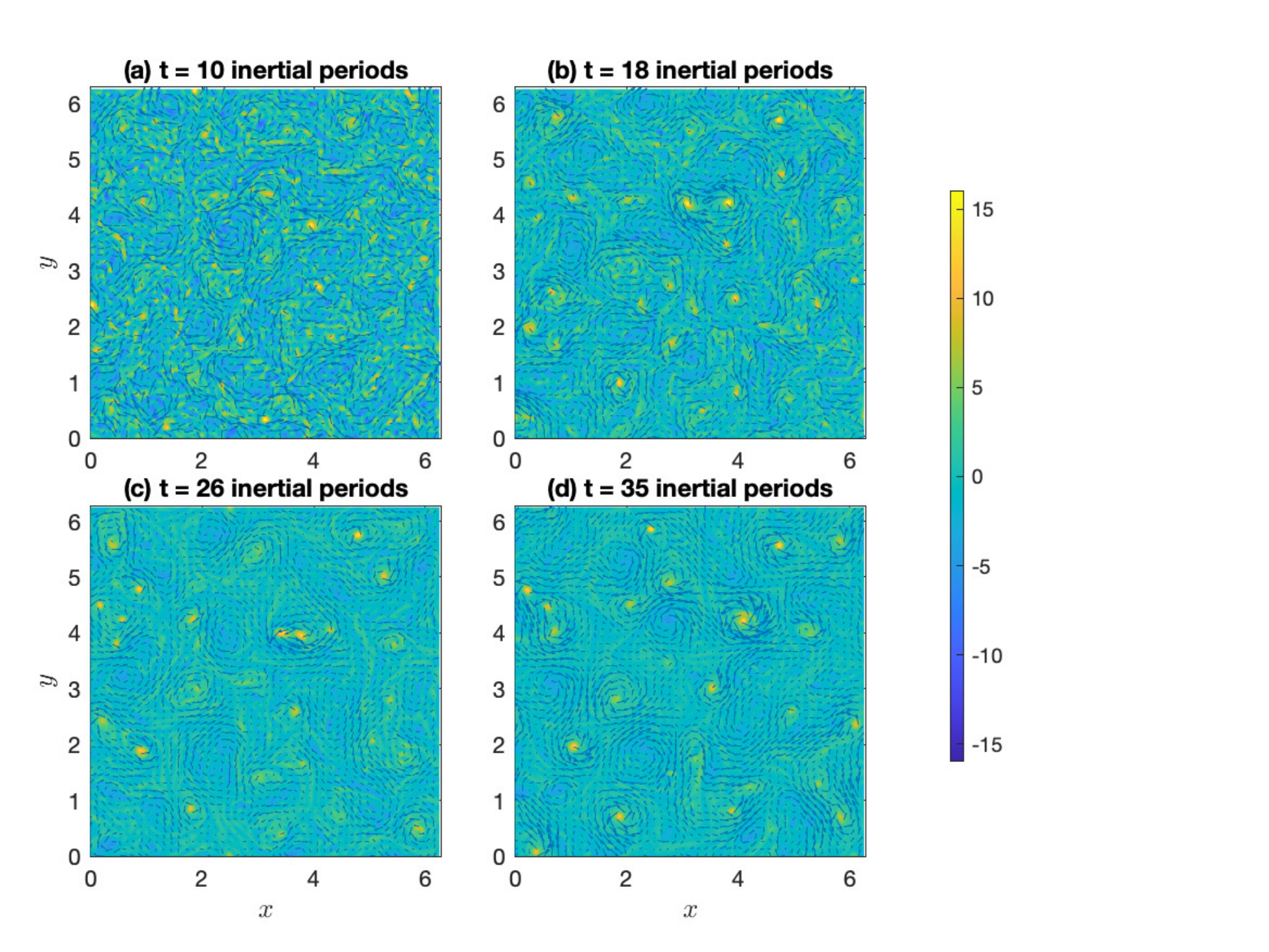}
	\end{minipage}
\caption{Contours of vertical vorticity and vector field indicating structure movement at $z=\pi$ for times $10\tau$, $18\tau$, $26\tau$, $35\tau$.}
\label{fig:complex}
\end{figure}

%%%%%%%%%%%%%%%%%%%%%%%%%%%%%%%%%%%%%%%%%%%%%%%%%%%%%%%%%%%%%%%%%%%%%%%%%
In contrast to the well organised initial state of the dancing dipoles in Section~\ref{SSec:DDip} or the merging monopoles of Section~\ref{SSec:MMpole}, this numerical test (rotation dominated turbulence) illustrates the formation of coherent structures from a set of random initial conditions. In doing this, one draws upon the tendency of geostrophic turbulence to form coherent structures~\cite{Cushman-RoisinBenoit2011QD}. Indeed, the time evolution is characterised by a transfer of energy from small to large scales with random initial conditions and initial energy in the vortical modes. This numerical test was initially considered in~\cite[Section 5.2]{Gerardo2014}. There the authors use the following Gaussian form of the initial vortical spectrum, given as a function of the wavenumber $k$,
\begin{equation} \label{eqn:RandICs_Gaussian}
F(k)=\epsilon_{f} \frac{\exp \left( -0.5(k-k_{f})^{2} / \gamma_{f} \right)}{\sqrt{2 \pi \gamma_{f}}}
\end{equation}
for $\epsilon_{f}=0.05$, $k_{f}=15$ and $\gamma_{f}=100$. Figure~\ref{fig:complex} plots the corresponding contours of vertical vorticity and velocity field at times $10\tau$, $18\tau$, $26\tau$ and $35\tau$ for height $z=\pi$~\footnote{Evolution of this vector field is shown in the supplementary file RandICsVectorField.avi.}. Our numerical analysis covers times $t=0,1\tau,\ldots,100 \tau$
with a time scale $\tau=0.72$ with $L = 0.42$ and $U=0.59$.

\begin{figure}[h!]
\begin{center}
{\includegraphics[width=0.9\textwidth]{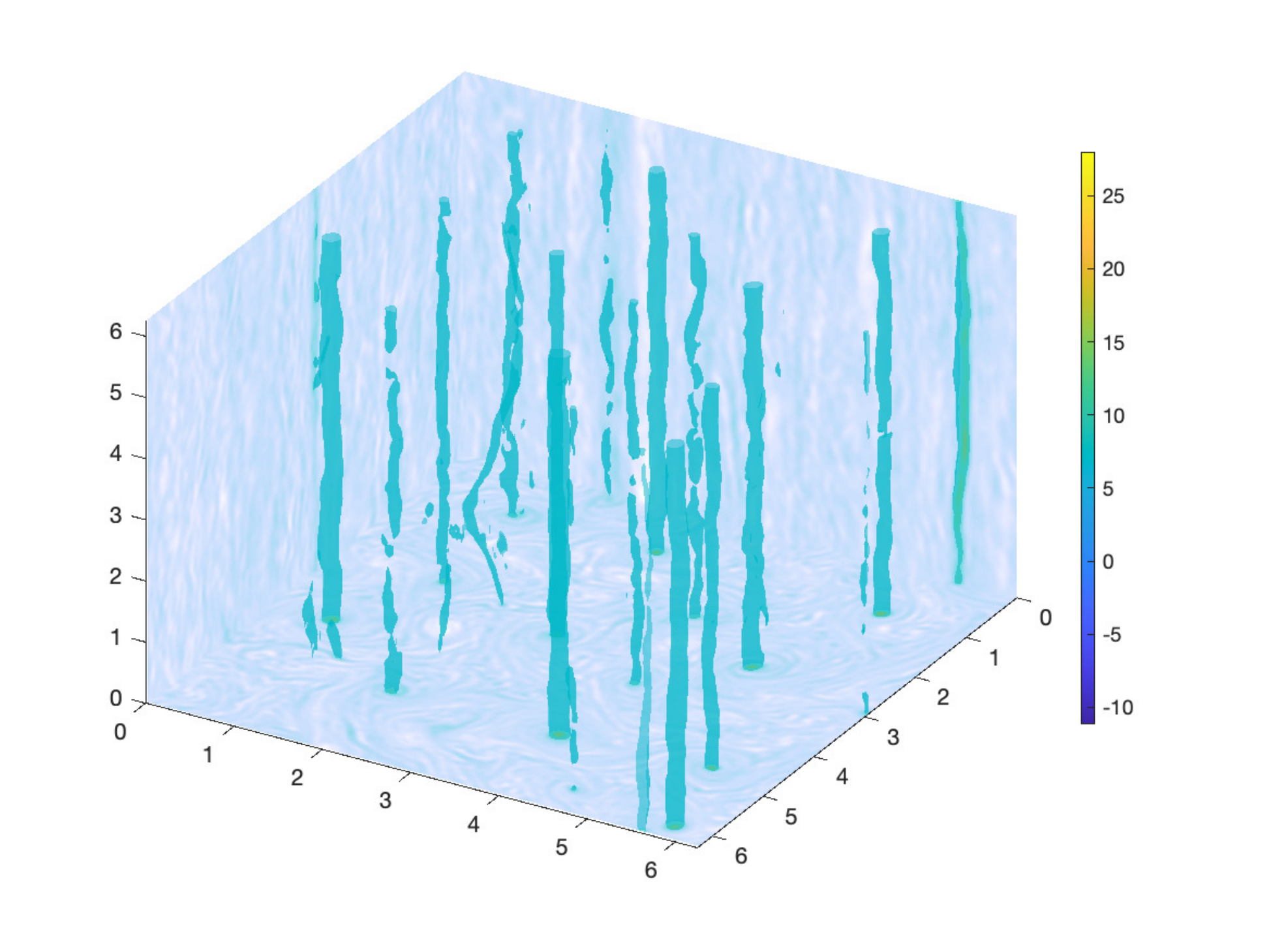}}
\end{center}
\caption{\label{fig:VertCoherence} Iso surfaces with constant vertical vorticity value $7$. The vorticity contours are displayed at the walls. }
\end{figure}

{
The analysis performed here is applied to two-dimensional data generated by a 3D model. Although the selected data is restricted to a fixed depth, the data has 3D dynamical effects. However, we are focusing on a parameter regime where rotation is strong with $f/N=10$, which is an appropriate scenario for the formation of vortices. Under these circumstances, one usually observes vertical coherence. Figure \ref{fig:VertCoherence} shows the iso-surface with constant vertical vorticity value 7, which is one quarter of the maximum value. One can observe vertically coherent vortices. Although 3D analysis might be more complete, our analysis shows that the merging of the vortices can be identified by restricting the study to horizontal 2D data at constant depth.} 

%%%%%%%%%%%%%%%%%%%%%%%%%%%%%%%%%%%%%%%%%%%%%%%%%%%%%%%%%%%%%%%%%%%%%%%%%
\FloatBarrier \section{Results and discussion}\label{Sec:Results}
%%%%%%%%%%%%%%%%%%%%%%%%%%%%%%%%%%%%%%%%%%%%%%%%%%%%%%%%%%%%%%%%%%%%%%%%%
Results from implementing the algorithms described in Section~\ref{Sec:algorithms} for the case of the periodically driven double well potential are presented in Section~\ref{SSec:DWP_model_results}, whilst those for the Boussinesq models are presented in Section~\ref{SSec:BEeModels_Results}. To focus on the dominant structures, we simplify the analysis by concentrating on the leading $4$ modes. 

In the results illustrated throughout this section, the equivariance mismatch under a given threshold is denoted as $\varsigma_{z}$. The beginning of a given lifespan is expected to be associated with the birth of a coherent structure or the entrance of a coherent structure into the patched region. We denote these time instances as $z_{\alpha}$. This is indicated by a green dot in the electronic version of the lifespan plots. In a similar manner the end, or death, of a lifespan is denoted by $z_{\omega}$. This is indicated by a red dot in the electronic version of the lifespan plots.  
Our analyses are limited to time windows comprised of $n=10$ matrices. As such, the superscript $(n)$ is henceforth omitted.

Whilst in the less complex models it is sufficient to simply identify lifespans of interest according to a hard threshold, when the dynamics are increasingly complex and the number of structures is greater, such as in those cases offered by the Boussinesq models, an additional algorithmic layer and a variable threshold are often required to identify lifespans of interest.
%%%%%%%%%%%%%%%%%%%%%%%%%%%%%%%%%%%%%%%%%%%%%%%%%%%%%%%%%%%%%%%%%%%%%%%%%
\FloatBarrier \subsection{Double well potential}\label{SSec:DWP_model_results}
%%%%%%%%%%%%%%%%%%%%%%%%%%%%%%%%%%%%%%%%%%%%%%%%%%%%%%%%%%%%%%%%%%%%%%%%%
As in~\cite{BlachutChantelle2020Atot}, we  partition $X$ into a grid of $2^{12}$ bins of equal volume. However, in contrast to the conditional inclusion~\footnote{In that case, conditional evolution included additional bins only when all image points from a single bin did not land in existing bins.} of that work, we now include all bins hit by the images of seeds at subsequent time steps. Our non-global Ulam matrices are built by seeding $Q=100$ uniformly distributed test points in certain bins, as dictated in Algorithm~\ref{alg:Seeding}. We seek to characterise time windows where $0 \le t \le t_{F}-n-1$, $t_{F}=500$ for discrete time flow maps approximated using Runge-Kutta numerical integration of a time dependent vector field flowed for $\tau=1$ using $20$ steps.

Sections~\ref{ssec:case_a} through~\ref{ssec:case_c} analyse the periodically forced double well potential by generating patches in three different regions of phase space. In the first example, Section~\ref{ssec:case_a}, we seed an area known to contain a coherent structure at initial time. In the second case, Section~\ref{ssec:case_b}, we patch a region where two structures are known to merge. In the third case, Section~\ref{ssec:case_c}, we seed a region that is known to be chaotic and is not visited by the centre of either well.
\vspace{-0.25cm}
%%%%%%%%%%%%%%%%%%%%%%%%%%%%%%%%%%%%%%%%%%%%%%%%%%%%%%%%%%%%%%%%%%%%%%%%%%%%%
\FloatBarrier \subsubsection{Seeding an area known to contain a coherent structure at initial time}\label{ssec:case_a} %Case a) (-2,0) a=b=1, n=10, depth=12
%%%%%%%%%%%%%%%%%%%%%%%%%%%%%%%%%%%%%%%%%%%%%%%%%%%%%%%%%%%%%%%%%%%%%%%%%%%%%
Let us first patch a circle of radius $1$ centred at $(-2,0)$. This patch corresponds to an area covered by the left well of the periodically driven double well potential at times $t=0 \pmod{100}$. One thus expects to detect the presence of coherent structures in the patch around such times. Our method begins with the implementation of Algorithms~\ref{alg:Seeding} and~\ref{alg:track_norm}, utilising the conservative threshold and $p=0.1$ as the most appropriate quasi-norm.

\begin{figure}[htb!] 
  \includegraphics[width=\textwidth]{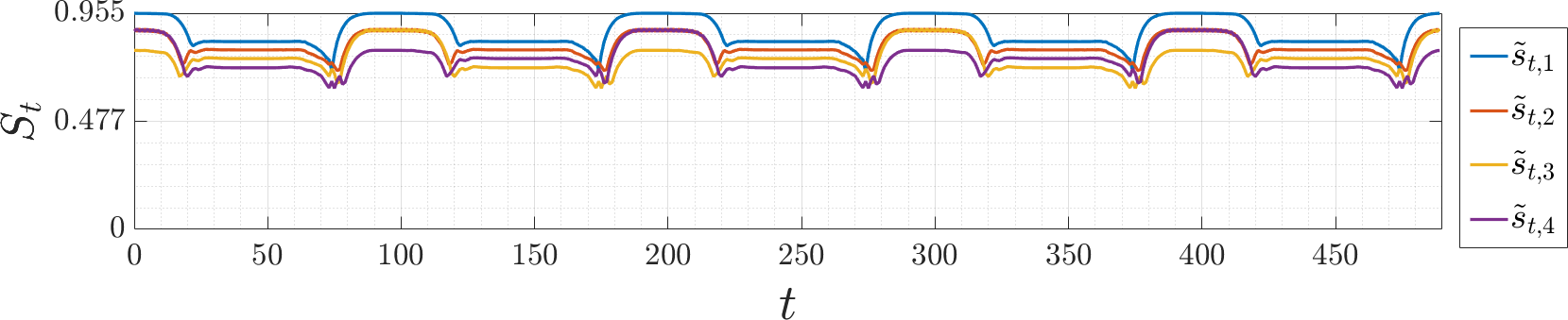}
\caption{Tracked paths for the rolling windows of singular values, where the initial seeding is concentrated on the left well for the double well potential case when $n=10$ and $p=0.1$.}
  \label{fig:case_DWP_1a_tracked}
\end{figure}

Figure~\ref{fig:case_DWP_1a_tracked} plots the paths of rolling windows of singular values that have been paired through time. This model is not subject to periodic boundary conditions, therefore a small amount of mass is expected to leave and enter the system over time. For this reason the singular values are close to $1$ but more generally they are slightly less than $1$. The leading singular value is closest to $1$ over time periods when we expect a coherent structure, that is the left-most well, to inhabit the patch.

For $t$ around times $0 \pmod{100}$, Figure~\ref{fig:case_DWP_1a_tracked} shows that the singular value paths are more clearly separated than those paths associated with time windows initialised around $t=50 \pmod{100}$. There are obvious changes around $t=25$ and $75$ when the singular values associated with each rolling window decrease and increase. These initial times coincide with merging and separation events (respectively).

\begin{figure}[H] 
\includegraphics[width=\textwidth]{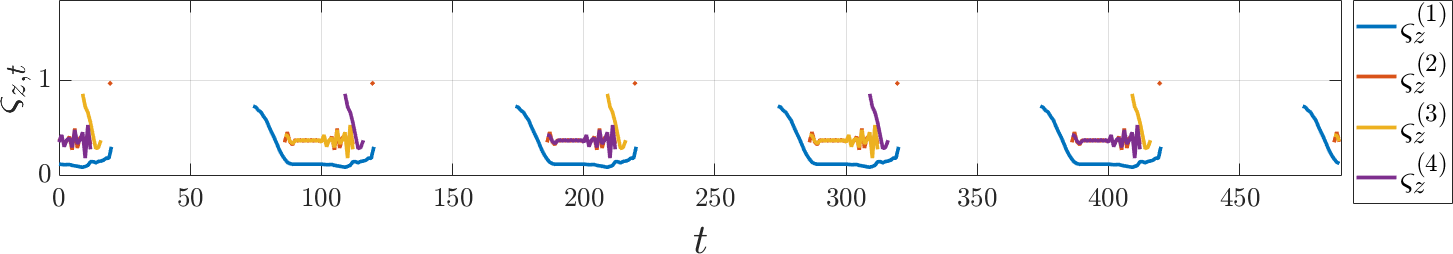}
\caption{Equivariance mismatch lifespans under a conservative threshold, as described in Algorithm~\ref{alg:lifespan}, for the case of seeding the left well in the double well potential when $n=10$ and vector pairing utilises $p=0.1$.}
\label{fig:case_DWP_1a_tracked_lifespan_eqMM_threshold}
\end{figure}

As described in Algorithm~\ref{alg:lifespan}, we track the lifespan of structures by considering the similarity of two vectors paired at neighbouring times under a conservative threshold. This is shown in Figure~\ref{fig:case_DWP_1a_tracked_lifespan_eqMM_threshold}. The beginning of a given lifespan $z_{\alpha}$, is expected to be associated with the birth of a coherent structure or the entrance of one into the patched region. In a similar manner, the end or death of a lifespan $z_{\omega}$ is expected to be associated with the loss of a coherent structure.

\begin{figure}[h!]
\setlength{\tabcolsep}{-0.5pt}
\begin{tabularx}{\columnwidth}{|p{0.65cm}| *3{>{\Centering}X}|}\hline
 \multicolumn{4}{|c|}{
\begin{minipage}[t][][b]{\columnwidth}
\centering
\vspace{-0.15cm}
\includegraphics[width=0.95\columnwidth]{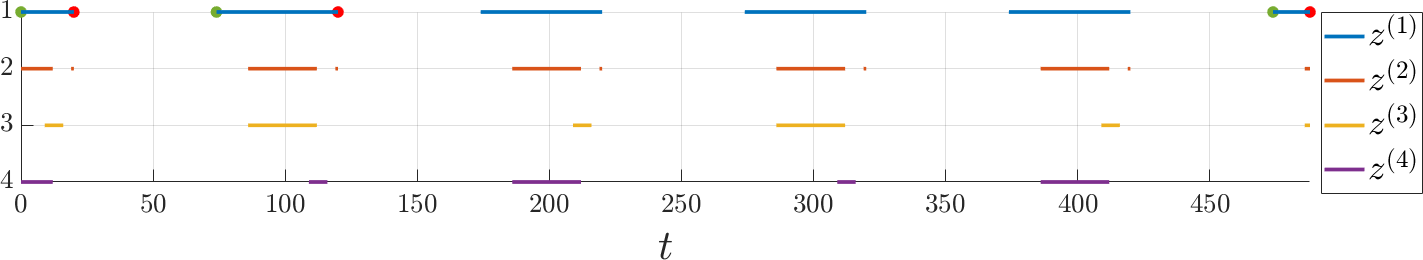}
\vspace{0.1cm}
\label{fig:case_DWP_1a_tracked_lifespan_threshold_type2_DWP_N10}
\end{minipage} }\\ 
\hline  
\multicolumn{4}{c}{\vspace{-0.51cm}} \\
\cline{2-4}
\multicolumn{1}{c|}{} &
\cellcolor{S1!60} &
\cellcolor{S1!60} &
\cellcolor{S1!60} \\[-1.em] 
\multicolumn{1}{c|}{} &
\cellcolor{S1!60} $z_{MinEq}$ for $[0,20]$&
\cellcolor{S1!60} $z_{MaxVarSV}$ for $[474,488]$& %%???????????? why this time
\cellcolor{S1!60} $z_{Eldest}$ for $[74,120]$\\
\hline
%\vspace{-0.1cm}
\begin{minipage}[t][][b]{0.03\textwidth}
\centering
\cellcolor{t4!60}
\vspace{1cm}
\footnotesize{$\; z_{\alpha}$}
\end{minipage}
&
\begin{minipage}[t][][b]{0.3\textwidth}
	\centering
	\vspace{-0.25cm}
	\includegraphics[height=2.75cm,center]{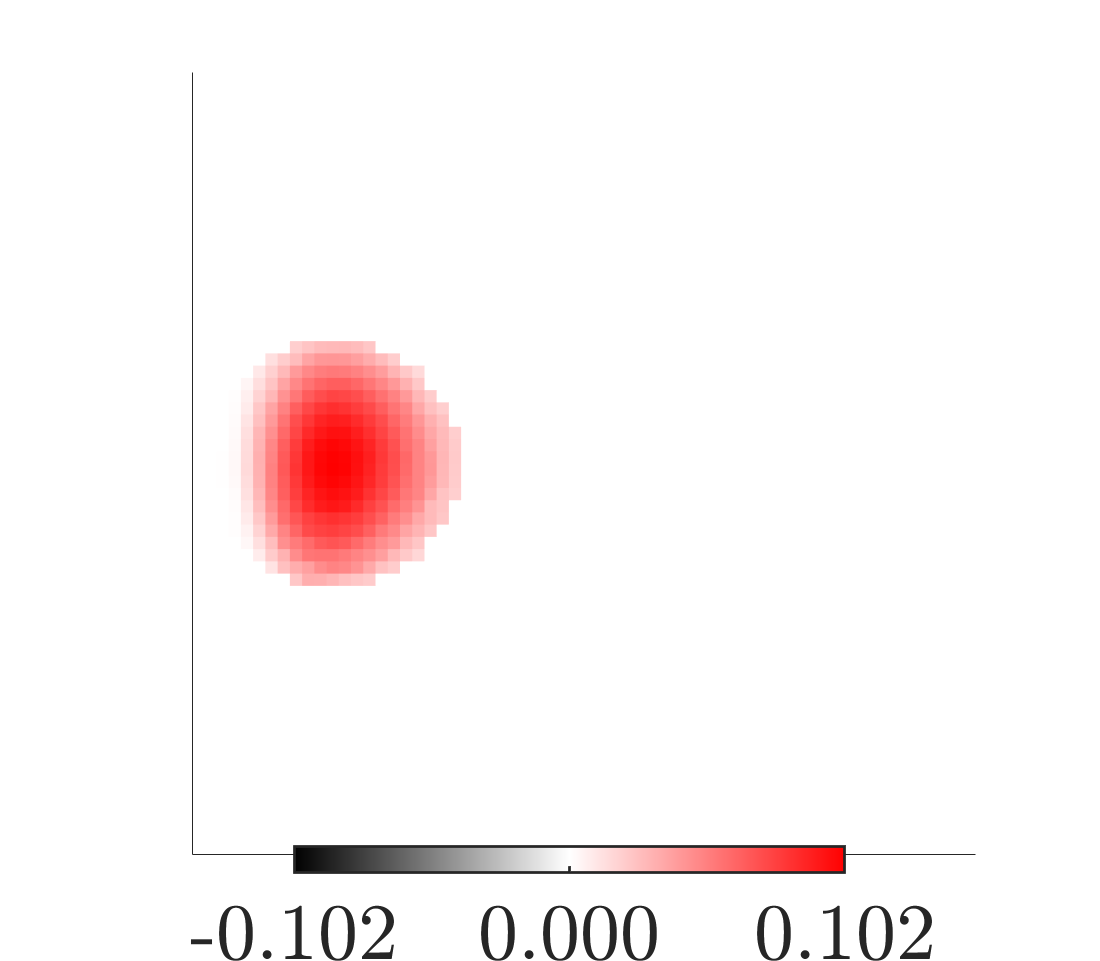} 
	\vspace{-0.25cm}
\end{minipage}
&
\begin{minipage}[t][][b]{0.3\textwidth}
	\centering
	\vspace{-0.25cm}
	\includegraphics[height=2.75cm,center]{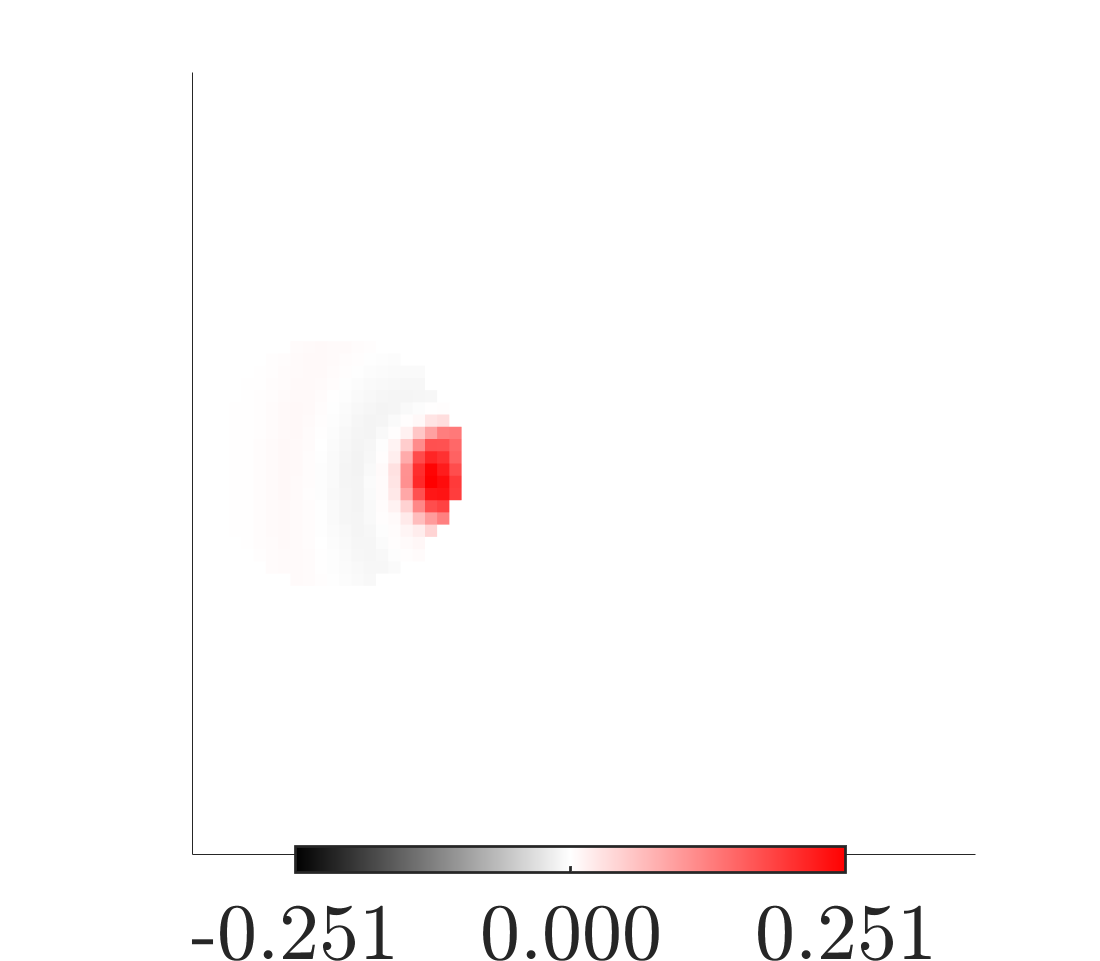} 
	\vspace{-0.25cm}
\end{minipage}
&
\begin{minipage}[t][][b]{0.3\textwidth}
	\centering
	\vspace{-0.25cm}
	\includegraphics[height=2.75cm,center]{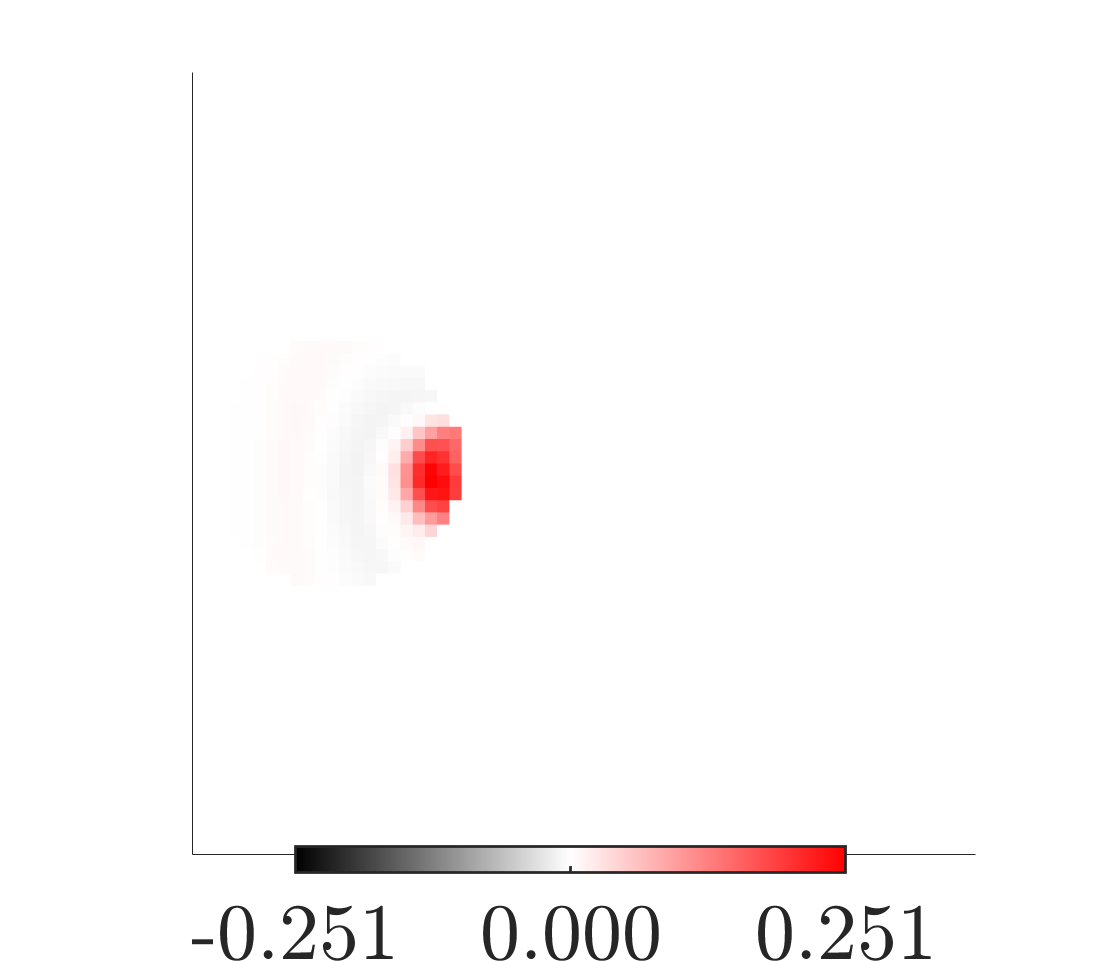} 
	\vspace{-0.25cm}
\end{minipage}
\\  \hline
\vspace{0cm}
\begin{minipage}[t][][b]{0.03\textwidth}
\centering
\cellcolor{t3!60}
\vspace{0.65cm}
\footnotesize{$\; z_{\omega}$}
\end{minipage}
&
% TRACKED mode 2 at various t_i
\begin{minipage}[t][][b]{0.3\textwidth}
	\centering
	\vspace{-0.25cm}
	\includegraphics[height=2.75cm,center]{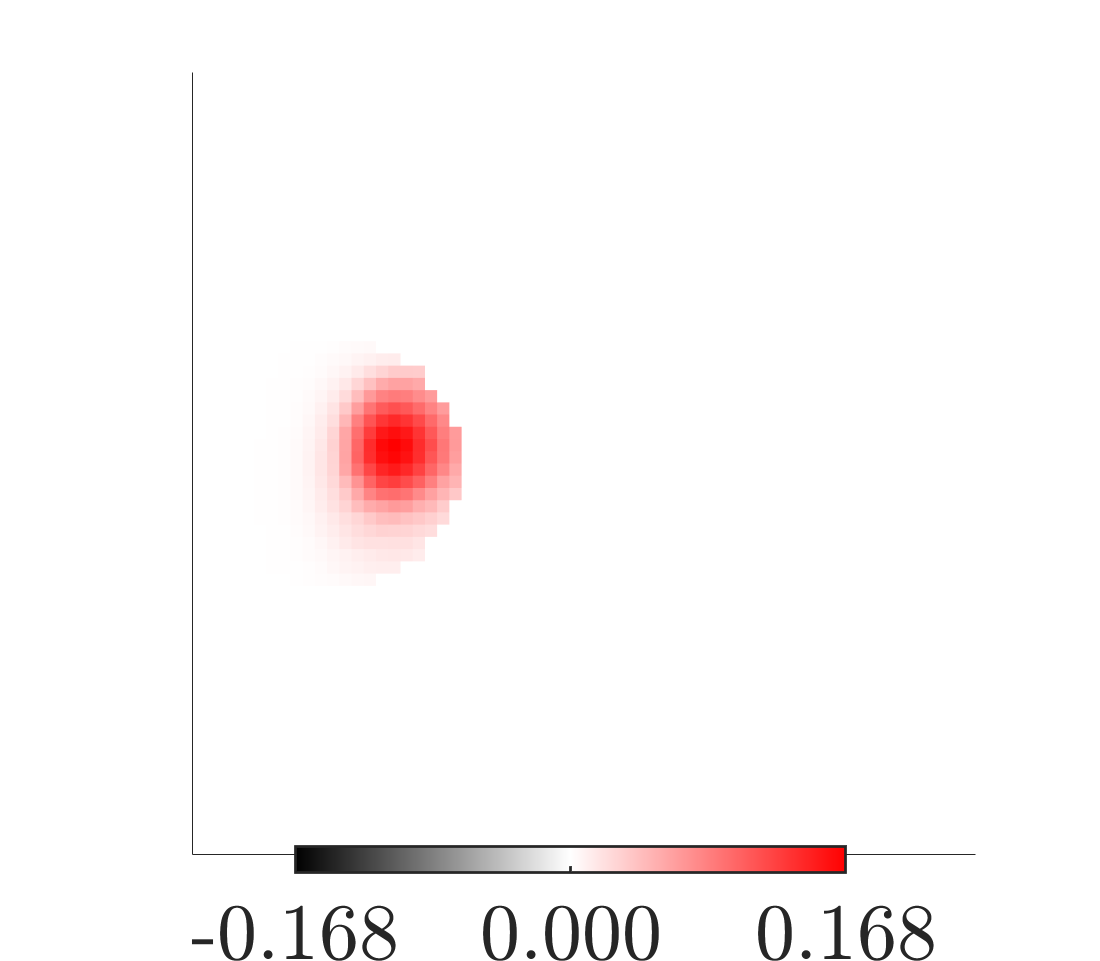} 
	\vspace{-0.25cm}
\end{minipage}
&
\begin{minipage}[t][][b]{0.3\textwidth}
	\centering
	\vspace{-0.25cm}
	\includegraphics[height=2.75cm,center]{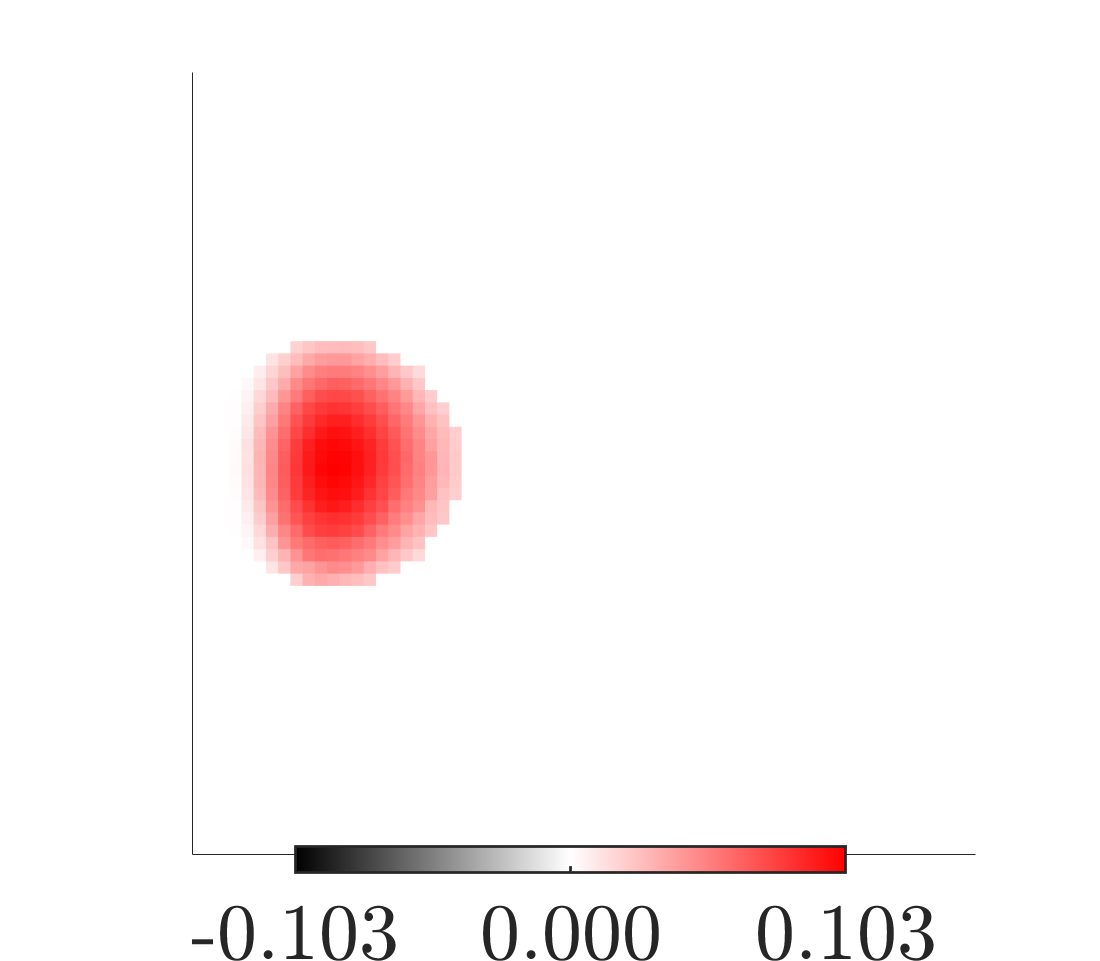} 
	\vspace{-0.25cm}
\end{minipage}
&
\begin{minipage}[t][][b]{0.3\textwidth}
	\centering
	\vspace{-0.25cm}
	\includegraphics[height=2.75cm,center]{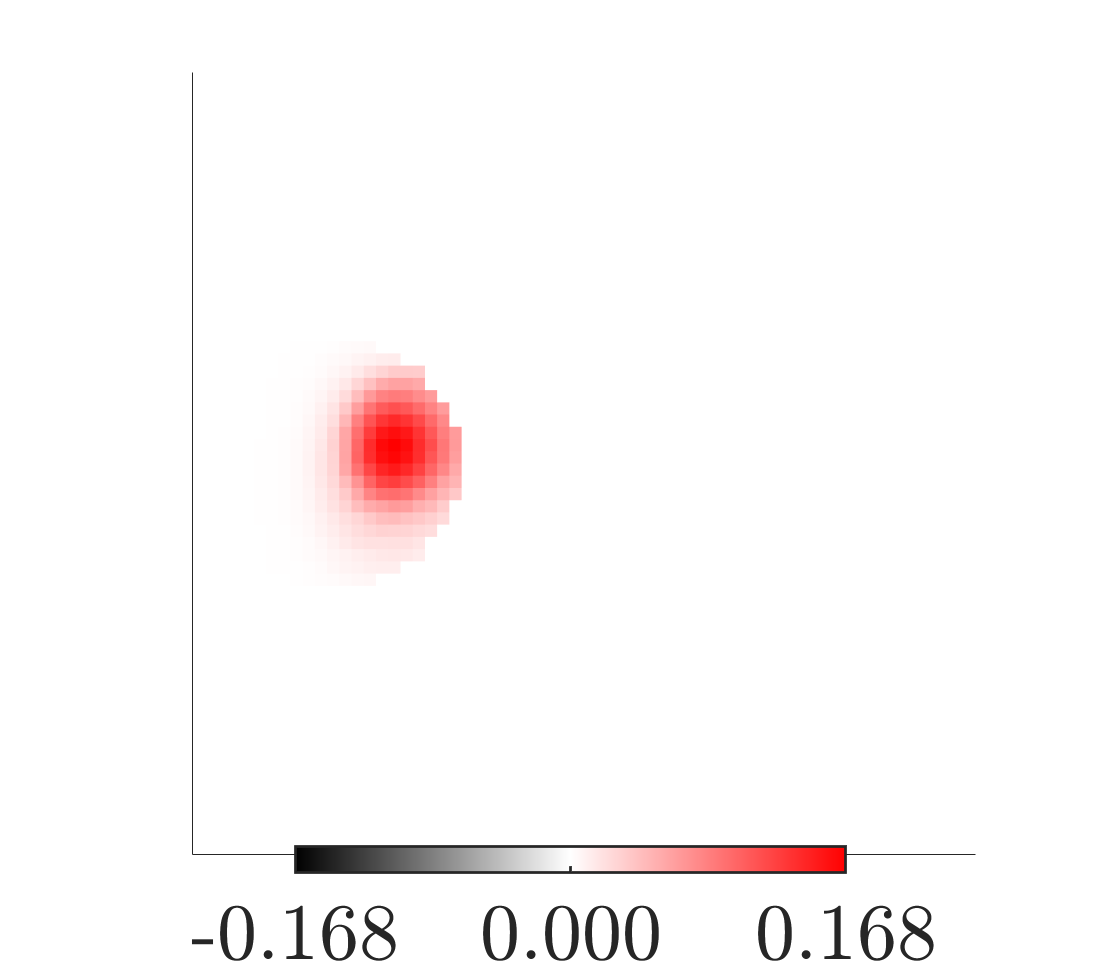} 
	\vspace{-0.25cm}
\end{minipage}
\tabularnewline 
\hline
\end{tabularx}
\caption{Lifespans detected using our three methods of identification of relevant modes; minimum equivariance, maximal variance and longest life. These results are for the case of seeding the left well when $n=10$ and $p=0.1$.}
\label{fig:main_table_case1a_n10}
\end{figure}

The results presented in Figure~\ref{fig:case_DWP_1a_tracked_lifespan_eqMM_threshold}, hint at a switching between the paths of lifespans for $\varsigma_{z}^{(3)}$ and $\varsigma_{z}^{(4)}$ (yellow and purple in electronic version) for the periodic case. %Given our task is to identify individual lifespans, as opposed to the possible connection between various lifespans through time, these anomalies do not impact our findings. 
An alternative visualisation of the lifespans detected in Figure~\ref{fig:case_DWP_1a_tracked_lifespan_eqMM_threshold} is presented in the leading row of Figure~\ref{fig:main_table_case1a_n10}. Presenting lifespans in this way allows one to concentrate on the dominance of the associated mode and the time over which a lifespan is said to exist. 
%%%%%%%%%%%%%%%%%%%%%%%%%%%%%%%%%%%%%%%%%%%%%%%%%%%%%%%%%%%%
% DWP Case 1a: n=10
%%%%%%%%%%%%%%%%%%%%%%%% 3 methods fig %%%%%%%%%%%%%%%%%%%%%
%%%%%%%%%%%%%%%%%%%%%%%%%%%%%%%%%%%%%%%%%%%%%%%%%%

We choose three methods for isolating lifespans of interest, $z_{MinEq}$, $z_{MaxVarSV}$ and $z_{Eldest}$, as defined in Algorithm~\ref{alg:3_lives}. Initial time (left) singular vectors for the three lifespans identified by these methods are presented in the two lowest rows of Figure~\ref{fig:main_table_case1a_n10}. 
{{
The leftmost column isolates a coherent structure that begins life in the left well as the time window associated with $z_{\alpha}$ opens, before shifting to the centre to merge with the right well as the time window associated with $z_{\omega}$ draws to a close.
}}
The middle column captures similar behaviour but in the reverse direction. The rightmost column identifies the time period over which a coherent structure enters the left well, inhabits that well for a period of time and later departs as the dynamics force it towards the right hand side of the domain.

These findings are in agreement with the singular value paths illustrated in Figure~\ref{fig:case_DWP_1a_tracked}. The leading path for rolling windows initialised near $t=0  \pmod{100}$ are clearly separated and nearly constant. This indicates that the identified structures are likely distinct yet neither shrinking or growing in time. In line with the singular value paths of Figure~\ref{fig:case_DWP_1a_tracked}, $z_{MinEq}$ is associated with a definitive, comparatively strong structure that persists until the associated singular values begin to fall and the lifespan draws to a close. As the singular value falls, the associated structure weakens and moves to the right, before exiting the associated patched region. This behaviour is repeated by $z_{Eldest}$. 

Additionally, $z_{Eldest}$ offers an insight into how the structure associated with the leftmost well behaves as these singular value paths begin to climb, before stabilising by $t=90$. An associated structure enters the patched region as the singular value begins to rise. As evidenced by the behaviour of $z_{MinEq}$ at $t=0$, this structure strongly persists whilst the associated singular values remain constant and near to one. As the associated structure weakens and exits the patched region the singular value paths begin to fall.
%%%%%%%%%%%%%%%%%%%%%%%%%%%%%%%%%%%%%%%%%%%%%%%%%%%%%%%%%%%%
\FloatBarrier \subsubsection{Seeding a patch over a region where two structures merge}\label{ssec:case_b} % case b) (0,0) a=b=1, n=10, depth=12
%%%%%%%%%%%%%%%%%%%%%%%%%%%%%%%%%%%%%%%%%%%%%%%%%%%%%%%%%%%%
To complement the case presented in Section~\ref{ssec:case_a}, the location of this patch centre is shifted to the origin, with all other parameters as in~\ref{ssec:case_a}. This patch characterises the region where two structures merge to one by times $50$, $150$ and so on. As per the previous example we begin our study by choosing the appropriate value for $p$. In this case both $0.2$ and $0.3$ achieve the minimum equivariance mismatch for a conservative threshold. Due to the increase in computational time that occurs as $p$ is decreased, $p=0.3$ is our preferred value. 

\begin{figure}[htb!] 
  \includegraphics[width=\textwidth]{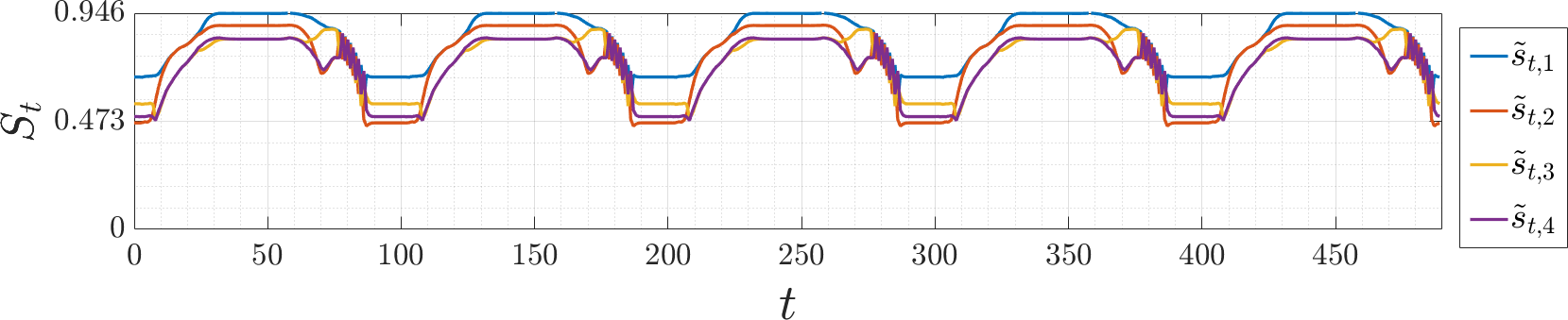}
  \caption{Rolling windows for a circular seeding of radius $1$ centred at $(0,0)$ as per Algorithm~\ref{alg:Seeding}. These are tracked using Algorithm~\ref{alg:Seeding}.}
  \label{fig:case_DWP_1b_tracked}
\end{figure}

A comparison of the peaks and troughs in rolling windows of singular values shown in Figures~\ref{fig:case_DWP_1b_tracked} and~\ref{fig:case_DWP_1a_tracked} further highlights the complementary nature of this case and the previous one. Of note are the potential signals offered by \textit{disappearing} and \textit{appearing} modes, associated with paths of rolling windows of singular values that shift towards or away from each other before evolving in parallel. The characteristic difference in this case, is that the merger is characterised by two rolling window paths that move closer before continuing in tandem. Once the merger has occurred these paths then separate and evolve parallel.

\begin{figure}[htb!]
\setlength{\tabcolsep}{-0.5pt}
\begin{tabularx}{\columnwidth}{|p{0.65cm}| *5{>{\Centering}X}|}\hline
 \multicolumn{6}{|c|}{
\begin{minipage}[t][][b]{\columnwidth}
\centering
\vspace{-0.15cm}
\includegraphics[width=0.95\columnwidth]{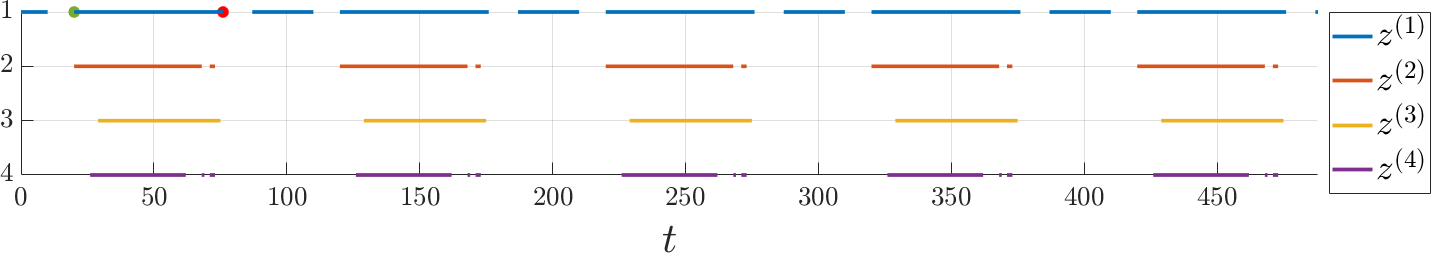}
\vspace{0.1cm}
\label{fig:case_DWP_1b_tracked_lifespan_threshold_type2}
\end{minipage} }\\ 
\hline  
\multicolumn{6}{c}{\vspace{-0.51cm} }\\
%\cline{2-4} 
\multicolumn{1}{|c|}{} &
\cellcolor{t4!60} &
\cellcolor{S1!60} &
\cellcolor{S1!60} &
\cellcolor{S1!60} &
\cellcolor{t3!60} \\[-1.em] 
\multicolumn{1}{|c|}{$t$} &
\cellcolor{t4!60} $20$ ($\tilde{u}_{20,1}$)&
\cellcolor{S1!60} $33$ ($\tilde{v}_{23,1}$)&
\cellcolor{S1!60} $39$ ($\tilde{v}_{29,1}$)&
\cellcolor{S1!60} $75$ ($\tilde{v}_{65,1}$)&
\cellcolor{t3!60} $86$ ($\tilde{v}_{76,1}$)\\
\hline
%\vspace{-0.1cm}
\begin{minipage}[t][][b]{0.03\textwidth}
\centering
\cellcolor{S1!60}
\vspace{0.85cm}
\footnotesize{$\; z^{(1)}$}
\end{minipage}
&
%\hspace{-2cm}
\begin{minipage}[t][][b]{0.2\textwidth}
	\centering
	\includegraphics[height=2.75cm,center]{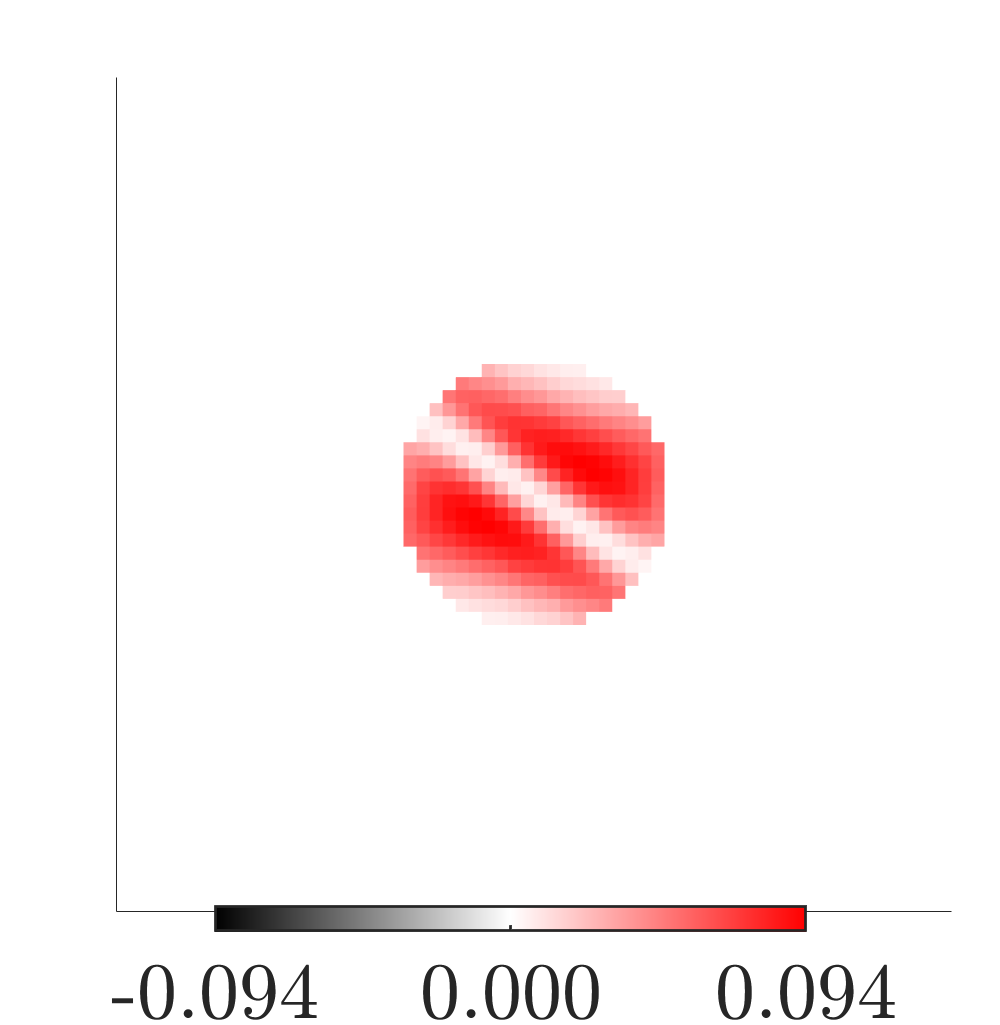} 
\vspace{-0.25cm}
\end{minipage}
&
\begin{minipage}[t][][b]{0.2\textwidth}
	\centering
	\includegraphics[height=2.75cm,center]{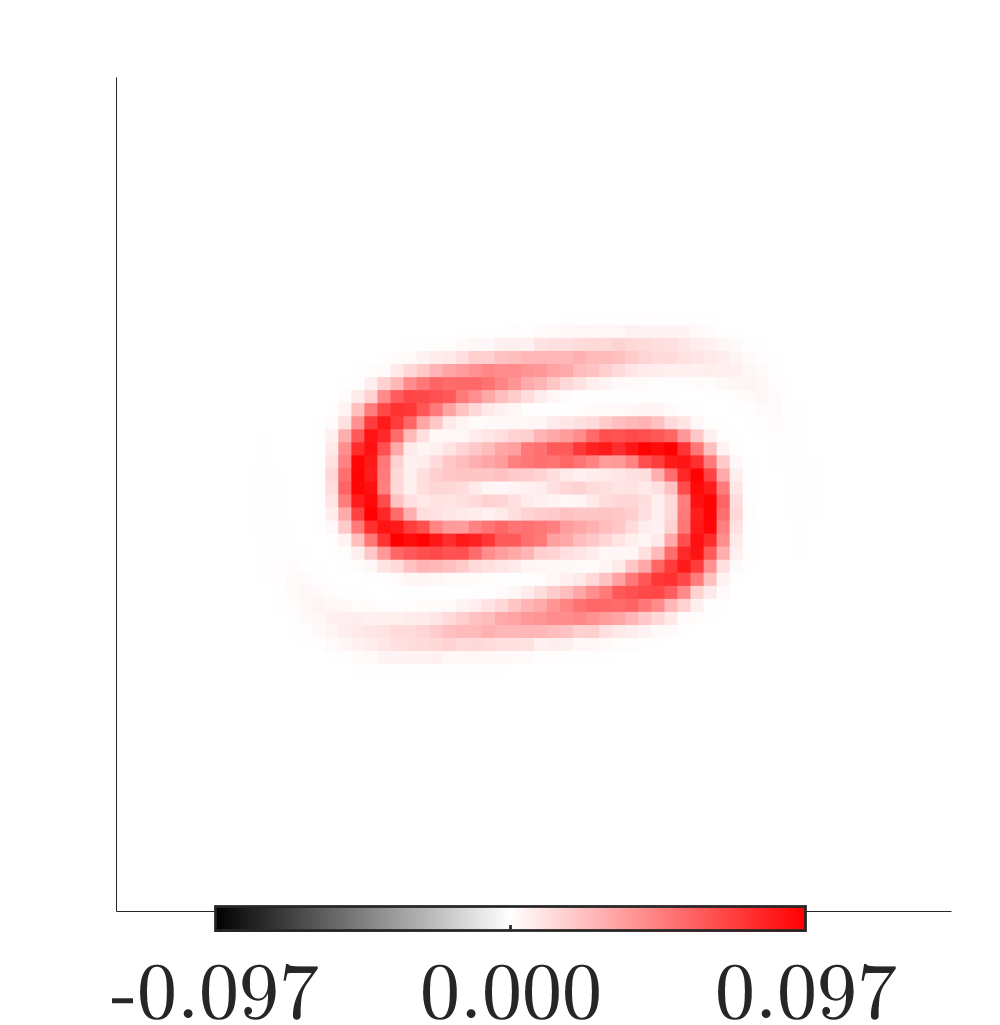} 
\vspace{-0.25cm}
\end{minipage}
&
\begin{minipage}[t][][b]{0.2\textwidth}
	\centering
	\includegraphics[height=2.75cm,center]{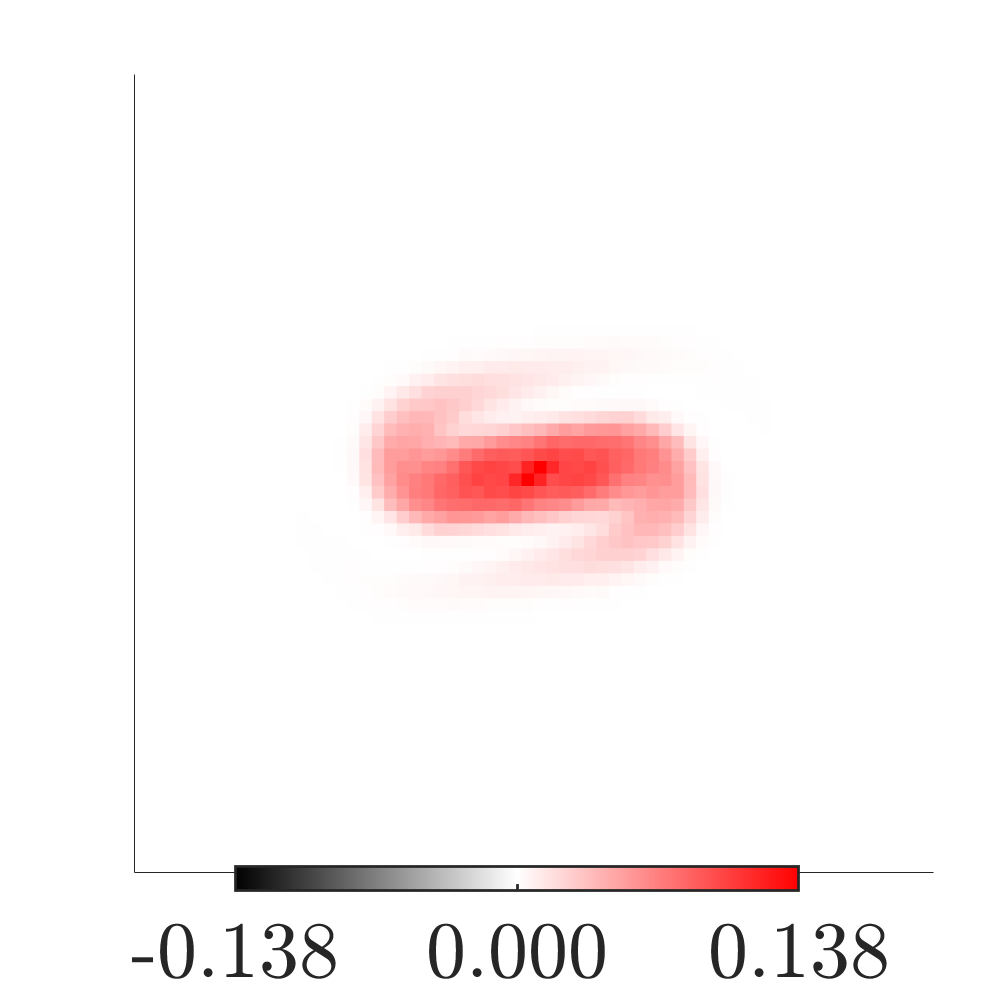} 
\vspace{-0.25cm}
\end{minipage}
&
\begin{minipage}[t][][b]{0.2\textwidth}
	\centering
	\includegraphics[height=2.75cm,center]{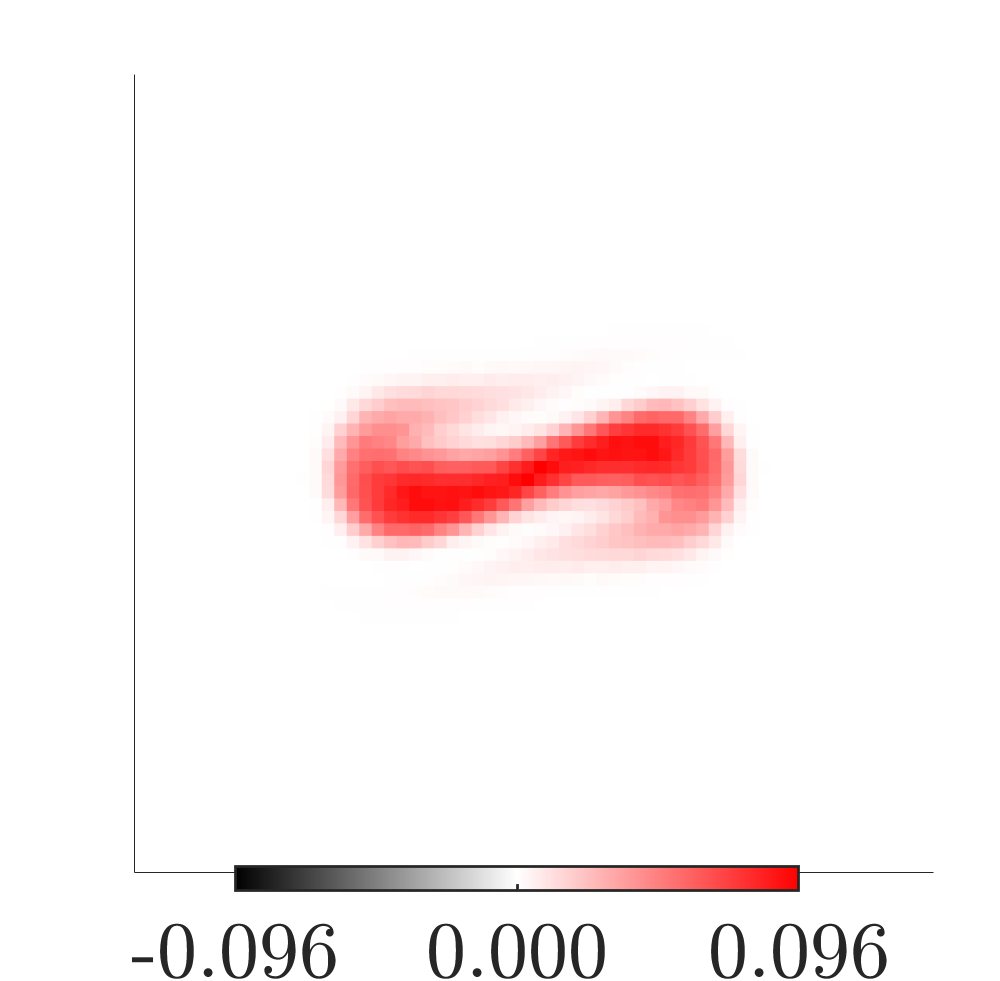} 
\vspace{-0.25cm}
\end{minipage}
&
\begin{minipage}[t][][b]{0.2\textwidth}
	\centering
	\includegraphics[height=2.75cm,center]{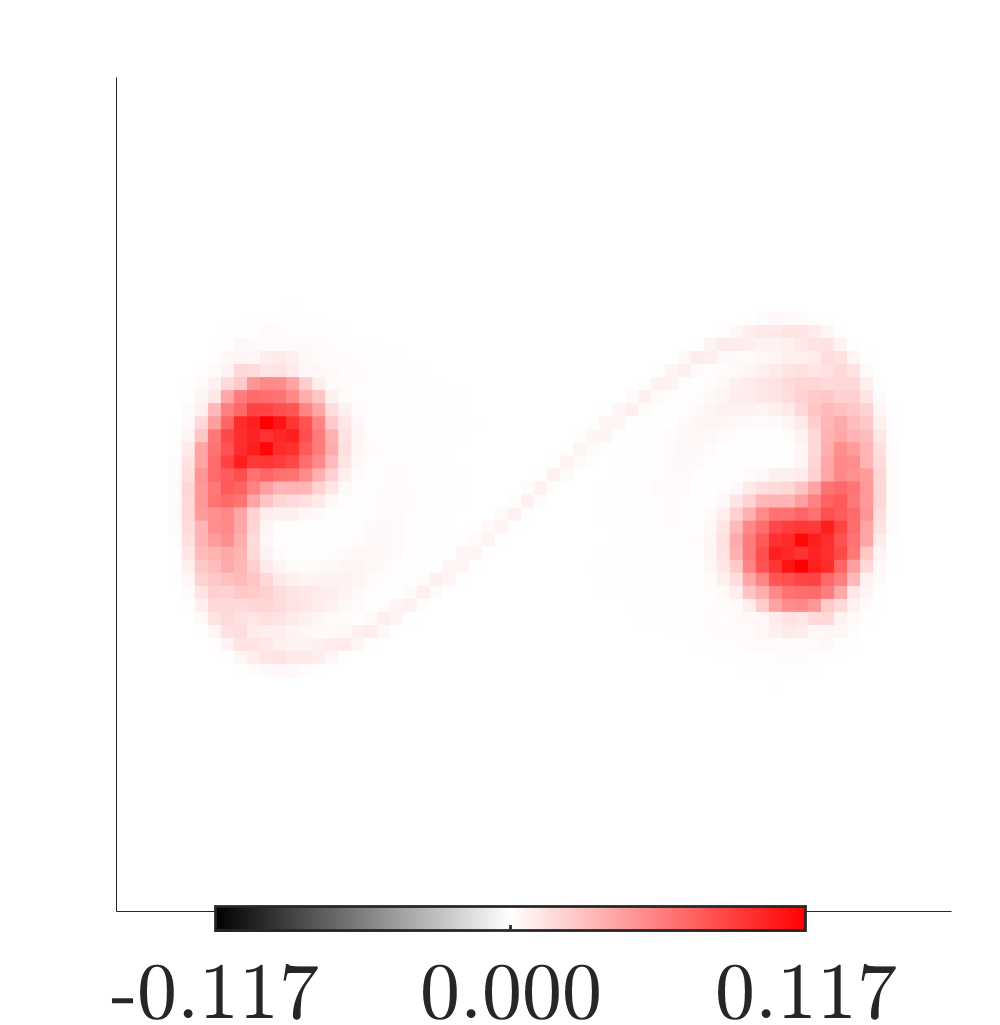} 
\vspace{-0.25cm}
\end{minipage}
\tabularnewline 
\hline
\end{tabularx}
\caption{Lifespans detected using our three methods of identification of relevant modes; minimum equivariance, maximal variance and longest life. These results are for the case of seeding the origin well when $n=10$ and $p=0.3$. One notes that $\tilde{u}_{76,1}$ is the left singular vector of associated with the time window initialised at $t=76$ for $n=10$ and the lifespan $z^{(1)}$. The corresponding right singular vector is denoted by $\tilde{v}_{76,1}$, this is associated with the dynamics at time $t=8{6}$.}
\label{fig:main_table_case1b_n10}
\end{figure}

Figure~\ref{fig:main_table_case1b_n10} presents the lifespans of interest identified by Algorithm~\ref{alg:3_lives}. In this case, all three methods identify the same period. This period is characterised by two structures entering the patched region from the left and right before merging into a single structure and then separating as they exit the patch.
%%%%%%%%%%%%%%%%%%%%%%%%%%%%%%%%%%%%%%%%%%%%%%%%%%%%%%%%%%%%%%%%%%%%%%%%%%%%%%%%%%%%%%%%%%%%%%%%%%%% DWP Case 1b: n=10
%%%%%%%%%%%%%%%%%%%%%%%%%%%%%%%%%%%%%%%%%%%%%%%%%%%%%%%%%%%%%%%%%%%%%%%%%%%%%%%%%%%%%%%%%%%%%%%%%%%%%%%%%%%%%%%%%%%%%%%%%%%% 3 methods fig %%%%%%%%%%%%%%%%%%%%%

As demonstrated by the final row in Figure~\ref{fig:main_table_case1b_n10}, our methods detect the entire lifespan of two structures that merge into one before again separating. In this sense the lifespan detected is a full lifecycle that includes the birth (following $z_{\alpha}$) and the death ($z_{\omega}$) of the associated structure formed by two previously distinct structures. This single structure is illustrated by $\tilde{v}_{{2}9,1}$ in Figure~\ref{fig:main_table_case1b_n10} whilst $\tilde{v}_{23,1}$ illustrates the beginnings of this merger between two entities. Completing the lifecycle is $\tilde{v}_{76,1}$, which illustrates how the newly formed structure separates into two.
%%%%%%%%%%%%%%%%%%%%%%%%%%%%%%%%%%%%%%%%%%%%%%%%%%%%%%%%%%%%%%%%%%%%%%%%%%%%%%%%
\FloatBarrier \subsubsection{Seeding a chaotic region}\label{ssec:case_c} % (-2,-2) a=b=1, n=10, depth=12
%%%%%%%%%%%%%%%%%%%%%%%%%%%%%%%%%%%%%%%%%%%%%%%%%%%%%%%%%%%%%%%%%%%%%%%%%%%%%%%%
One does not expect to locate meaningful coherent structures in regions characterised by extreme mixing or chaotic behaviour. For comparison with previous seedings, let us examine the lifespans identified when one seeds an area where neither well is centred for any length of time. We centre a patch at $(-2,-2)$ of radius $1$ for rolling windows of length $n=10$. For all $p \in \mathcal{P}$ this case never achieves an average equivariance mismatch below the conservative threshold. The consequence of this is that no lifespans are detected. This means that our methods indicate no meaningful coherent structures are identified when this circular region  is seeded and flowed for $10$ time steps, as expected.

%%%%%%%%%%%%%%%%%%%%%%%%%%%%%%%%%%%%%%%%%%%%%%%%%%%%%%%%%%%%%%%%%%%%%%%%%
\FloatBarrier \subsection{Boussinesq equation models}\label{SSec:BEeModels_Results}
%{Modon eddy pairs (Dipoles)}
%%%%%%%%%%%%%%%%%%%%%%%%%%%%%%%%%%%%%%%%%%%%%%%%%%%%%%%%%%%%%%%%%%%%%%%%%
%\label{SSec:RICs}\label{SSec:MMpole}\label{SSec:DDip}
Let us now test these methods on the more complex datasets generated by the Boussinesq equations. This allows for greater insight into the detection of lifespans and associated signals in more complex environments. Again we begin with the construction of localised Ulam matrices, as described in Algorithm~\ref{alg:Seeding}, setting $Q=100$. The discrete time flow maps are approximated using Runge-Kutta numerical integration, in space and time, for $\tau$ equal to the time scales noted in Sections~\ref{SSec:DDip},~\ref{SSec:MMpole} and~\ref{SSec:RICs} over $10$ steps. Our numerical analysis is limited to the midplane $z=\pi$ and calculations are performed on velocity fields with a grid of dimension $192 \times 192$ in the $x$, and $y$ directions at each inertial period. Given the finer dynamics of these models, we increase resolution to a depth of $14$ whilst continuing to analyse $n=10$ matrices.

Section~\ref{ssec:case_dipoles} considers the evolution of dipole pairs, as described in Section~\ref{SSec:DDip}. This model is utilised to patch a region of phase space where coherent structures are known to visit. Section~\ref{SSec:MMs} considers the case of merging monopoles, as described in Section~\ref{SSec:MMpole}. In this case, the characteristic merging event is easily detected using lifespans obtained from Algorithms~\ref{alg:3_lives} and~\ref{alg:CSorNot}. Section~\ref{SSec:RICs_Results} explores the more chaotic example characterised by interacting, evolving objects that take a variety of shapes, as described in Section~\ref{SSec:RICs}. In spite of the increased complexity of this model, Algorithms~\ref{alg:3_lives} and~\ref{alg:CSorNot} clearly identify a merging of structures in the presence of large amounts of background noise.

%%%%%%%%%%%%%%%%%%%%%%%%%%%%%%%%%%%%%%%%%%%%%%%%%%
\FloatBarrier \subsubsection[Dancing dipole pairs]{Dancing dipole pairs}\label{ssec:case_dipoles}
\begin{figure}[!htbp]
\centering
\begin{subfigure}[b]{\textwidth}
\includegraphics[width=\textwidth]{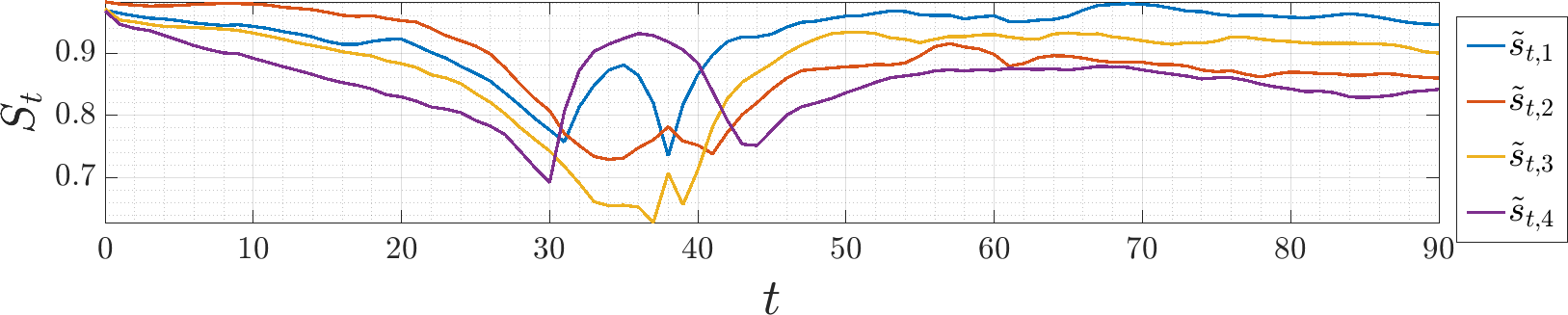}
\caption{\label{fig:case_DipG_c} Leading $4$ rolling windows of singular values for the dancing dipoles model. These are paired through time using Algorithm~\ref{alg:Seeding} with a circular patch of radius $1/2$ centred at $(5.55,3.75)$ for $n=10$. Tracking is achieved using $p=0.8$.}
\end{subfigure}
\hfill
\begin{subfigure}[b]{\textwidth}
\includegraphics[width=\textwidth]{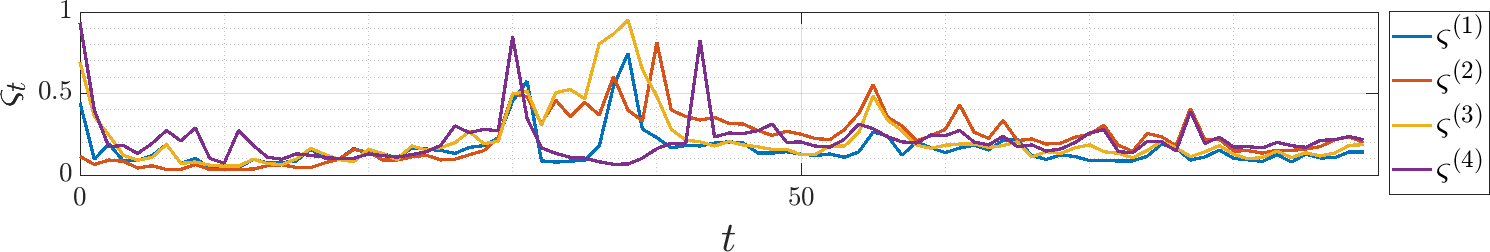}
\caption{\label{fig:case_DipG_c_EqMM} Equivariance mismatch, as described in Algorithm~\ref{alg:lifespan}, for the case associated with Figure~\ref{fig:case_DipG_c} where tracking was achieved using $p={0.8}$. Figure~\ref{fig:case_DipG_c} and~\ref{fig:case_DipG_c_EqMM} both give $t$ in terms of inertial periods.}
\end{subfigure}
\caption{Rolling windows of singular value paths and equivariance mismatch for the dancing dipoles numerical test.}
\end{figure}

To detect the presence of coherent structures in a patched region, we seed a circular patch of radius $1/2$ centred at $(5.55,3.75)$. It is evident from Figure~\ref{fig:TwoDipoles}, that this region is likely to include passing structures around time $30$. The corresponding rolling windows of singular value paths for {$p=0.8$} are shown in Figure~\ref{fig:case_DipG_c}. Figure~\ref{fig:case_DipG_c_EqMM} plots the associated equivariance mismatch.
%%%%%%%%%%%%%%%%%%%%%%%%%%%%%%%%%%%%%%%%%%%%%%%%%%
%%%%%%%%%%%%%%%%%%%%%%%%%%%%%%%%%%%%%%%%%%%%%%%%%%

%%%%%%%%%%%%%%%%%%%%%%%% 3 methods fig %%%%%%%%%%%%%%%%%%%%%
\begin{figure}[!htbp]
\setlength{\tabcolsep}{-0.75pt}
\begin{tabularx}{\columnwidth}{|p{0.65cm}| *4{>{\Centering}X}|}\hline
\multicolumn{5}{|c|}{
\begin{minipage}[t][][b]{\columnwidth}
\centering
\vspace{-0.15cm}
\includegraphics[width=0.95\columnwidth]{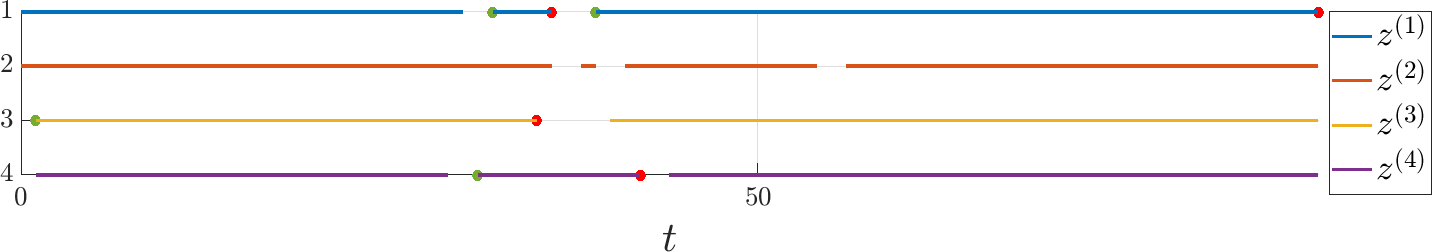}
\vspace{0.1cm}
\label{fig:case_DWP_1a_tracked_lifespan_threshold_type2_Dip_c}
\end{minipage} }
\\ \hline  
% row 2
\multicolumn{1}{c|}{} &
{\cellcolor{S1!60}}&
{\cellcolor{S4!60}}&
{\cellcolor{S3!60}}&
{\cellcolor{S1!60}}\\[-1.em] 
\multicolumn{1}{c|}{} &
{\cellcolor{S1!60} $z_{MinEq,1}$ for $t \in [32,36]$} &
{\cellcolor{S4!60} $z_{MinEq,2}$ for $t \in [31,42]$} & 
{\cellcolor{S3!60} $z_{MaxVarSV}$ } for $t \in [1,35]$ & 
{\cellcolor{S1!60} $z_{Eldest}$ } for $t \in [39,88]$\\  
\cline{2-5} \multicolumn{5}{c}{\vspace{-0.5cm}}\\
\hline
% row 4
\begin{minipage}[t][][b]{0.03\textwidth}
\centering
\vspace{1.9cm}
\footnotesize{$\;\tilde{u}$}
\end{minipage}
&
\begin{minipage}[t][][c]{0.24\textwidth}
	\centering
	\includegraphics[height=3cm,center]{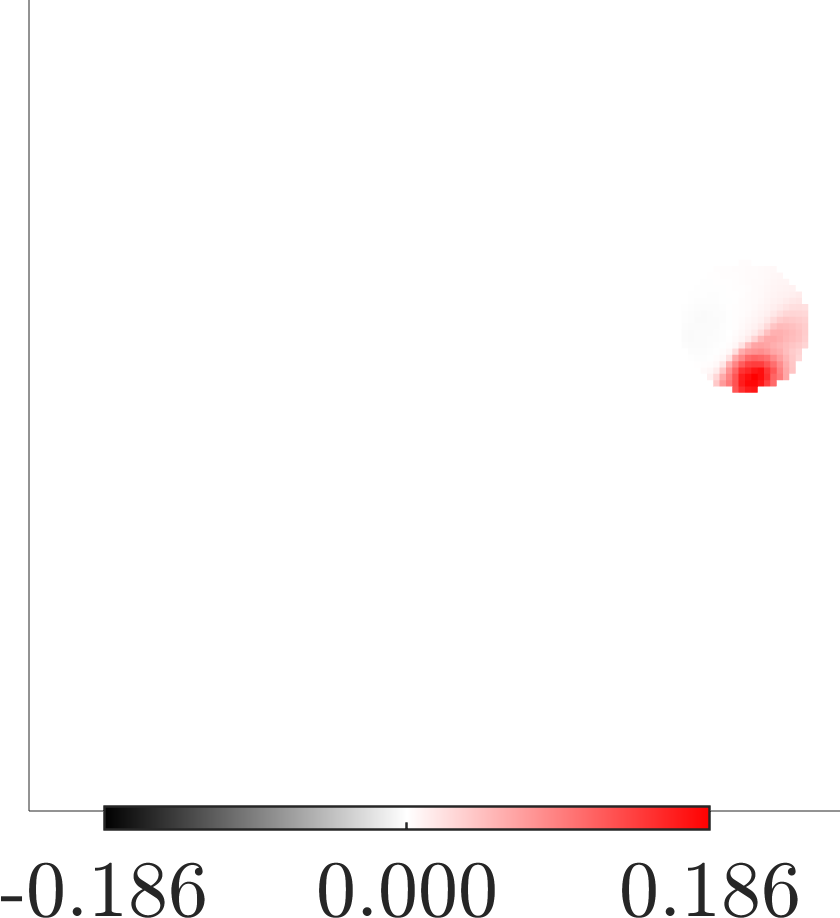} 
	\vspace{-0.2cm}
\end{minipage}
&
\begin{minipage}[t][][b]{0.24\textwidth}
	\centering
	\includegraphics[height=3cm,center]{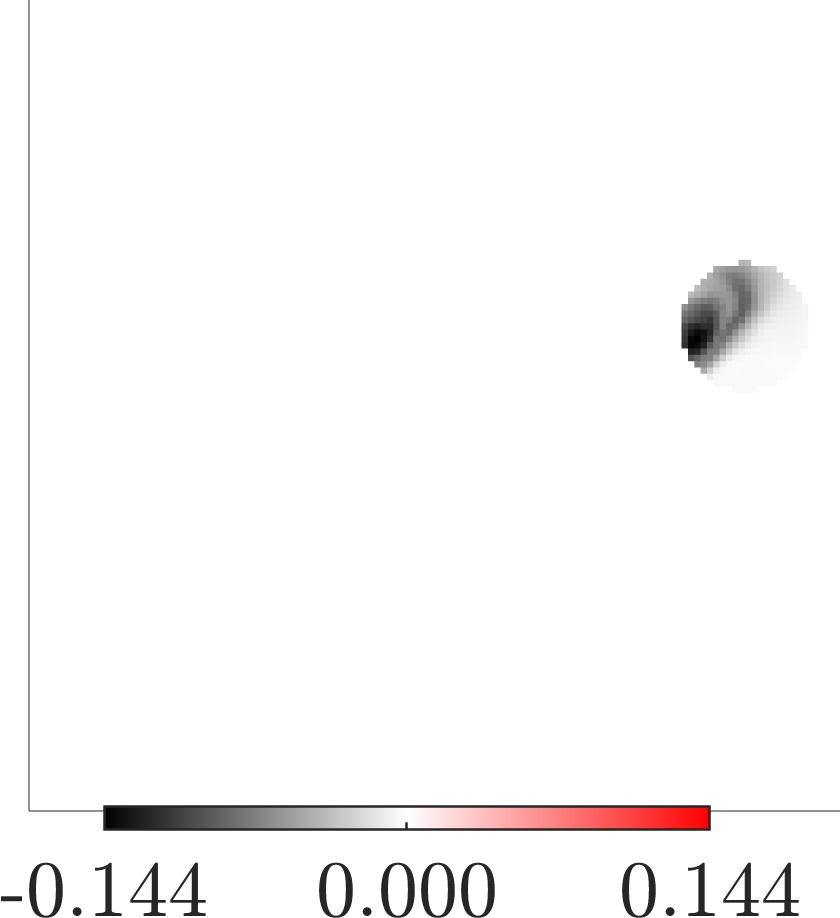} 
\end{minipage}
&
\begin{minipage}[t][][c]{0.24\textwidth}
	\centering
	\includegraphics[height=3cm,center]{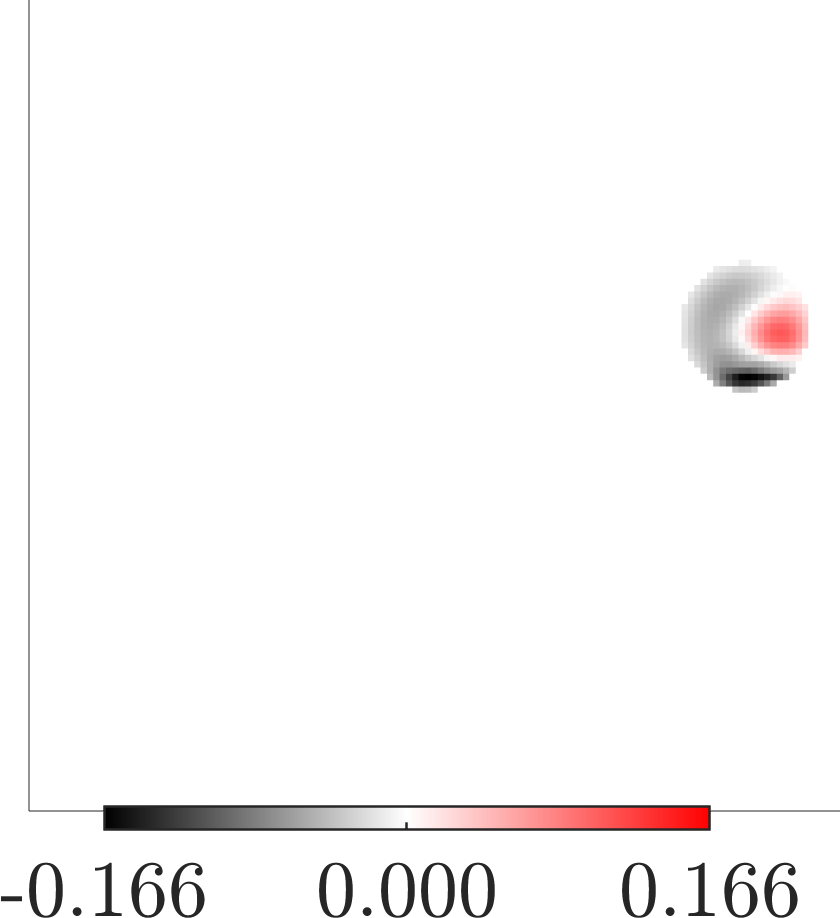}
	\vspace{-0.2cm}
\end{minipage}\hfil
&
\begin{minipage}[t][][b]{0.24\columnwidth}
	\centering
	\includegraphics[height=3cm,center]{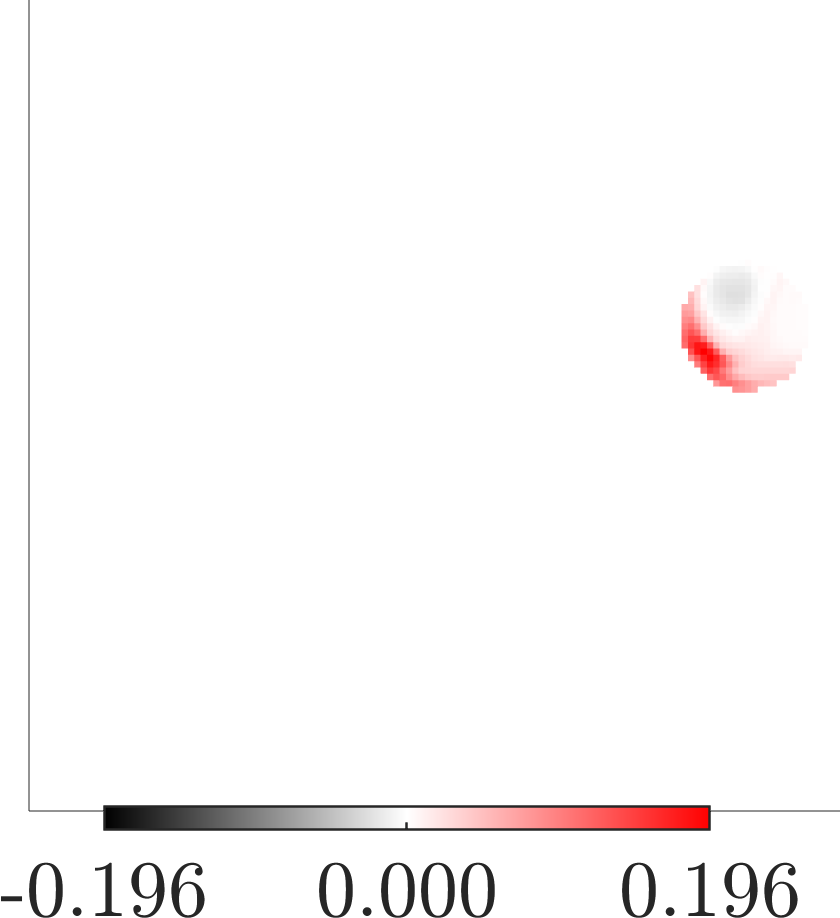}
\end{minipage}
\\  \hline
\begin{minipage}[t][][c]{0.03\textwidth}
\centering
\vspace{1.9cm}
\footnotesize{$\;\tilde{v}$}
\end{minipage}
&
\begin{minipage}[t][][c]{0.24\textwidth}
	\centering
	\includegraphics[height=3cm,center]{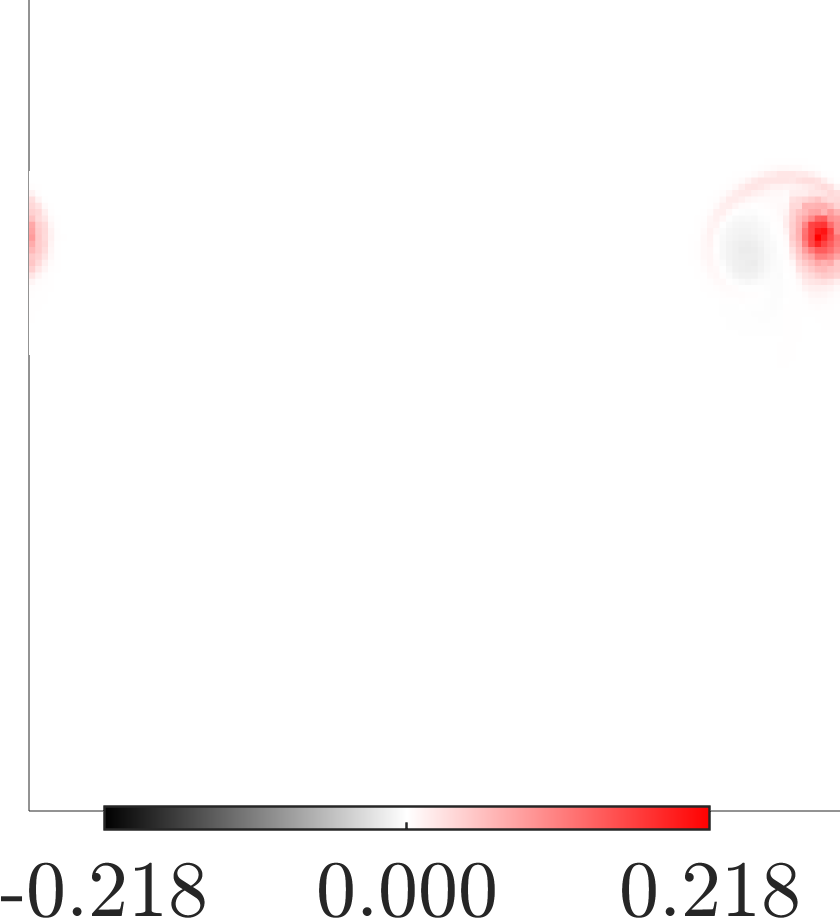}
	\vspace{-0.2cm}
\end{minipage}\hfil
&
\begin{minipage}[t][][b]{0.24\columnwidth}
	\centering
	\includegraphics[height=3cm,center]{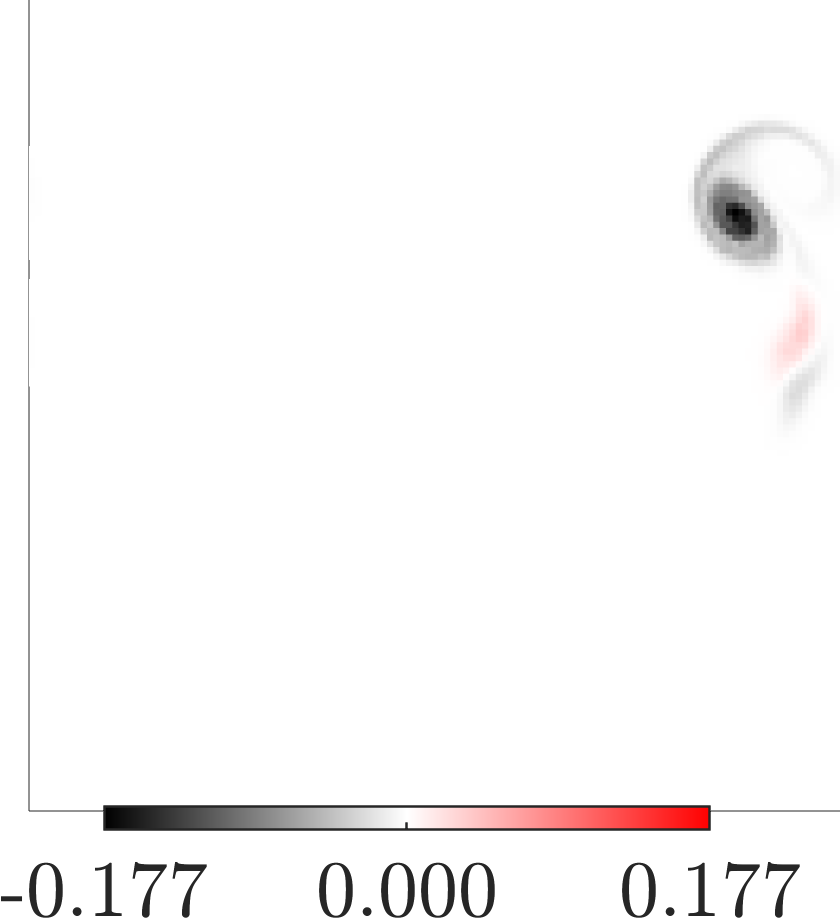} 
\end{minipage}
&
\begin{minipage}[t][][c]{0.24\textwidth}
	\centering
	\includegraphics[height=3cm,center]{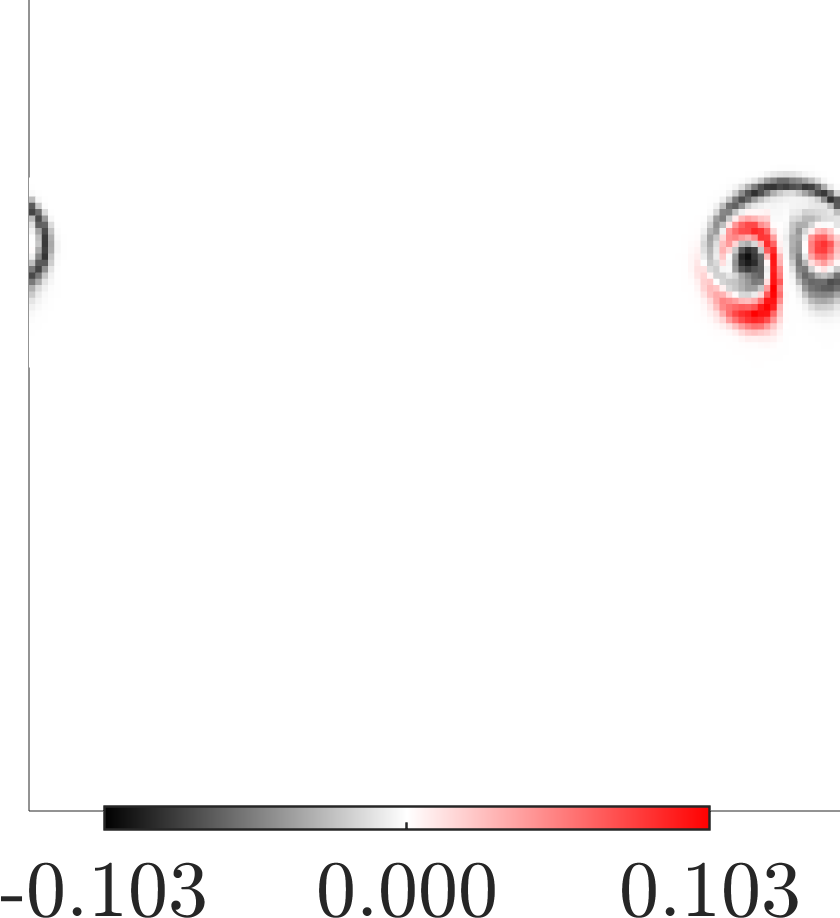}
	\vspace{-0.2cm}
\end{minipage}\hfil
&
\begin{minipage}[t][][b]{0.24\columnwidth}
	\centering
	\includegraphics[height=3cm,center]{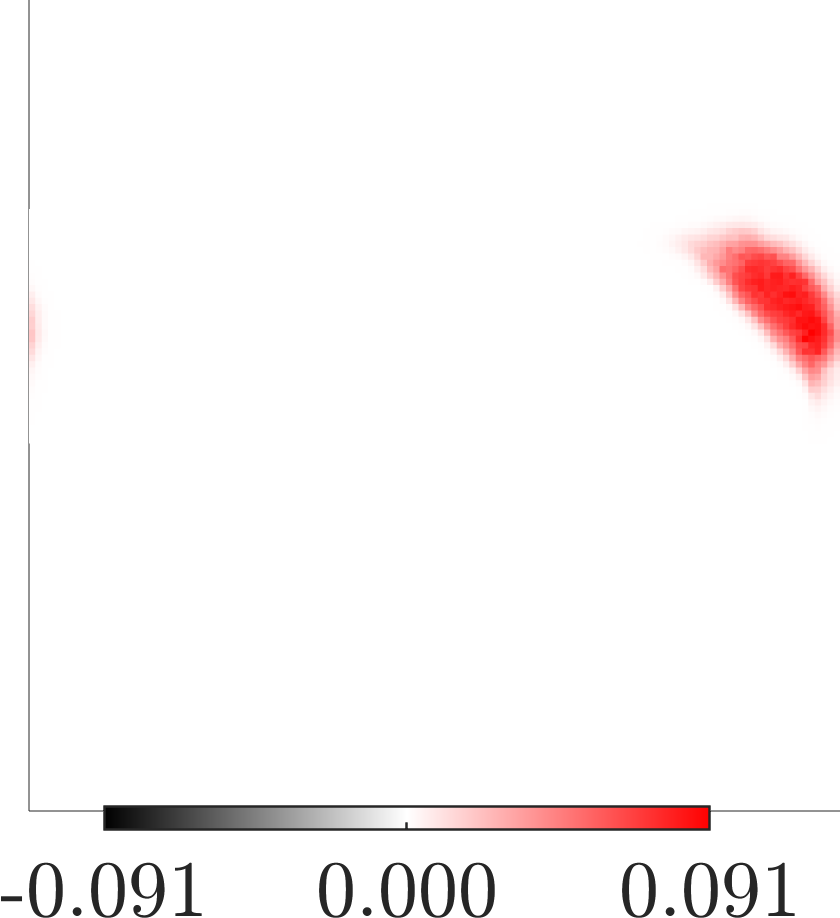}
\end{minipage}
\\ \hline
\end{tabularx}
\caption{General lifespans (top row) as well as the lifespans $z_{MinEq}$, $z_{MaxVarSV}$ and $z_{Eldest}$, as detected in the dancing dipoles model using Algorithms~\ref{alg:lifespan} and~\ref{alg:3_lives}, for $t$ in terms of inertial periods. These results are for the case presented in Figure~\ref{fig:case_DipG_c}. Relevant left singular vectors of initial time windows ($\tilde{u}$) are plotted in the middle row ($z_{\alpha}$) with right singular vectors ($\tilde{v}$) plotted in the final row for $z_{\omega}$.}
\label{fig:main_table_DipG_c_n10}
\end{figure}
%%%%%%%%%%%%%%%%%%%%%%%%%%%%%%%%%%%%%%%%%%%%%%%%%%
%%%%%%%%%%%%%%%%%%%%%%%%%%%%%%%%%%%%%%%%%%%%%%%%%%
\begin{figure}[!htbp]
\setlength{\tabcolsep}{-0.5pt}
\begin{tabularx}{\columnwidth}{|p{0.65cm}| *5{>{\Centering}X}|}
\cline{2-6} \multicolumn{6}{c}{\vspace{-0.5cm}}\\
\multicolumn{1}{c}{} &
\multicolumn{1}{|c}{\cellcolor{S4!60}} &
\multicolumn{1}{|c}{\cellcolor{S4!60}} &
\multicolumn{1}{|c}{\cellcolor{S4!60}} &
\multicolumn{1}{|c}{\cellcolor{S4!60}} &
\multicolumn{1}{|c|}{\cellcolor{S4!60}} \\[-1.em] 
\multicolumn{1}{c}{} &
\multicolumn{1}{|c}{\cellcolor{S4!60} $\tilde{u}_{30,4}$}&
\multicolumn{1}{|c}{\cellcolor{S4!60} $\tilde{u}_{32,4}$}&
\multicolumn{1}{|c}{\cellcolor{S4!60} $\tilde{u}_{42,4}$}&
\multicolumn{1}{|c}{\cellcolor{S4!60} $\tilde{v}_{32,4}$}&
\multicolumn{1}{|c|}{\cellcolor{S4!60} $\tilde{v}_{41,4}$}\\
\hline
%\vspace{-0.1cm}
\begin{minipage}[t][][b]{0.03\textwidth}
\centering
\cellcolor{S4!60}
\vspace{0.85cm}
\footnotesize{$\; \tilde{s}_{t,4}$}
\end{minipage}
&
%\hspace{-2cm}
\begin{minipage}[t][][b]{0.2\textwidth}
	\centering
	\includegraphics[height=2.75cm,center]{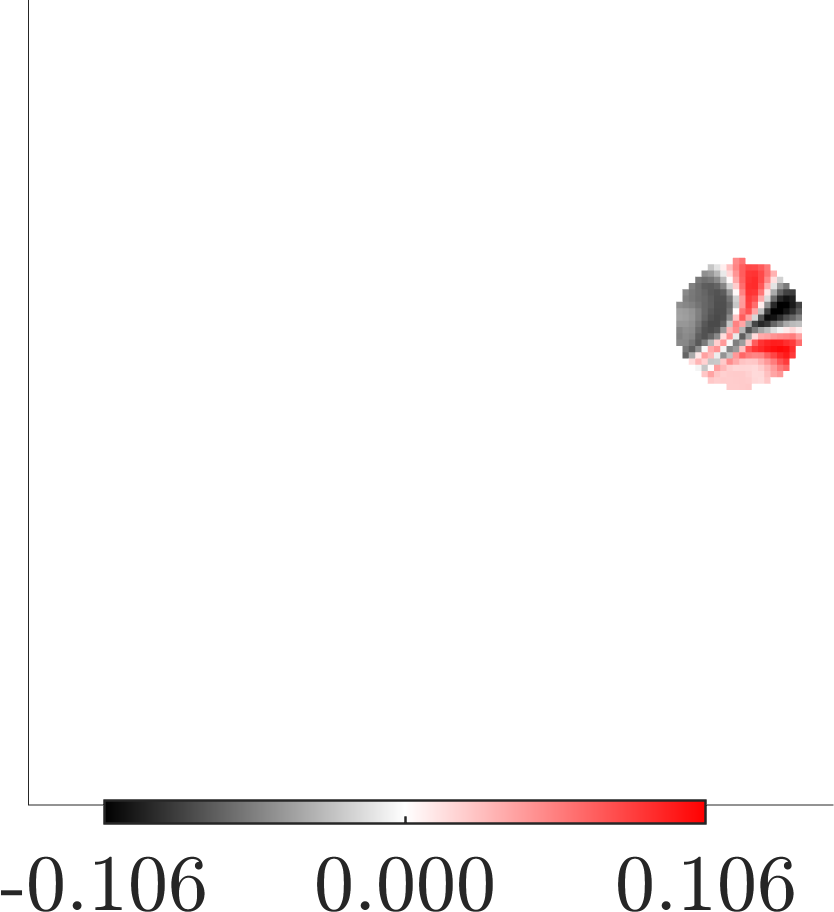} 
\vspace{-0.25cm}
\end{minipage}
&
\begin{minipage}[t][][b]{0.2\textwidth}
	\centering
	\includegraphics[height=2.75cm,center]{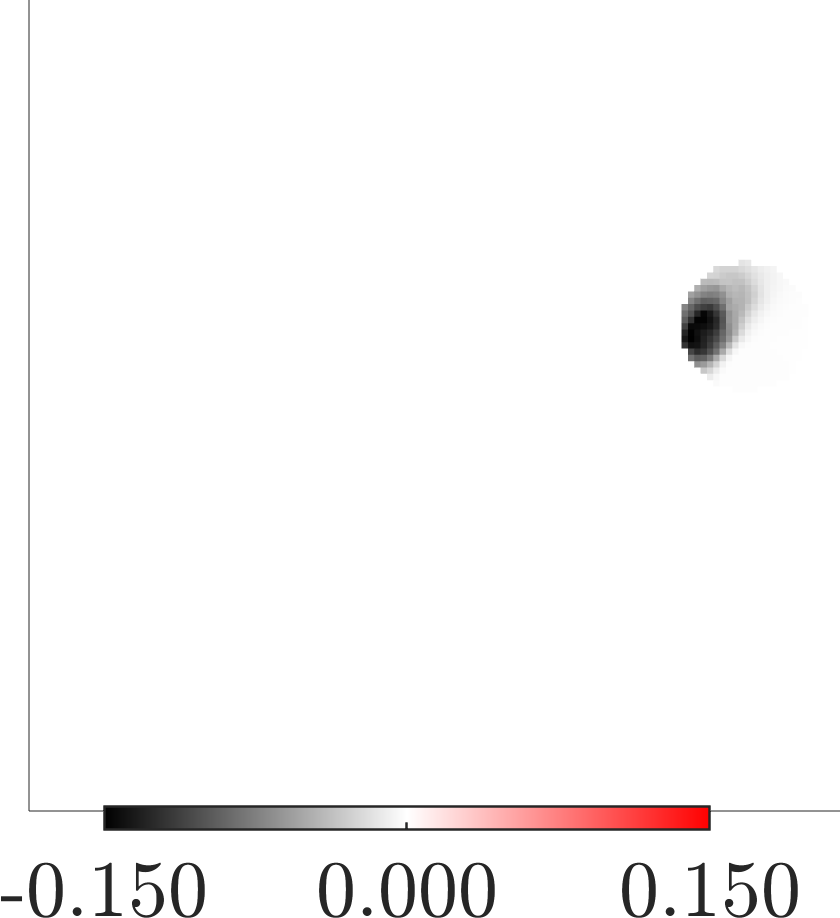} 
\vspace{-0.25cm}
\end{minipage}
&
\begin{minipage}[t][][b]{0.2\textwidth}
	\centering
	\includegraphics[height=2.75cm,center]{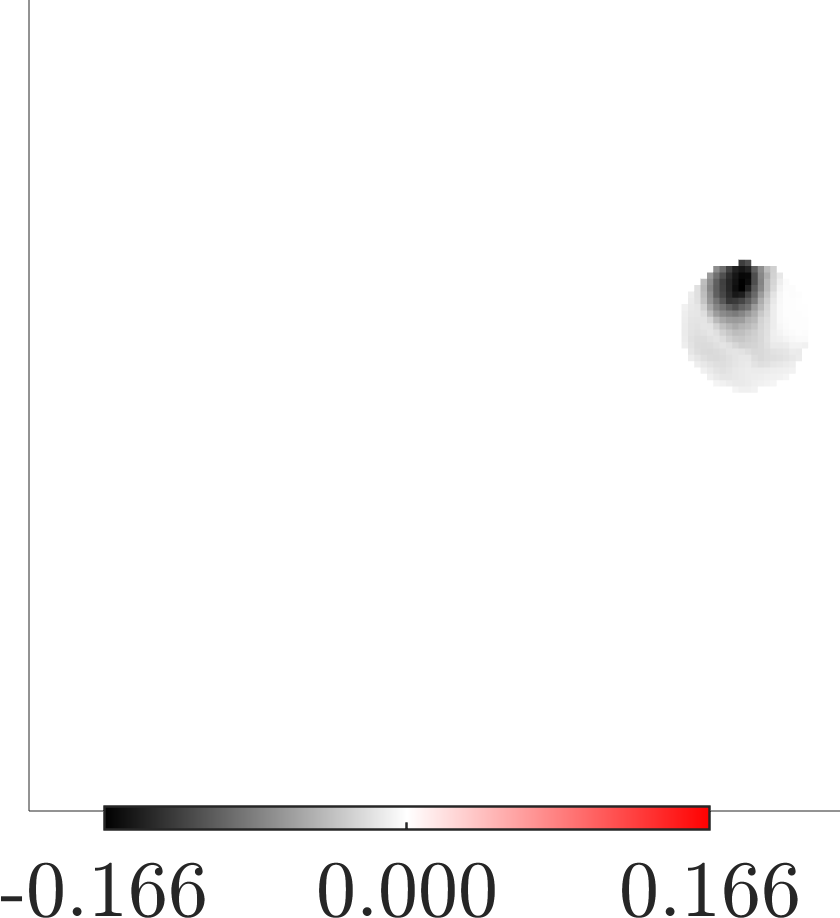} 
\vspace{-0.25cm}
\end{minipage}
&
\begin{minipage}[t][][b]{0.2\textwidth}
	\centering
	\includegraphics[height=2.75cm,center]{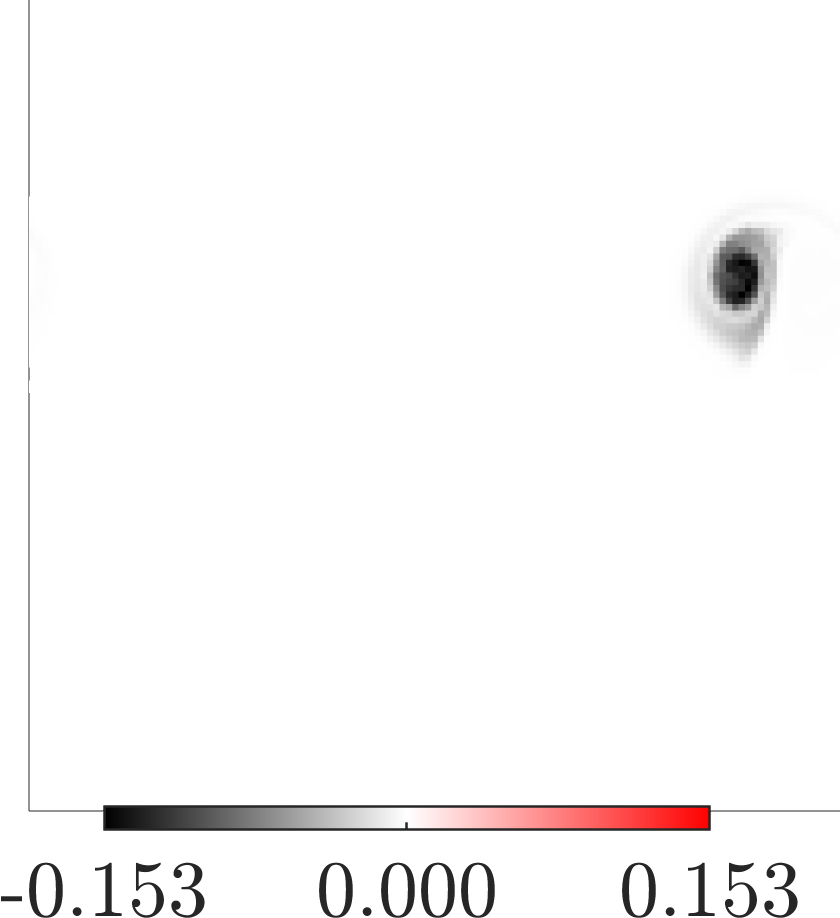} 
\vspace{-0.25cm}
\end{minipage}
&
\begin{minipage}[t][][b]{0.2\textwidth}
	\centering
	\includegraphics[height=2.75cm,center]{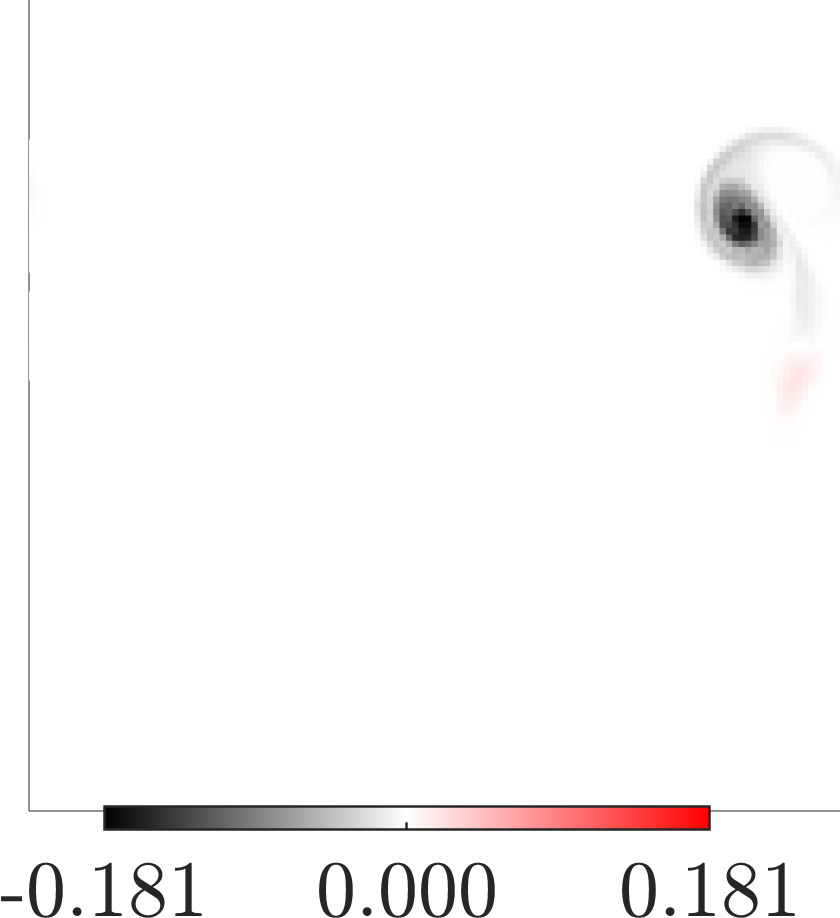} 
\vspace{-0.25cm}
\end{minipage}
\tabularnewline 
\hline
\end{tabularx}
\caption{Detailed view of the behaviour of $\{ \tilde{u}_{t,4}\}$ and $\{ \tilde{v}_{t,4}\}$ for a variety of $t$ inertial periods. This figure partially aligns with the peak in $\{ \tilde{s}_{t,4} \}$ of Figure~\ref{fig:case_DipG_c} and the lifespan identified by $z_{MinEq,2}$.}
\label{fig:minEq_caseStRotn_n10}
\end{figure}
%%%%%%%%%%%%%%%%%%%%%%%%%%%%%%%%%%%%%%%%%%%%%%%%%%
It is clear from the peaks in Figure~\ref{fig:case_DipG_c} that the paths of singular values are again capable of detecting the presence of dynamic structures within a patch. One notes that both peaks in $\{ \tilde{s}_{t,{1}} \}$ and $\{ \tilde{s}_{t,4} \}$ arise as coherent structures enter the patched region and dissipate as they leave. Figure~\ref{fig:case_DipG_c_EqMM} shows that vectors $\{\tilde{v}\}$ are paired somewhat consistently on either sides of these peaks. Moreover, as the peaks occur $\{\tilde{v}_{t,{1}}\}$ and $\{\tilde{v}_{t,4}\}$ achieve comparatively low values of equivariance mismatch $\varsigma$ despite the concomitant changes occurring in $\{ \tilde{s}_{t,{1}} \}$ and $\{ \tilde{s}_{t,4} \}$. This shows that fundamental changes are occurring in the patched region and that $\{\tilde{v}_{t,{1}}\}$ and $\{\tilde{v}_{t,4}\}$ are consistently paired as these changes occur.

Figure~\ref{fig:main_table_DipG_c_n10} illustrates the identified lifespans alongside a visualisation of two lifespans of interest. These were detected using our general method, as described in Algorithms~\ref{alg:lifespan} and~\ref{alg:3_lives}, using a conservative threshold. 
It is inferred, from examining the evolution of the singular vectors shown in Figures~\ref{fig:main_table_DipG_c_n10} and~\ref{fig:minEq_caseStRotn_n10} in association with the vector field, that the larger peak ($\{ \tilde{s}_{t,4} \}_{ 30 \le t \le 43}$) is associated with the passage of the newly merged upper poles through the patched region, whilst the smaller peak ($\{ \tilde{s}_{t,1} \}_{31 \le t \le 38}$) is predominantly associated with the movement of the rightmost lower pole through the same area.
Furthermore, the lifespan $z_{MinEq,1}$ captures the early movement of this lower pole through the patched region. In addition, the lifespan $z_{MinEq,2}$ is associated with the evolution of the two upper poles following their merger. That is, the peak in $\{ \tilde{s}_{t,4} \}$ of Figure~\ref{fig:case_DipG_c}. Figure~\ref{fig:minEq_caseStRotn_n10} illustrates that as the path defined by $\{ \tilde{s}_{t,4} \}$ begins to separate from the other paths traced by singular values in Figure~\ref{fig:case_DipG_c}, the mode associated with the lifespan $z_{MinEq,2}$ is distinguished and better isolates the dynamical behaviour of interest. 

Conversely, as shown in the final two columns of Figure~\ref{fig:main_table_DipG_c_n10}, the lifespans identified with $z_{MaxVarSV}$ and $z_{Eldest}$ are more spurious and less informative. This is influenced by the concentrated dynamical behaviour characterising this model, which sees large portions of the domain remain stagnant for the $t$ considered.

\begin{figure}[!htbp]
\centering
\begin{subfigure}[c]{\columnwidth}
\hspace{0.8cm}
\includegraphics[width=0.8\textwidth]{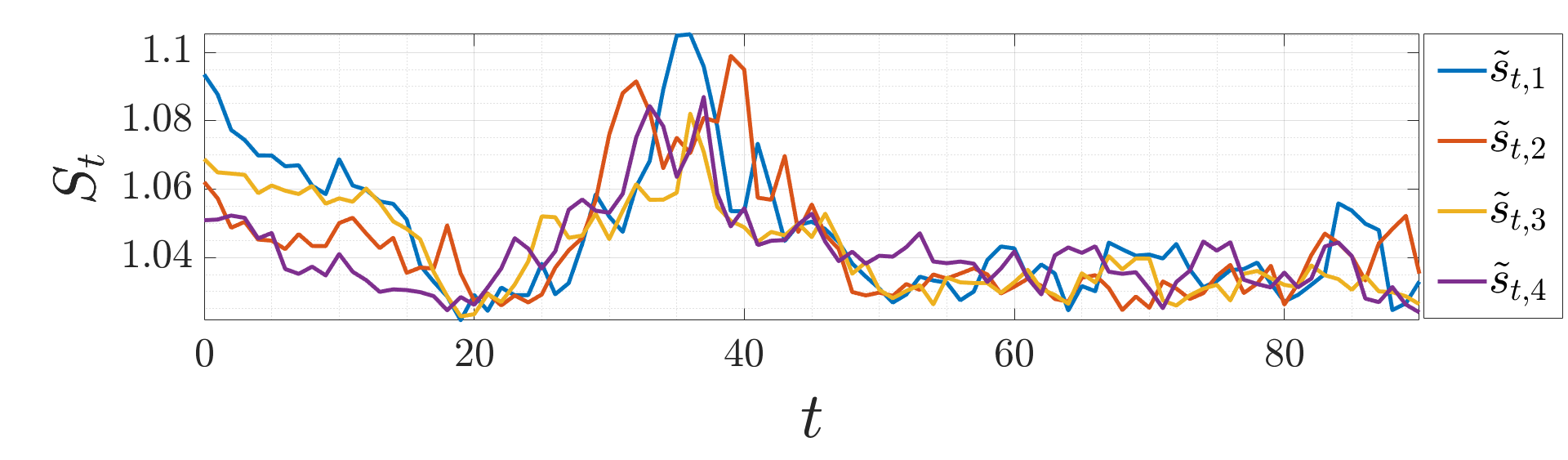}
\caption{\label{fig:case_Dip_FullTO} Leading $4$ rolling windows of singular values for the dancing dipoles model for the full transfer operator. These are paired through time using Algorithm~\ref{alg:Seeding} for $n=10$ and $p=2$.}
\end{subfigure}
\hfill
\begin{subfigure}[c]{\columnwidth}
\hspace{0.8cm}
\includegraphics[width=0.8\textwidth]{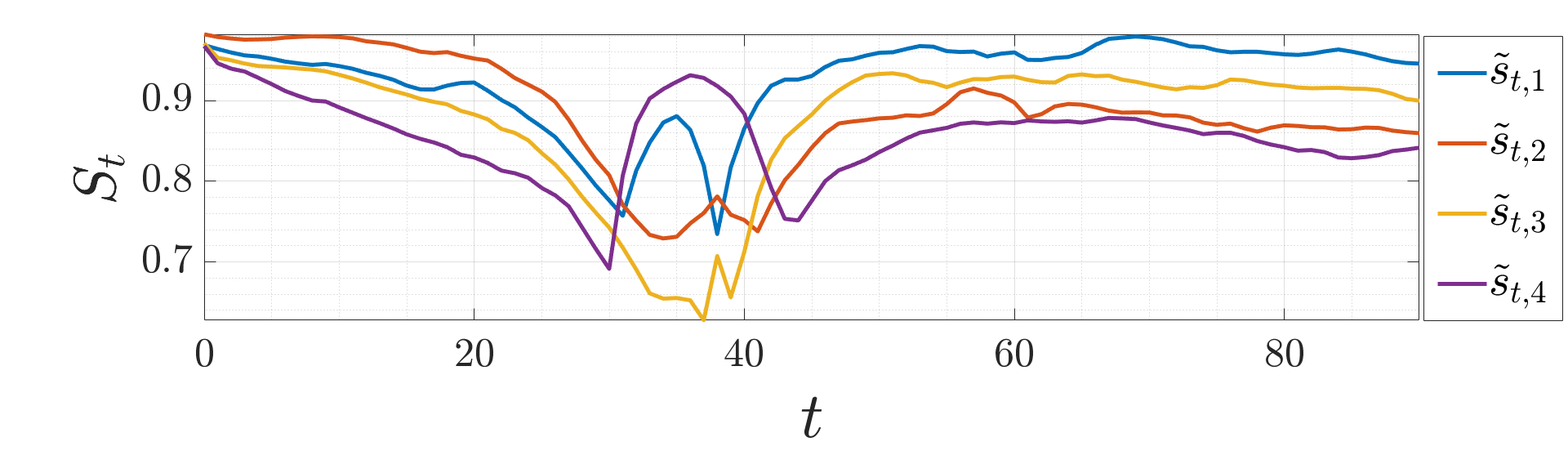}
\caption{\label{fig:case_Dip_FullTO_restrictedInit} Leading $4$ rolling windows of singular values for the dancing dipoles model for the full transfer operator restricted to the initial seeding of bins whose centres lie within a circular patch of radius $1/2$ centred at $(5.55,3.75)$. These are paired through time using Algorithm~\ref{alg:Seeding} for $n=10$ and $p=0.8$.}
\end{subfigure}
\hfill
\begin{subfigure}[c]{\columnwidth}
\hspace{0.8cm}
\includegraphics[width=0.8\textwidth]{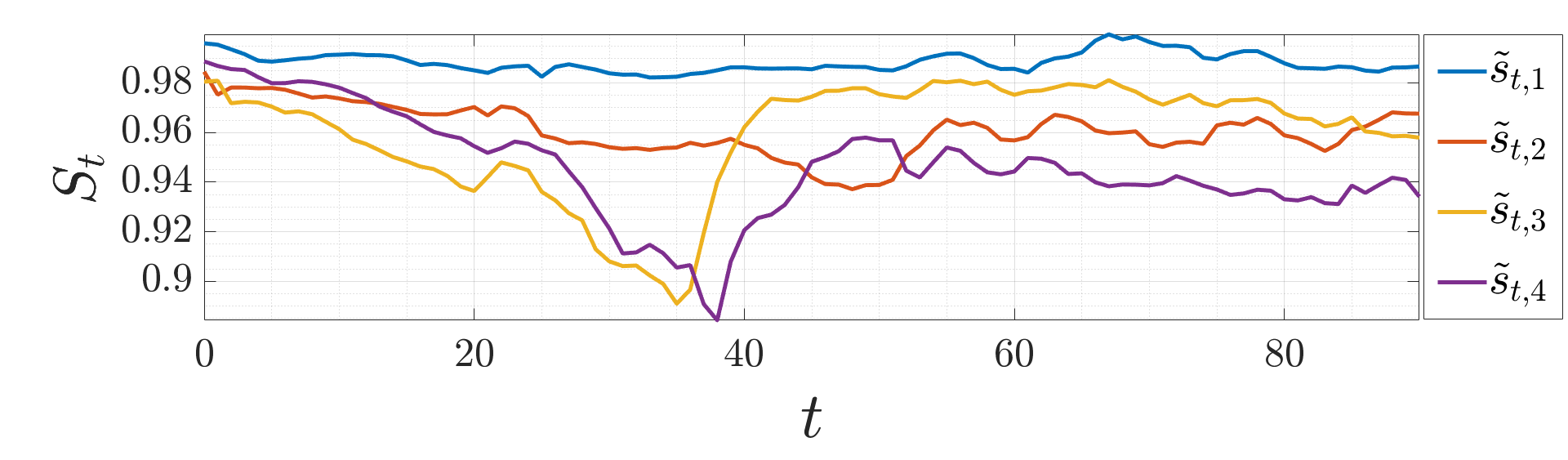}
\caption{\label{fig:case_Dip_FullTO_restrictedFin} Leading $4$ rolling windows of singular values for the dancing dipoles model for the full transfer operator restricted to the image of the initial seeding described in Figure~\ref{fig:case_Dip_FullTO_restrictedInit}. Pairing of instances through time is achieved using Algorithm~\ref{alg:Seeding} for $n=10$ and $p=0.5$.}
\end{subfigure}
\caption{\label{fig:case_Dip_FullTO_All} Rolling windows of singular value paths for the dancing dipoles model for relevant incarnations of the full transfer operator.}
\end{figure}

%%%% adddition of TO elaboration
{{Of interest at this point is an evaluation of how these results compare to those obtained using the {\textit{full}} transfer operator. A selection of results for the case of comparable, discrete time compositions of numerical approximations to the full transfer operator are presented in Figure~\ref{fig:case_Dip_FullTO_All}. Figure~\ref{fig:case_Dip_FullTO} presents the rolling window results for the full transfer operator defined in terms of the parameters utilised throughout this section. Figures~\ref{fig:case_Dip_FullTO_restrictedInit} and~\ref{fig:case_Dip_FullTO_restrictedFin} incorporate a restriction of the transfer operator to either the initial patch of radius $1/2$ or to the image of this patch, respectively. 
It is evident that Figure~\ref{fig:case_DipG_c}, which was obtained using a localised approach, and Figure~\ref{fig:case_Dip_FullTO_restrictedInit} produce similar results despite their distinct approaches.
}}

{{Indeed, our method moves beyond simply restricting the transfer operator to an initial or final region. This method restricts both the domain and target to build a truncated operator that does not require complete knowledge of the dynamics to detect coherently evolving structures. 
Whilst in general one may not expect such structures to coincide with global coherent structures of the flow, these results suggest our methods 
could be relevant in the analysis of phenomena where the localisation of some quantity is well understood at an initial time.
Furthermore, it is reasonable to expect that structures identified by our algorithms are meaningful for the global operator, although possibly not dominant, given the initially seeded region (or some sub-region of it) is not mixing significantly with its surrounds.
}}

%%%%%%%%%%%%%%%%%%%%%%%%%%%%%%%%%%%%%%%%%%%%%%%%%%%%%%%%%%%%%%%%%%%%%%%%%
\FloatBarrier \subsubsection{Merging monopoles}\label{SSec:MMs}
\begin{figure}[!htbp]
  \includegraphics[width=\textwidth]{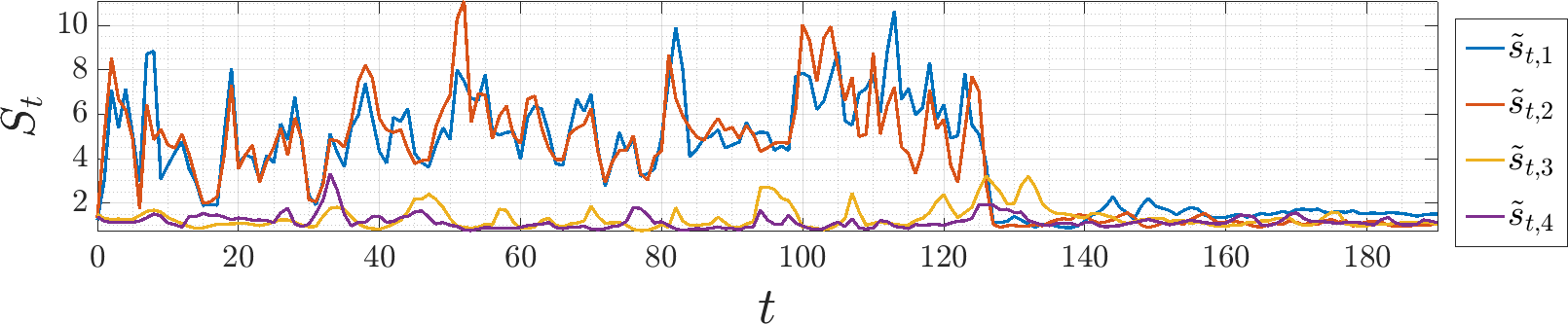}
  \caption{\label{fig:case_Monopoles_RWs} Rolling windows of singular values for the leading $4$ modes, tracked using Algorithm~\ref{alg:Seeding} with $p=0.1$ and $t$ is given in terms of inertial periods. The patched region is a circle of radius $1$ centred at $(\pi,\pi)$.}
\end{figure}

%%%%%%%%%%%%%%%%%%%%%%%%%%%%%%%%%%%%%%%%%%%%%%%%%%%%%%%%%%%%%%%%%%%%%%%%%
We now turn to the case of the merging monopoles, where we seek to determine if patching an area characterised by the collision of two poles provides details about the lifespan of coherent structures. This collision is illustrated in Figure~\ref{fig:Monopoles}. We choose to instantiate a circular patch of radius $1$ centred at $(\pi,\pi)$ for $p=0.1$ and the $95\%$ threshold described in Algorithm~\ref{alg:lifespan}. Implementing Algorithm~\ref{alg:Seeding} gives the tracked paths for rolling windows of singular values as shown in Figure~\ref{fig:case_Monopoles_RWs}.

Firstly, one notes the dramatic decline in the difference between the leading $4$ paths that occurs around $t=130$ inertial periods. There one notes that time windows associated with $t$ near to $130$ are associated with a merging of the two main poles. This can be seen in Figure~\ref{fig:Monopoles}, where one notes that by $t=133$ the two poles have begun to merge. The final plot on the (bottom) right in Figure~\ref{fig:Monopoles} shows that by $t=138$ inertial periods this merging is complete. 

Secondly, the fall in $\{ \tilde{s}_{t,1} \}$ and $\{ \tilde{s}_{t,2} \}$ that occurs alongside a peak in $\{ \tilde{s}_{t,4} \}$ near to $t=30$, is noteworthy. In terms of the dynamics, one notes the development of two weaker and much smaller structures occurs around this time. These structures are located on the outer parts of the larger poles and can be seen in panel (b) of Figure~\ref{fig:Monopoles}. These smaller structures begin to form around time $20$, becoming more stable in size between $30$ and $40$ inertial periods. They continue to exist until around $125$ inertial periods, just prior to the merging of the main poles.

\begin{figure}[!htbp]
  \includegraphics[width=\textwidth]{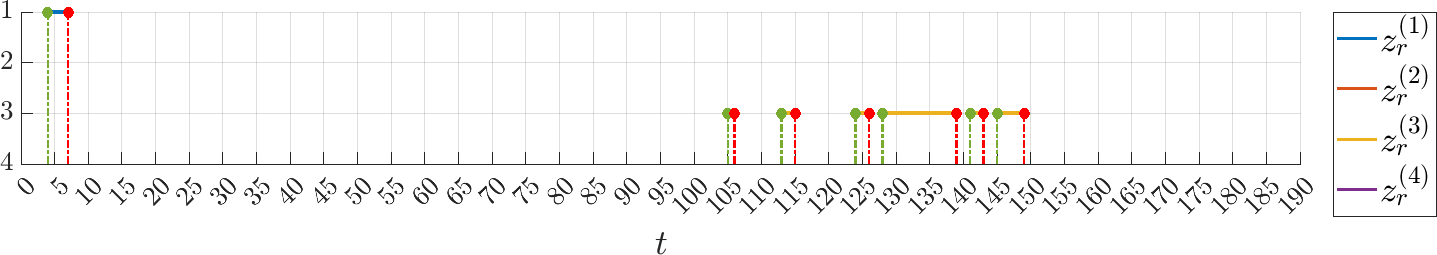}
  \caption{\label{fig:case_MONOPOLEs_lifespans_reg} Regular coherent structures (Algorithm~\ref{alg:CSorNot}) for pairings of length $1$ or more.}
\end{figure}

We again use Algorithm~\ref{alg:lifespan} alongside our three methods from Algorithm~\ref{alg:3_lives} to identify possible lifespans of interest. The longest regularised lifespan is then identified using
Algorithm~\ref{alg:CSorNot} with our standard isoperimetric threshold of $0.85$. {These results are shown in Figure~\ref{fig:case_MONOPOLEs_lifespans_reg}.} In this case, the longest regularised lifespan is a subset of $z_{Eldest}$ for $t \in [128,139]$. This lifespan, associated with  $z_{r}^{(3)}$, identifies a merging of the two main poles.

\begin{figure}[!htbp]
\setlength{\tabcolsep}{0.2pt}
\begin{tabularx}{\columnwidth}{|p{0.525cm}| *5{>{\Centering}X}|}
\cline{2-6} \multicolumn{6}{c}{\vspace{-0.5cm}}\\
\multicolumn{1}{c}{} &
\multicolumn{1}{|c}{\cellcolor{S3!60}} &
\multicolumn{1}{|c}{\cellcolor{S3!60}} &
\multicolumn{1}{|c}{\cellcolor{S3!60}} &
\multicolumn{1}{|c}{\cellcolor{S3!60}} &
\multicolumn{1}{|c|}{\cellcolor{S3!60}} \\[-1.em] 
\multicolumn{1}{c}{} &
\multicolumn{1}{|c}{\cellcolor{S3!60} $\tilde{u}_{128,3}$}&
\multicolumn{1}{|c}{\cellcolor{S3!60} $\tilde{u}_{129,3}$}&
\multicolumn{1}{|c}{\cellcolor{S3!60} $\tilde{u}_{130,3}$}&
\multicolumn{1}{|c}{\cellcolor{S3!60} $\tilde{u}_{133,3}$}&
\multicolumn{1}{|c|}{\cellcolor{S3!60} $\tilde{u}_{138,3}$}\\
\hline
\begin{minipage}[t][][b]{0.03\textwidth}
\centering
\cellcolor{S3!60}
\vspace{0.85cm}
\footnotesize{$\,z^{(3)}\;$}
\end{minipage}
&
\begin{minipage}[t][][b]{0.2\textwidth}
	\centering
	\includegraphics[height=2.725cm,center]{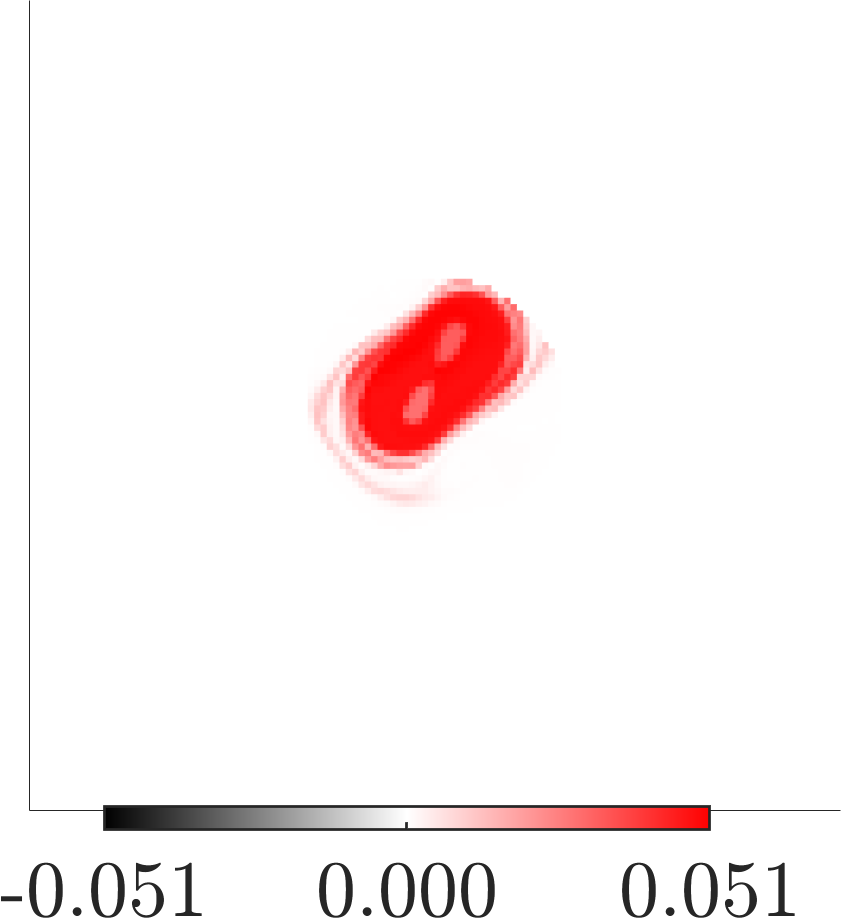} 
\vspace{-0.25cm}
\end{minipage}
&
\begin{minipage}[t][][b]{0.2\textwidth}
	\centering
	\includegraphics[height=2.725cm,center]{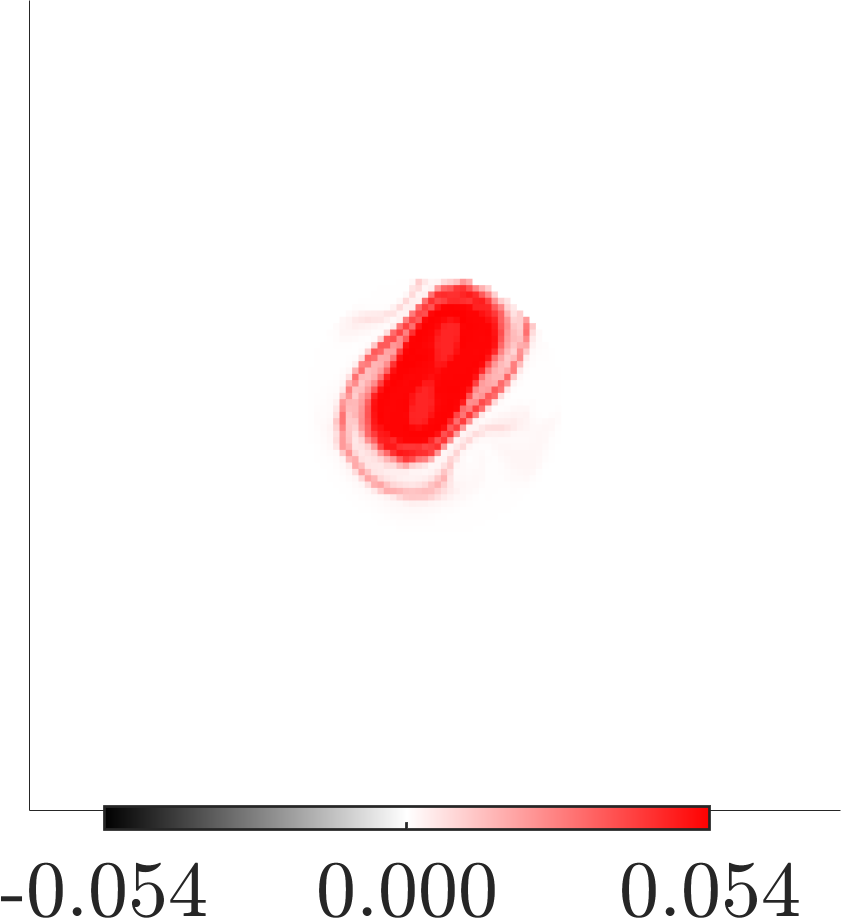} 
\vspace{-0.25cm}
\end{minipage}
&
\begin{minipage}[t][][b]{0.2\textwidth}
	\centering
	\includegraphics[height=2.75cm,center]{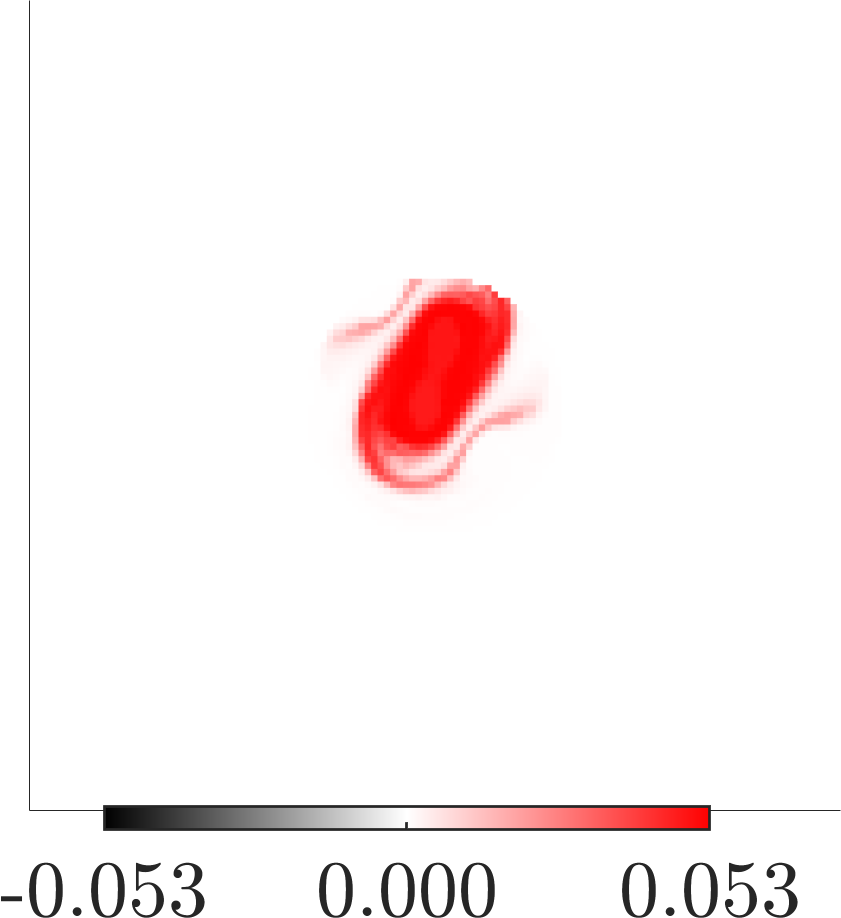} 
\vspace{-0.25cm}
\end{minipage}
&
\begin{minipage}[t][][b]{0.2\textwidth}
	\centering
	\includegraphics[height=2.725cm,center]{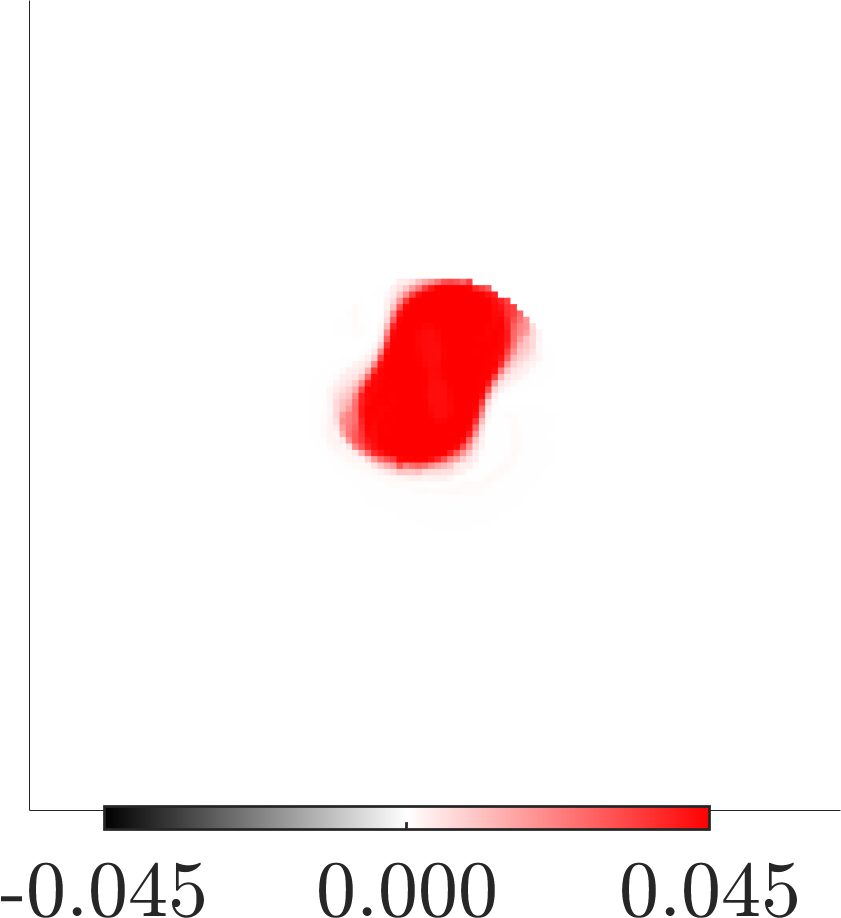} 
\vspace{-0.25cm}
\end{minipage}
&
\begin{minipage}[t][][b]{0.2\textwidth}
	\centering
	\includegraphics[height=2.725cm,center]{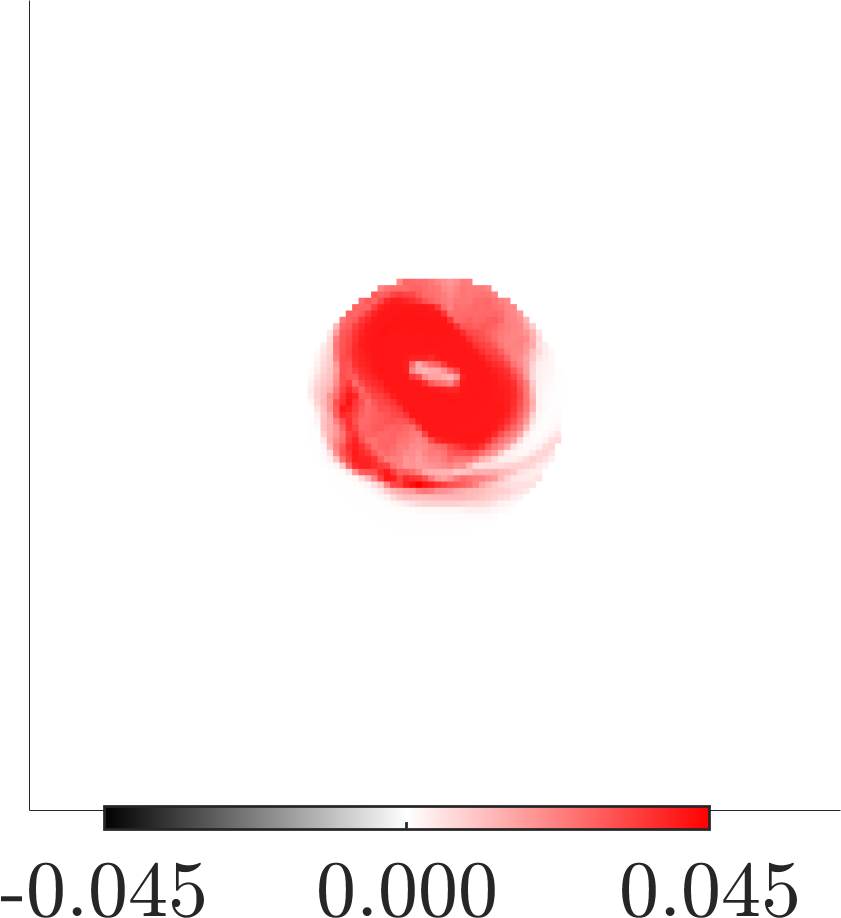} 
\vspace{-0.25cm}
\end{minipage}
\tabularnewline 
\hline
\end{tabularx}
\caption{Selected illustrations of the behaviour of the lifespan $z_{Eldest}$ in the patched region at various initial times for a variety of time windows. The coherent structures identified by these dynamics corresponds to the merging of the poles illustrated in Figure~\ref{fig:Monopoles}.}
\label{fig:maxVar_table_caseMM_n10_z3}
\end{figure}

Moreover, Figure~\ref{fig:maxVar_table_caseMM_n10_z3} presents a detailed exploration of the initial-time singular vectors $\tilde{u}$ associated with $z_{Eldest}$ and the merging event of interest. As expected from the illustrations in Figure~\ref{fig:Monopoles}, this merging is shown to occur following $130$ inertial periods. This also corresponds to the dramatic fall in the variance of the leading paths of rolling windows of singular values shown in Figure~\ref{fig:case_Monopoles_RWs}. Furthermore, our previous examination of Figure~\ref{fig:Monopoles} suggested that the merger was expected to be complete by $138$ inertial periods. This corresponds well with the regularised lifespan $z_{r}^{(3)}$ that exists for time windows initialised when $t \in [128,139]$.

%\vspace{-0.5cm}
%%%%%%%%%%%%%%%%%%%%%%%%%%%%%%%%%%%%%%%%%%%%%%%%%%%%%%%%%%%%%%%%%%%%%%%%%
\FloatBarrier \subsubsection{Random initial conditions}\label{SSec:RICs_Results} %(3.75,4) a=1.5, b=3/4
\begin{figure}[H]% $(3.75,4)$, $a=1$, $b=3/4$ $n=10$ tracked
\includegraphics[width=\textwidth]{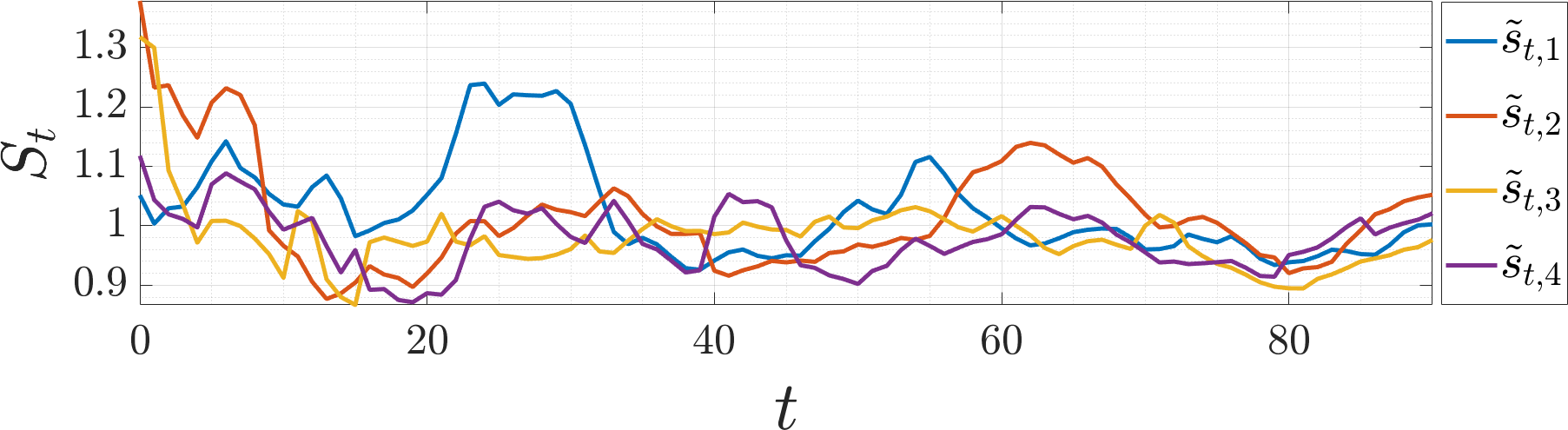}
  \caption{Rolling windows of singular values, tracked using $p=0.1$ in Algorithm~\ref{alg:Seeding} for depth$=14$, $n=10$ and an elliptical patch centred at $(3.75,4)$ with semi-major and semi-minor axes $a=3/2$ and $b=3/4$. Here $t$ is given in terms of inertial periods.}
  \label{fig:case_RandICsG_a}
\end{figure}

%%%%%%%%%%%%%%%%%%%%%%%%%%%%%%%%%%%%%%%%%%%%%%%%%%%%%%%%%%%%%%%%%%%%%%%%%
The final numerical test we consider is one generated using random initial conditions, as described in Section~\ref{SSec:RICs}. Here it is known that two of the more dominant structures will merge in an area above the centre of phase space that is slightly to the right. This is illustrated in Figure~\ref{fig:complex}. Using this information, we seed an elliptical patch centred at $(3.75,4)$ of semi-major and minor axes lengths $3/2$ and $3/4$ and utilise $p=0.1$ alongside the standard threshold of Algorithm~\ref{alg:lifespan}. Rolling windows of singular values for this case are presented in Figure~\ref{fig:case_RandICsG_a}. 

The clearly separated peak that emerges in $\{ \tilde{s}_{t,1} \}$ for $t$ between $20$ and $30$ inertial periods is of immediate interest. Whilst this peak is well separated, it does not grow to a value dramatically greater than one. This suggests a dynamically meaningful mode is likely associated with this path, rather than simply a point attractor. Two smaller peaks with similar characteristics also develop in $\{ \tilde{s}_{t,1} \}_{50 \le t \le 60}$ and $\{ \tilde{s}_{t,2} \}_{50 \le t \le 70}$.

We identify lifespans of interest using the three methods described in Algorithm~\ref{alg:3_lives}. In this case, the peaks between times $50$ and $70$ are not identified as being associated with any lifespan of immediate interest. On the other hand, the peak in {$\{ \tilde{s}_{t,1} \}$} that develops between times $20$ and $30$ is identified by $z_{MaxVarSV,2}$ which exists for $t \in [16,29]$. Further exploration of this mode is undertaken using Algorithm~\ref{alg:CSorNot}. 

\begin{figure}[H]
\setlength{\tabcolsep}{0.2pt}
\begin{tabularx}{\columnwidth}{|p{0.525cm}| *5{>{\Centering}X}|}
\cline{2-6} \multicolumn{6}{c}{\vspace{-0.5cm}}\\
\multicolumn{1}{c}{} &
\multicolumn{1}{|c}{\cellcolor{S1!60}} &
\multicolumn{1}{|c}{\cellcolor{S1!60}} &
\multicolumn{1}{|c}{\cellcolor{S1!60}} &
\multicolumn{1}{|c}{\cellcolor{S1!60}} &
\multicolumn{1}{|c|}{\cellcolor{S1!60}} \\[-1.em] 
\multicolumn{1}{c}{} &
\multicolumn{1}{|c}{\cellcolor{S1!60} $\tilde{u}_{19,1}$}&
\multicolumn{1}{|c}{\cellcolor{S1!60} $\tilde{u}_{21,1}$}&
\multicolumn{1}{|c}{\cellcolor{S1!60} $\tilde{u}_{22,1}$}&
\multicolumn{1}{|c}{\cellcolor{S1!60} $\tilde{v}_{21,1}$}&
\multicolumn{1}{|c|}{\cellcolor{S1!60} $\tilde{v}_{22,1}$}\\
\hline
%\vspace{-0.1cm}
\begin{minipage}[t][][b]{0.03\textwidth}
\centering
\cellcolor{S1!60}
\vspace{0.85cm}
\footnotesize{$\,z^{(1)}$}
\end{minipage}
&
\begin{minipage}[t][][b]{0.2\textwidth}
	\centering
	\includegraphics[height=2.725cm,center]{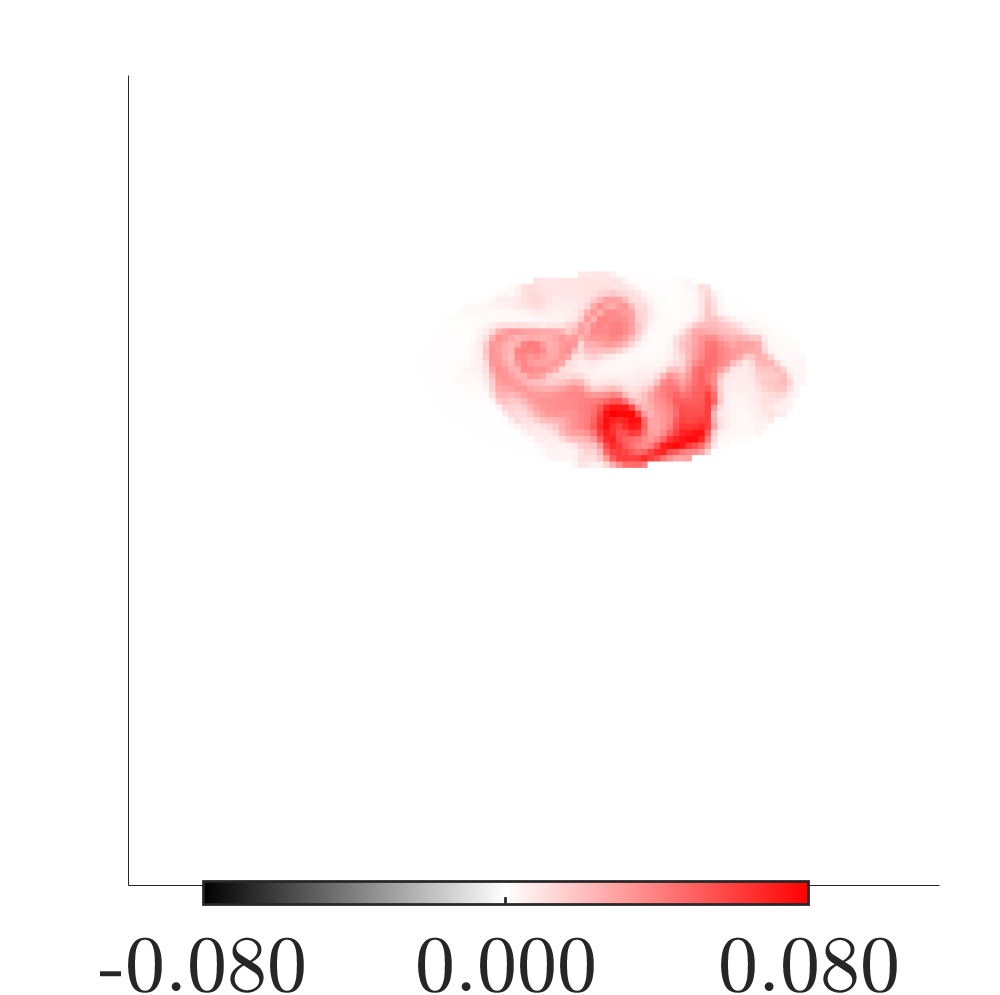} 
\vspace{-0.25cm}
\end{minipage}
&
\begin{minipage}[t][][b]{0.2\textwidth}
	\centering
	\includegraphics[height=2.725cm,center]{f1c.png} 
\vspace{-0.25cm}
\end{minipage}
&
\begin{minipage}[t][][b]{0.2\textwidth}
	\centering
	\includegraphics[height=2.725cm,center]{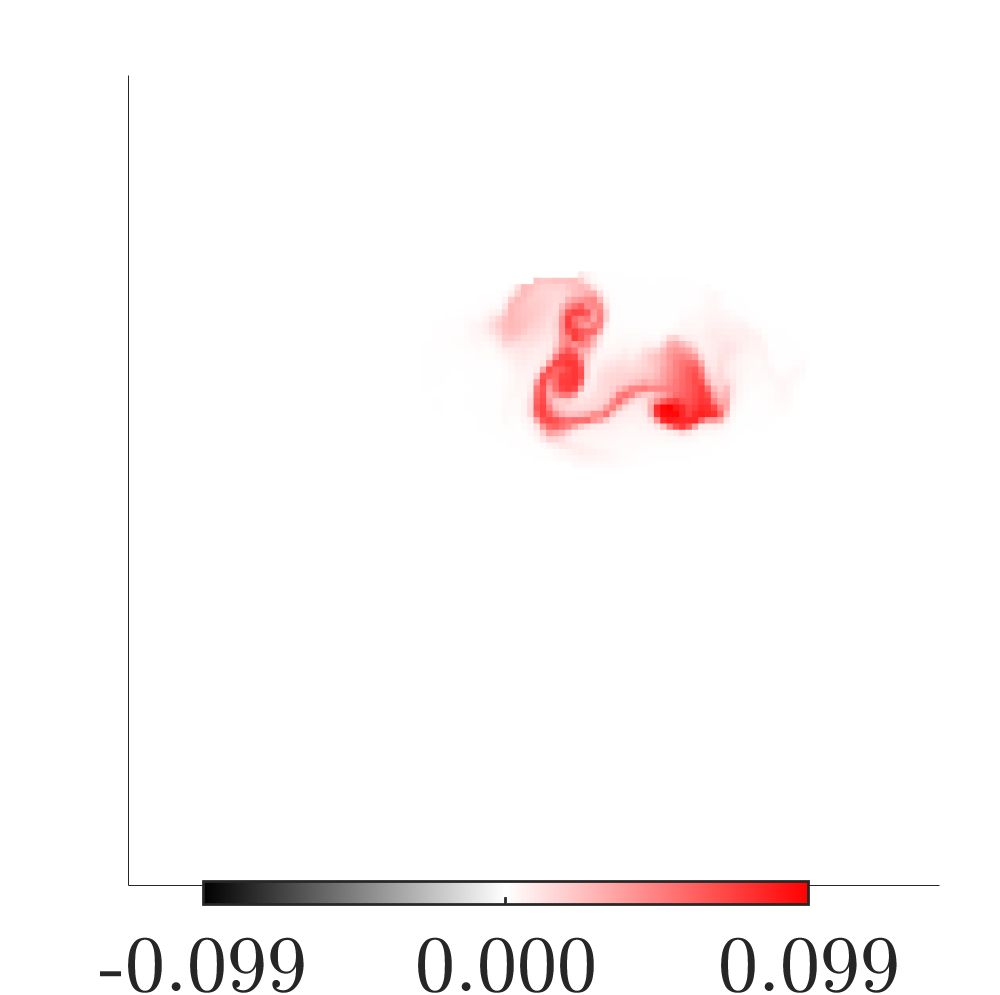} 
\vspace{-0.25cm}
\end{minipage}
&
\begin{minipage}[t][][b]{0.2\textwidth}
	\centering
	\includegraphics[height=2.725cm,center]{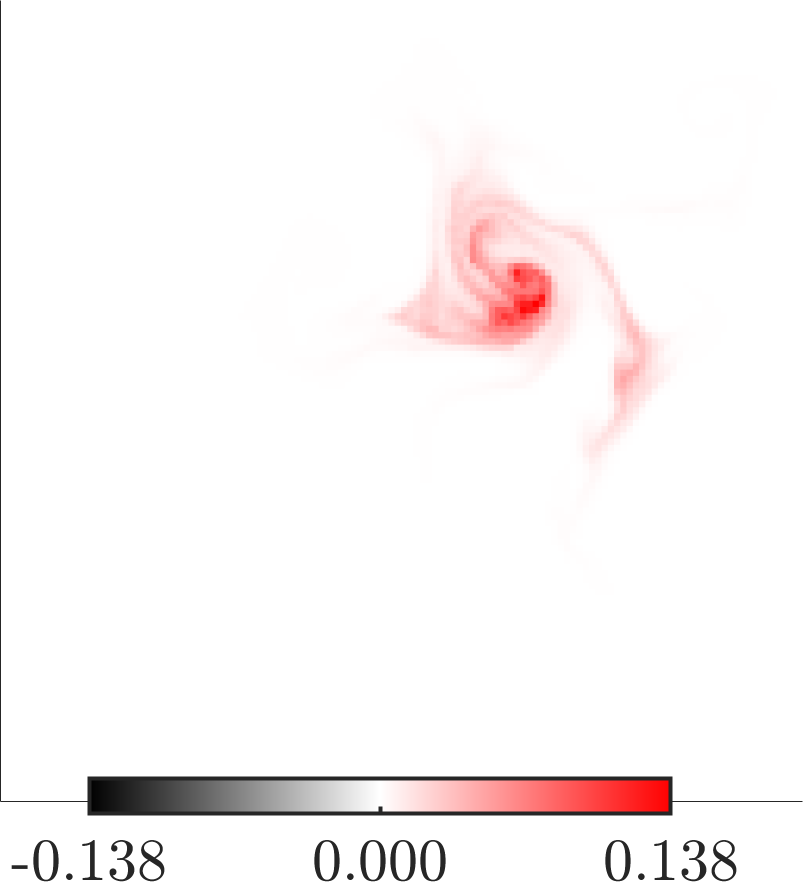} 
\vspace{-0.25cm}
\end{minipage}
&
\begin{minipage}[t][][b]{0.2\textwidth}
	\centering
	\includegraphics[height=2.725cm,center]{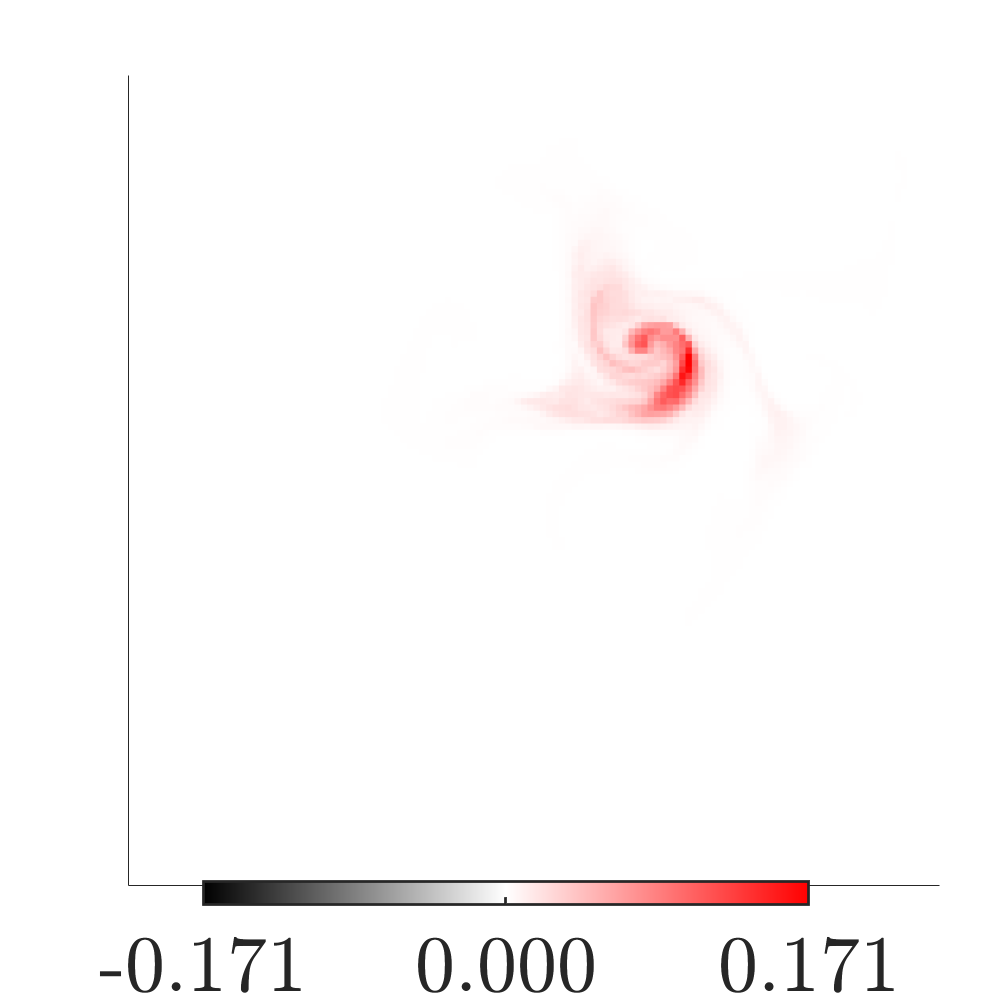}
\vspace{-0.25cm}
\end{minipage}
\tabularnewline 
\hline
\end{tabularx}
\caption{Selected illustrations of the behaviour of the lifespan $z_{MaxVarSV,2}$ in the patched region at various times for a variety of time windows. The coherent structures identified by these dynamics corresponds to the isolated peak in $\{\tilde{s}_{t,1}$\} illustrated in Figure~\ref{fig:case_RandICsG_a}.}
\label{fig:maxVar_table_RICG_a_n10_z1}
\end{figure}

Figure~\ref{fig:maxVar_table_RICG_a_n10_z1} presents a detailed view of various left singular vectors associated with time windows that comprise the regularised components $z_{r}^{(1)}$ of the lifespan $z_{MaxVarSV,2}$. This lifespan is clearly associated with two increasingly entwined structures that eventually merge to form a new structure by the time the lifespan ceases.

\begin{figure}[!htbp]
  \includegraphics[width=\textwidth]{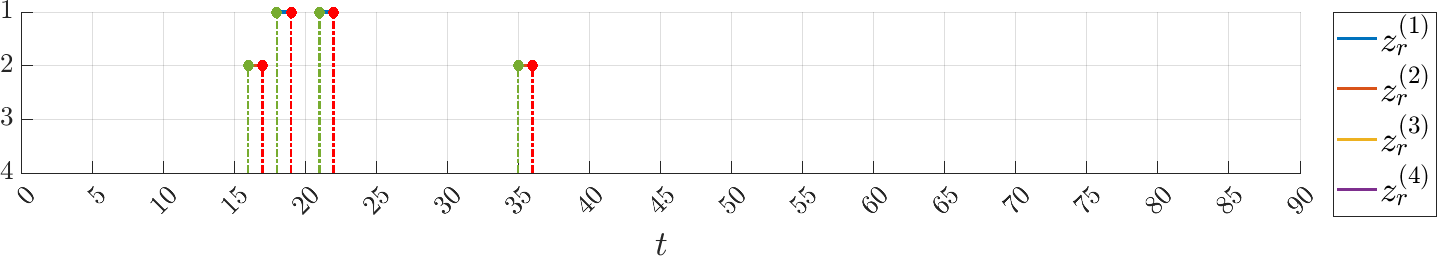}
  \caption{Time windows associated with regular coherent structures (Algorithm~\ref{alg:CSorNot}) for pairings of length unity or greater.}
  \label{fig:case_RandICsG_a_lifespans_reg}
\end{figure}

Applying the additional layer provided by Algorithm~\ref{alg:CSorNot} to the four lifespans revealed in Figure~\ref{fig:maxVar_table_RICG_a_n10_z1}, results in the {regularised lifespans} shown in Figure~\ref{fig:case_RandICsG_a_lifespans_reg}. Utilising our standard regularity threshold of $0.85$ results in only certain instances of $z_{MaxVarSV,2}$ being present in the regularised lifespans.

Furthermore, the regularised lifespan, $z_{r}^{(1)}$ is split into {two} smaller components. One notes from Figure~\ref{fig:maxVar_table_RICG_a_n10_z1}, that time $\tilde{u}_{19,1}$ corresponds to a subtle change, which includes a greater emphasise on the tail of the merging structures. On the other hand, $\tilde{u}_{21,1}$  and $\tilde{u}_{22,1}$ are associated with a more balanced emphasis on the structures involved in the merging event. The outcome of the merger is illustrated by {the associated right singular vectors} in Figure~\ref{fig:maxVar_table_RICG_a_n10_z1}. As such, this division of the dominant lifespans into much smaller components better clarifies the exact dynamics of interactions that characterise particular lifespans. 

%%%%%%%%%%%%%%%%%%%%%%%%%%%%%%%%%%%%%%%%%%%%%%%%%%%%%%%%%%%%%%%%%%%%%%%%%
\FloatBarrier \section{Concluding remarks}
\label{Sec:Conclusions}
%%%%%%%%%%%%%%%%%%%%%%%%%%%%%%%%%%%%%%%%%%%%%%%%%%%%%%%%%%%%%%%%%%%%%%%%%
The algorithms developed in this work have been useful in isolating dynamically meaningful objects and their associated lifespans. Such objects included not only coherently evolving structures, but also coherent structures that interact and experience fundamental structural changes, such as merging and separation events. Furthermore, the optional, additional layers of complexity provided by Algorithms~\ref{alg:3_lives} and~\ref{alg:CSorNot} were especially useful in the identification of particular lifespans, such as those associated with the longest lived or the most dynamically variant structures.

The key tool in this study was a numerical approximation to the localised (non-global) transfer operator. Utilising this methodology allowed us to isolate dynamically meaningful information, such as the presence of coherent structures in localised regions. This strategy, of limiting one's analysis to localised regions of phase space, has important ramifications. Firstly, it allows for a rapid analysis when the localised target area can be approximated. Secondly, it allows for more concise and focused results that clearly determine the presence of coherent structures or dynamic events. Moreover, it allows for effective results to be garnered without the complications introduced by large amounts of background noise, such as multiple dynamic events that involve the merging or separation of vortices across various regions.

Rolling windows of singular values associated with compositions of conditional (localised) matrices also exhibited important signals regarding the dynamical behaviour of coherent structures. For example, in the patched regions of our analyses, a well defined separation of rolling windows of singular values was often critical to the identification of the presence of coherent structures. More importantly, rolling windows of singular values for the conditional Ulam matrices could always be utilised to identify the presence of coherent structures in a patched region. Distinct changes in singular values signalled either the entry of a coherent structure into a patched region, or the exit of one out of a patched region. When a merger occurred, similar signals were observed as when a coherent structure was exiting the patched region. That is, relative to its previous trajectory, the associated path was observed to fall rapidly in range and dominance over time. 

Not only did our algorithms detect the presence of coherent structures in patched regions, they also identified periods when no coherent structures were present. Likewise, these algorithms allowed us to differentiate between spurious structures and dynamically meaningful information. Whilst the exact efficacy of our methods requires further investigation, these methods were shown useful in the analysis of various dynamical settings. For example, our methods identified the presence of coherent structures in the simple setting provided by a periodically shifting double well potential and fundamental structural changes in the more complex systems modelled using the Boussinesq equations. 

Our methods utilised the equivariance mismatch between pairs of vectors to develop those algorithms aimed at the detection of lifespans. These lifespans were refined using additional algorithmic layers. In particular, structures that were associated with lifespans exhibiting the largest variance of singular values were found to be particularly informative, as were those lifespans characterised by an associated singular vector that contained at least one coherent structure of a certain isoperimetric regularity. Algorithm~\ref{alg:CSorNot} was found to be especially effective in isolating subperiods of dynamical interest as well as separating particular behaviour from other anomalies in the dynamics. Importantly, despite their alternate approaches, both algorithms were found to isolate similar objects.

Finally, whilst singular values signalled the presence of coherent structures in patched regions by exhibiting behaviour that was distinct from previous patterns in terms of separation and variance, further research into the exact meaning of particular signals is necessary. For example, our results show that paths of rolling windows of singular values move in tandem before separating, or the converse, around merging and separation events. However, delving deeper into the behaviour associated with particular signals could be bolstered by the development of new insights into the choice of the most appropriate time window length $n$. Furthermore, one must not forget that the consequences of strong vortex interactions extend well beyond direct mergers or separations. For example, when a smaller structure approaches a stronger one there is also the possibility that this smaller structure will be stretched into a vortex sheet that then wraps itself around the larger structure~\cite{MEUNIER2005431}. This type of complicated behaviour, as well as the possible incorporation of our approach with other algorithms such as sparse eigenbasis approximation~\cite{FroylandRossSakeralliou}, could also be considered in future work that seeks to characterise the lifespan of coherent structures and how this is impacted by dynamical events.

\section*{Acknowledgments}
The authors thank P. Koltai, K. Padberg-Gehle for their input on an earlier version of this manuscript.
C. Blachut and C. Gonz\'alez-Tokman have been partially supported by the Australian Research Council and by the University of Queensland's Promoting Women Fellowship program. C. Blachut has been
supported by an Australian Government Research Training Program Stipend Scholarship at The University of Queensland and partially supported by funding from the Australian Research Council (Grant DP200101764).
G. Hernandez-Duenas was supported in part by grants UNAM-DGAPA-PAPIIT IN113019 \& Conacyt A1-S-17634. Some simulations were performed at the Laboratorio Nacional de Visualización Científica Avanzada at UNAM Campus Juriquilla, and G. Hernandez-Duenas received technical support from Luis Aguilar, Alejandro De León, and Jair García from that lab.
{
\section*{Data availability}
The datasets generated during and/or analysed during the current study are available from the corresponding author on reasonable request.
}

{
\section*{Conflict of Interest}
The authors declare that they have no conflict of interest.
}
%%%%%%%%%%%%%%%%%%%%%%%%%%%%%%%%%%%%%%%%%%%%%%%%%%%%%%%%%%%%%%%%%%%%%%%%%

% ***************************************************
% Bibliography
%****************************************************
\normalem
\bibliographystyle{abbrv}
\bibliography{./PITS9_bib.bib}

\appendix
\section{Algorithms}\label{appd:algs}
%%%%%%%%%%%%%%%%%%% Alg 1
\begin{algorithm}[H] 
\caption{Building non-global, conditional Ulam matrices}
\label{alg:Seeding}
\begin{algorithmic}[1]
\Require numbers $t_{i}$, $t_{F}$, $n$, $\mathcal N$, $Q$ and $\vars{depth}$, vector field $\vect{u}$ and region $\vars{patch}$
\State Initialise $m \gets 2^{\vars{depth}}$ and $I,J \gets \emptyset$
\For{$t \gets t_{i}$ to $t_{F}-n$} 
\State {$\mathcal B$}  $=\{ B_{1}, B_{2}, \ldots, B_{m} \}$ with associated \textit{centres} $\{ c_{1}, c_{2}, \ldots, c_{m} \}$\label{alg:Seeding:partition}

\For {$i \gets 1$ to $m$} 
\If {$c_{i}$ {lies inside} $\vars{patch}$} 
\State Uniformly distribute $Q$ test points $x_{i,1}, x_{i,2}, \ldots, x_{i,Q}\in B_{i}$ \label{alg:Seeding:unifDistn}

\For {$q \gets 1$ to $Q$} 
\State $J \gets J \cup \{j\}$ with $T_{t,1} (x_{i,q}) \in B_{j}$ obtained via integration of $\vect{u}$ \label{alg:Seeding:flow} 
\EndFor

\State $I \gets I \cup \{i\}$ \label{alg:Seeding:I1}

\EndIf
\EndFor

\State $(P{(t,1)})_{i,j} = \frac{1}{Q} \sum\limits_{q=1}^{Q} \mathbbm{1}_{B_{j}}(T_{t,1} (x_{i,q}))$ where $i \in I$ and $j \in J$
\State $I \gets J$,  $J \gets \emptyset$ \label{alg:Seeding:part2}
\For{$\tilde{n} \gets 2$ to $n$} 
\For {$i \in I$} 
\State Uniformly distribute $x_{i,1}, x_{i,2}, \ldots, x_{i,Q}\in B_{i}$ \label{alg:Seeding:unifDistn2}
\For {$q \gets 1$ to $Q$} 
\State $J \gets J \cup \{j\}$ with $T_{t,\tilde{n}} (x_{i,q}) \in B_{j}$ obtained via integration of $\vect{u}$ \label{alg:Seeding:endPts2}
\EndFor
\EndFor

\State $(P{(t,\tilde{n})})_{i,j} = \frac{1}{Q} \sum\limits_{q=1}^{Q} \mathbbm{1}_{B_{j}}(T_{t,\tilde{n}} (x_{i,q}))$ where $i \in I$ and $j \in J$
\State $I \gets J$,  $J \gets \emptyset$

\EndFor 

\State $P^{(n)}_{t} \gets P({t,1}) \cdot P({t,2}) \cdots P({t,n})$ \label{alg:Seeding:cond_ccycle}
\State $[U_{t}^{(n)}, \, S_{t}^{(n)}, \,V_{t}^{(n)}] \gets {svds}(P^{(n)}_{t}, \; \mathcal N)$\label{alg:Seeding:SVD}

\EndFor 
\Ensure $\{ P_{t}^{(n)} \}_{ t_{i} \le t \le t_{F}-n}$, $\{ P({t,\tilde{n}}) \}_{ t_{i} \le t \le t_{F}-n, 1 \le \tilde{n} \le n}$, 
\newline $\{ S_{t}^{(n)} \}_{ t_{i} \le t \le t_{F}-n}$,  $\{ U_{t}^{(n)} \}_{ t_{i} \le t \le t_{F}-n}$, $\{ V_{t}^{(n)} \}_{ t_{i} \le t \le t_{F}-n}$
\end{algorithmic}
\end{algorithm}
%%%%%%%%%%%%%%%%%%%

%%%%%%%%%%%%%%%%%%% Alg 2
\begin{algorithm}[H] 
\caption{Tracking modes through time using singular vectors}
\label{alg:track_norm}
\begin{algorithmic}[1]
\Require $t_{i}$, $t_{F}$, $n$, $\mathcal N$, $p$ and collections $\{ S_{t}^{(n)} \}$, $\{ V_{t}^{(n)} \}$, $\{ P({t,\tilde{n}}) \}$ from Alg.~\ref{alg:Seeding}
\State Define the initial mode association $\hat{S}_{t_{i},j} \gets S_{t_{i},j}^{(n)}$ for $j\in\{1,\ldots,\mathcal{N}\}$
\For {$t \gets t_{i}$ to $t_{F}-n-1$}
\State Define initial sets characterising all possible transitions $j',j'' \gets \{1, \ldots, \mathcal{N} \}$
\While{$j' != \emptyset$ } 
\State \vspace{-0.45cm}
\begin{quote}
{\em{Lift}} $v_{t,j'}^{(n)}$ and $P({t+1,n})v_{t+1,j''}^{(n)}$ to a common dimension,
determined by the union of their supports,
by augmenting with $0$s.\label{alg:state:lift}
\end{quote}
\State \hspace{-0.45cm} \vars{dist} $ \gets \min_{j',j''}{ \left( \; {\left\| {v_{t,j'}^{(n)}} \pm P({t+1,n})v_{t+1,j''}^{(n)}/ {\| P({t+1,n})v_{t+1,j''}^{(n)} \|}_{2} \; \right\| }_{p}\; \right)} $\label{alg:track_norm:min}
\State \hspace{-0.45cm} \vars{modes} $ \gets \arg\min_{j',j''}\left(\vars{dist}\right)$
\State \hspace{-0.45cm} Set 
$j' \gets j' \setminus \vars{modes}(1)$ and 
$j'' \gets j'' \setminus \vars{modes}(2)$% from j' to j''
\State \hspace{-0.45cm} Create the new mode association $\hat{S}_{t+1,\vars{modes}(1)} \gets S_{t+1,\vars{modes}(2)}^{(n)}$
\EndWhile
\EndFor
\State Characterise each of the $j$ paths by average value over all $t$ \par $\bar{S}_{j} \gets \frac{1}{t_{F}-n-1-t_{i}} \sum_{t=t_{i}}^{t_{F}-n-1} \hat{S}_{t,j}$
\State Re-sort paths in (descending) order of average value \par $[ \, \sim \, , \; \{\vars{{sorted\_modes}}\}] \gets \textrm{sort}(\bar{S}, \; \textit{`descending'})$
\State The singular value paths of ordered modes are determined by $\{\tilde{S}_{t}^{(n)} \}$ 
where \par
$\{\tilde{s}_{t,j}^{(n)} \}
\gets
\{\hat{S}_{t,{\vars{sorted\_modes}(j)}} \}$
\State The ordered collection of $\mathcal{N}$ left and right singular vectors associated with these paths are $\{ \tilde{U}_{t}^{(n)} \}$ and $\{ \tilde{V}_{t}^{(n)} \}$ 
\Ensure Paired singular value paths $\{\tilde{s}_{t,j}^{(n)} \}_{ t_{i} \le t \le t_{F}-n, 1 \le j \le \mathcal{N}}$, \newline
\indent  the associated left singular vector paths $\{\tilde{u}_{t,j}^{(n)} \}_{ t_{i} \le t \le t_{F}-n, 1 \le j \le \mathcal{N}}$ \newline 
\indent and the paired right singular vectors $\{\tilde{v}_{t,j}^{(n)} \}_{ t_{i} \le t \le t_{F}-n, 1 \le j \le \mathcal{N}}$.
\end{algorithmic}
\end{algorithm}
%%%%%%%%%%%%%%%%%%%

%%%%%%%%%%%%%%%%%%% Alg 3
%\setcounter{algorithm}{2} % force to ignore first 2
\begin{algorithm}[H]
\caption{Equivariance and pairing mismatch and the determination of lifespans}
\label{alg:lifespan}
\begin{algorithmic}[1]
\Require $n$, $ t_{i}$, $t_{F}$, $\mathcal{N}$, \texttt{threshold_c}, \texttt{threshold_up}, \texttt{threshold_down}, \texttt{threshold_p}, $\{ P({t,\tilde{n}}) \}$ defined as per Alg.~\ref{alg:Seeding} and $\{ \tilde{v}_{t,j}^{(n)} \}$ from Alg.~\ref{alg:track_norm} 
\For{$j \gets 1$ to $\mathcal{N}$}
\State $ \varsigma_{t}^{(j)} \gets \min \left(
{ \left\| {\tilde{v}_{t,j}^{(n)}} \pm P({t+1,n})\tilde{v}_{t+1,j}^{(n)}/ {\| P({t+1,n})\tilde{v}_{t+1,j}^{(n)} \|}_{2} \; \right\| }_{2} \right) /{\sqrt{2}} $ \label{alg:lifespan:pm}
\newline \indent where $\tilde{v}_{t,j'}^{(n)}$ and $P({t+1,n})\tilde{v}_{t+1,j''}^{(n)}$ are augmented as per Alg.~\ref{alg:track_norm}, Op.~\ref{alg:state:lift}.
\For{$t \gets t_{i}$ to $t_{F}-n-2$}
\State $ \varsigma_{t+1}^{(j)} \gets \min \left(
{ \left\| {\tilde{v}_{t+1,j}^{(n)}} \pm P({t+2,n})\tilde{v}_{t+2,j}^{(n)}/ {\| P({t+2,n})\tilde{v}_{t+2,j}^{(n)} \|}_{2} \; \right\| }_{2} \right) /{\sqrt{2}} $
\newline \indent \indent 
where $\tilde{v}_{t+1,j}^{(n)}$ and $P({t+2,n})\tilde{v}_{t+2,j}^{(n)}$
are augmented as per Alg.~\ref{alg:track_norm}, Op.~\ref{alg:state:lift}.
\If{choosing conservative threshold}
\State {lifespan} $z_{j,t}$ exists unless $\varsigma_{t}^{(j)} \ge $ \texttt{threshold_c} 
\ElsIf{choosing percentage change of $\varsigma^{(j)}$ in either time direction}
\State {lifespan}  $z_{j,t}$ exists unless
$\left( \varsigma_{t}^{(j)}>\texttt{threshold_up}\right)$
\textbf{ or }
\newline \indent \indent 
$\left( \varsigma_{t}^{(j)}>\texttt{threshold_down}\textbf{ \& }\frac{ |\varsigma_{t}^{(j)} - \;\varsigma_{t+1}^{(j)}|}{ \min \left(\varsigma_{t}^{(j)}, \; \varsigma_{t+1}^{(j)} \right)}>\texttt{threshold_p} \right)$\label{alg:lifespan:IfEnd}
\EndIf
\EndFor
\EndFor
\Ensure Lifespans $\{z_{j,t}\}$ where $j \in \{1, \ldots, \mathcal{N} \}$ and $t \in \{t_{i}, t_{i}+1, \ldots, t_{F}-n-2 \}$ \newline
\indent and equivariance mismatch $\{  \varsigma_{t}^{(j)}  \}_{t_{i} \le t \le t_{F}-n-1, 1 \le j \le \mathcal{N}}$
%\EndFunction
\end{algorithmic}
\end{algorithm}
%%%%%%%%%%%%%%%%%%%

%%%%%%%%%%%%%%%%%%% Alg 4
\begin{algorithm}[H]
\caption{Identifying types of lifespans of dynamical interest}
\label{alg:3_lives}
\begin{algorithmic}[1]
\Require $n$, $\mathcal{N}$, $\{\tilde{s}_{t,j}^{(n)} \}$ from Alg.~\ref{alg:track_norm}, $\{z_{j,t}\}$ and $\{  \varsigma_{t}^{(j)}  \}$ from Alg.~\ref{alg:lifespan} 
\State $z_{Eldest}$, $z_{MinEq}$, $z_{MaxVarSV} \gets \emptyset$, \texttt{age}, \texttt{SV\_var}$\gets 0$ and \texttt{ME}$\gets \infty$
\For{$j \gets 1$ to $\mathcal{N}$}
\State $z^{(j)}$ contains all lifespans associated with $j$-th component of $\{z_{j,t}\}$
\For{ each individual lifespan in $z^{(j)}$}
\State \texttt{age_z}$=\#\{z_{\alpha},\ldots,z_{\omega}\}$
\State \texttt{ME_z}$=\frac{1}{\texttt{age_z}}\sum_{t=z_{\alpha}}^{z_{\omega}} \varsigma_{j,t}$
\State \texttt{SV\_var_z}$=\frac{1}{\texttt{age_z}-1}\sum_{t=z_{\alpha}}^{z_{\omega}} | \tilde{s}_{t,j}^{(n)} - \overline{s}|^{2} $ 
where $\overline{s}=\frac{1}{\texttt{age_z}}\sum_{t=z_{\alpha}}^{z_{\omega}} \tilde{s}_{t,j}^{(n)} $
\State $\texttt{age}=\max \left( \texttt{age}, \texttt{age_z}\right)$ \label{alg:line:age} for maximising mode $j'$, birth $z'_{\alpha}$, death $z'_{\omega}$\label{alg:line:EqMM}
\State $\texttt{ME}=
\min \left( \texttt{ME}, \texttt{ME_z}\right)$
with minimising $j''$, $z''_{\alpha}$, $z''_{\omega}$ \label{alg:line:SValVar}
\State $\texttt{SV\_var}=
\max \left( \texttt{SV\_var},\texttt{SV\_var_z}\right)$
with maximising $j'''$, $z'''_{\alpha}$, $z'''_{\omega}$
\EndFor
\EndFor
\Ensure Characteristic lifespans $z_{Eldest}=\{z_{j'_{},z'_{\alpha}},\ldots,z_{j',z'_{\omega}}\}$,
\newline \indent $z_{MinEq}=\{z_{j''_{},z''_{\alpha}},\ldots,z_{j'',z''_{\omega}}\}$ and $z_{MaxVarSV}=\{z_{j'''_{},z'''_{\alpha}},\ldots,z_{j''',z'''_{\omega}}\}$.
\end{algorithmic}
\end{algorithm}
%%%%%%%%%%%%%%%%%%%

%%%%%%%%%%%%%%%%%%% Alg 5
%\setcounter{algorithm}{4} % force to ignore first 
\begin{algorithm}[!htbp]%RWmoviesTop4Vec_RECTMATRIX_DWP_LFSpan_Regular.m
\caption{Additional layer for the detection of regular coherent structures}
\label{alg:CSorNot}
\begin{algorithmic}[1]
\Require $n$, $ t_{i}$, $t_{F}$, $\mathcal{N}$, $0 \le \texttt{iso_thresh} \le 1$, $\{ \tilde{v}_{t,j}^{(n)} \}$ of Alg.~2 and $\{z_{j,t}\}$ of Alg.~\ref{alg:lifespan}
\For{$j \gets 1$ to $\mathcal{N}$}
\For{$t \gets t_{i}$ to $t_{F}-n-2$}
\State \texttt{col} $\gets$ linear interpolation of $\tilde{v}_{t,j}^{(n)}$ to estimate values at bin corners 
\State threshold \texttt{col} into level sets with ascending contours $C = \{ C_{1}$, $C_{2}$, $C_{3} \}$
\If{ $\texttt{sign}\left( C_{1} \right)~!=~\texttt{sign} \left( C_{3} \right)$}
\State Collect elements of \texttt{col} either side of $C_{1}$ and $C_{3}$ 
\ElsIf{$\texttt{sign} \left( C_{3} \right) > 0$}
\State Collect elements of $\texttt{col}>C_{2}$
\Else
\State Collect elements of $\texttt{col}<C_{2}$
\EndIf
\State For each connected component in new collection, calculate \newline
\indent \indent area $A$, perimeter length $L$ and isoperimetric ratio, $\mathcal{I}=4A\pi/L^{2}$, 
\begin{flushright} {\Comment{Default the value for $1$ pixel to $\mathcal{I}=1$}} \end{flushright}  
\State Let $\mathcal{I}_{max}$ be largest isoperimetric ratio of connected components
\If{$\mathcal{I}_{max}>\texttt{iso_thresh}$}
\State regularised lifespan $z_{r}$ exists for this value of $t$ and $j$
\Else
\State regularised lifespan $z_{r}$ does not exist for this $t$ and $j$
\EndIf
\EndFor
\EndFor
\Ensure Regularised lifespans $\{(z_{r})_{j,t}\}$
\end{algorithmic}
\end{algorithm}
%%%%%%%%%%%%%%%%%%%

\end{document}